\documentclass[reqno,11pt]{amsart}
\usepackage{amsfonts}
\usepackage{graphicx}
\usepackage{amsfonts,amsmath, amssymb}
\usepackage{bbm,mathrsfs,mathscinet}
\usepackage{cancel}
\usepackage{hyperref}
\usepackage{marginnote}
\usepackage{color}
\allowdisplaybreaks[4]

\numberwithin{equation}{section}

\newcommand{\R}{\mathbb{R}}

\newtheorem{theorem}{Theorem}[section]

\newtheorem{lemma}[theorem]{Lemma}
\newtheorem{proposition}[theorem]{Proposition}
\newtheorem{remark}[theorem]{Remark}

\def\v{\varepsilon}

\def\t{\theta}

\def\s{\sigma}

\def\i{\infty}
\def\f{\frac}
\def\pa{\partial}

\def\dis{\displaystyle}
\def\la{\langle}
\def\ra{\rangle}
\def\MX{\mathbf{X}}
\def\MY{\mathbf{Y}}
\def\MBX{\mathbb{X}}
\def\MBY{\mathbb{Y}}
\def\bfu{\mathbf{u}}
\def\bfv{\mathbf{v}}
\def\bfq{\mathbf{q}}

\begin{document}
	
\title[Steady Prandtl Layer for Compressible Fluids]{Steady Prandtl Layer for Compressible Fluids}

\author[Y. Guo]{Yan Guo}
\address[Y. Guo]{The Division of Applied Mathematics, Brown University, Providence, USA.}
\email{yan\_guo@brown.edu}

\author[Y. Wang]{Yong Wang}
\address[Y. Wang]{Academy of Mathematics and Systems Science, Chinese Academy of Sciences, Beijing 100190, China }
\email{yongwang@amss.ac.cn}

\begin{abstract}
Despite its importance, there have been few PDE results to
investigate Prandtl layers for compressible fluids, in which the thermal
boundary layer for the temperature field interacts with the classical velocity
Prandtl boundary layer in a strong nonlinear fashion. We establish local-in-$x$
well-posedness of Prandtl layers for steady compressible fluids and construct
Prandtl layer expansion in strong Sobolev norms, a foundation for the validity
of Prandtl layer theory \cite{Guo-Wang}. Our method is based on a
combination of the recent quotient estimate for the velocity field and an
introduction of a pseudo entropy estimate for the temperature field.
\end{abstract}

\date{\today}
\maketitle
	
\setcounter{tocdepth}{1}
\tableofcontents
	
\thispagestyle{empty}

\section{Introduction}
We study the steady compressible Navier-Stokes (NS) equations in the 2-dimensional domain $(x,Y)\in (0,L)\times \R_+$: 
{\small
	\begin{align}\label{1.1}
		\begin{cases}
			\dis \mbox{\bf div}(\rho^{NS} {U}^{NS})=0,\\[2mm]
			\dis  \rho^{NS} ( {U}^{NS}\cdot \mathbf{\nabla}) {U}^{NS} +   {\bf\nabla} p^{NS}   = \mu \, \v {\bf\Delta} {U}^{NS} + \lambda\, \v\, {\bf\nabla}\mbox{\bf div}{U}^{NS},\\
			\dis c_v \rho^{NS} (U^{NS}\cdot \nabla) T^{NS} + p^{NS} \mbox{\bf div}{U}^{NS} = \kappa \v {\bf\Delta} {T}^{NS} + 2\mu \v |S(U^{NS})|^2 + (\lambda-\mu) \v  |\mbox{\bf div}{U}^{NS}|^2,
		\end{cases}
\end{align}}
where $\mbox{\bf div}=\pa_x + \pa_Y$, $\nabla=(\pa_x, \pa_Y)^t$ and ${\bf\Delta}=\pa_{xx}+ \pa_{YY}$,  $p^{NS}=\rho^{NS}T^{NS} $, ${U}^{NS}=(u^{NS}, v^{NS})$,  $S(U)=\frac12 \big(\nabla U + \nabla U^T\big)$. Here the  parameter $0<\v\ll1$ is  the order of the inverse of Reynolds number,  $\mu$ and $\lambda$ are the {\it normalized} viscous coefficients, $\kappa$ is the {\it normalized} heat conductivity. We also normalize the specific heat $c_v=1$  for simplicity. In this paper, we always assume  that $\mu>0$ and $\kappa>0$.

For  \eqref{1.1}, we impose the no-slip boundary condition for velocity
\begin{align}\label{1.2}
	(u^{NS}, v^{NS})|_{Y=0}=0,
\end{align} 
and for the temperature either   Neumann  boundary condition  (NBC)
\begin{align}\label{1.3-0}
	\pa_Y T^{NS}|_{Y=0}&=0,
\end{align}
or Dirichlet boundary condition (DBC)
\begin{align}\label{1.3-1}
	T^{NS}|_{Y=0}&=\mathscr{T}_b(x)>0.
\end{align}

\medskip

\noindent{$\bullet$ \it The Euler background with mismatch:}  For $\v=0$ in \eqref{1.1}, we have the following Euler equations
\begin{align}\label{Euler}
	\begin{cases}
		\dis \mbox{\bf div}(\rho_e U_e)=0,\\
		\dis \rho_e (U_e\cdot\nabla)U_e + \nabla P_e=0,\\
		\dis \rho_e (U_e\cdot \nabla) T_e + P_e \mbox{\bf div}{U}_e =0,\\
		\dis v_e|_{Y=0}=0.
	\end{cases}
\end{align}
We assume a far-field  Euler shear flow as $[\rho_e^0(Y), u_e^0(Y), 0, T_e^0(Y)]$, where the functions $\rho_e^0(\cdot), u_e^0(\cdot)$, and $T_e^0(\cdot)$ are arbitrary smooth  satisfying
\begin{align}\label{1.4}
	p_e^0 = \rho_e^0(Y) T_e^0(Y) = 1.
\end{align}
It is clear that any such shear flow constitutes a solution to  the compressible Euler equations  \eqref{Euler}. For later use, we  assume 
\begin{align}\label{1.3}
	\begin{cases} 
		\dis 0<u_e^0(0)\leq u_e^0(Y) \leq u_e^+,   \\
		\dis  u_e^{+}:=\lim_{Y\rightarrow\infty} u_e^0(Y)>0, \\
		\dis |u_e^0(Y)-u_e^{+}|\leq |u_e^{+}-u^0_e(0)|  e^{-Y},  \\
		\dis |\pa_{Y}^k u_e^0(Y)|\lesssim [u_e^+-u^0_e(0)] e^{-Y},
	\end{cases}
	\& \quad 
	\begin{cases} 
		\dis\f12\leq  T_e^0(Y) \leq 2, \\
		\dis  1:=\lim_{Y\to +\infty}T^0_e(Y)>0,\\
		\dis|T_e^0(Y)-1| \lesssim |1-T_e^0(0)|e^{-Y},\\
		\dis |\pa_{Y}^k T_e^0(Y)|\lesssim [1-T^0_e(0)] e^{-Y}.
	\end{cases}
\end{align}
It is well-known that generically such an Euler flow has severe mismatches with the physical boundary conditions \eqref{1.2}-\eqref{1.3-1}:
\begin{align}\label{0.9}
u_e^0(0)\neq 0\equiv u^{NS}(x,0),\quad T_e^0(0)\neq T^{NS}(x,0).
\end{align}

\smallskip

\noindent{$\bullet$  \it Incompressible Prandtl boundary layer:} To resolve the mismatch for velocity, in 1904, Prandtl developed the Prandtl boundary layer theory for the simplified  incompressible Navier-Stokes equations (formally, the Mach number $\mathscr{M}_a=0$ and $\rho^{NS}=1$ in $\eqref{1.1}_{1,2}$):
\begin{align}
\begin{cases}
\mbox{\bf div} U^{NS}=0,\\
({U}^{NS}\cdot \mathbf{\nabla}) {U}^{NS} +   {\bf\nabla} p^{NS}   = \mu \, \v {\bf\Delta} {U}^{NS}.
\end{cases}
\end{align}
He derived the celebrated Prandtl equations as follows:
\begin{align}\label{1.10-P}
\begin{cases}
\pa_x\bar{u}_p + \pa_y \bar{v}_p=0,\\
\bar{u}_p \pa_x\bar{u}_p + \bar{v}_p \pa_y \bar{u}_p  + \pa_xp_{e}(x,0) =\mu \pa_{yy} \bar{u}_p,
\end{cases}
\end{align}
where the stretched variable $y:=\f{Y}{\sqrt{\v}}$. Such a Prandtl boundary layer characterized a sharp transition of the velocity field within a thin layer of $\sqrt{\v}$ near the boundary $\{Y=0\}$. The Prandtl theory states that in the incompressible fluids with small viscosity (i.e. large Reynolds number), the behavior near the boundary is described by the Prandtl boundary layer equations \eqref{1.10-P}, while the region away from the boundary is approximated by the inviscid Euler equations, i.e.,
\begin{align}\label{0.11}
	\mbox{Incompressible NS}=\mbox{Euler} + \mbox{Prandtl layer} + O(\sqrt{\v}).
\end{align}

\smallskip

\noindent{$\bullet$  \it Thermal boundary layer:}  
Clearly, the classical Prandtl layer fails to address the mismatch of the temperature in \eqref{0.9}, i.e., $T_e^0(0)\neq T^{NS}(x,0)$. Therefore the classical Prandtl boundary layer theory for incompressible fluids (formally, the Mach number $\mathscr{M}_a=0$ and $\rho^{NS}=1$ in $\eqref{1.1}_{1,2}$) with large Reynolds number
becomes inadequate to model compressible fluids with $\mathscr{M}_a\neq0$ and
a sharp temperature variation within a thin layer of a solid body (thermal
boundary layer) \cite{Schlichting}, at which the temperature
satisfies naturally either Neumann (isothermal body) or Dirichlet boundary
conditions. It is well-known that such a thermal boundary layer for the temperature field interacts
with the classical velocity Prandtl boundary layer in a strong nonlinear
fashion through the compressible Navier-Stokes equations,  playing a crucial role in many important applications in Aeronautics  and Astronautics.

\smallskip

We formally expand  the solution $(\rho^{NS},  u^{NS}, v^{NS}, T^{NS})$ of compressible Navier-Stokes equations \eqref{1.1} in $\v$:
{\small
	\begin{align}\label{3.4-0}
		\begin{split}
			\rho^{NS}&=\rho_e^0(Y)  + \rho_p^0(x,y) + \sum_{i=1}^{N} \sqrt{\v}^{i} [\rho_{e}^{i}(x,Y) + \rho_{p}^{i}(x,y)] + \v^{N_0}\rho^{(\v)}:=\rho_s + \v^{N_0}\rho^{(\v)},\\
			u^{NS}&= u_e^0(Y) + u_p^0(x,y) + \sum_{i=1}^{N} \sqrt{\v}^{i} [u_{e}^{i}(x,Y) + u_{p}^{i}(x,y)] + \v^{N_0}u^{(\v)} :=u_s + \v^{N_0}u^{(\v)},\\
			v^{NS}& =\sqrt{\v}\Big\{ [v_e^1(x,Y) + v_p^0(x,y)] + \sum_{i=1}^{N-1} \sqrt{\v}^{i} [v_{e}^{i+1}(x,Y) + v_{p}^{i}(x,y)] + \sqrt{\v}^{N} v_p^{N} + \v^{N_0} v^{(\v)}\Big\}\\
			&:=\sqrt{\v} \big\{v_s + \v^{N_0}v^{(\v)}\big\},\\
			T^{NS}& = T_e^0(Y) + T_p^0(x,y) + \sum_{i=1}^{N} \sqrt{\v}^{i} \, [T_{e}^{i}(x,Y) + T_{p}^{i}(x,y)] + \v^{N_0}T^{(\v)}:=T_s + \v^{N_0}T^{(\v)}.
		\end{split}
\end{align}}
First, the leading order Prandtl layer $(\rho^0_p, u^0_p, v^0_p, T^0_p)$ satisfies the nonlinear compressible Prandtl equations \eqref{3.5-0}-\eqref{3.6-0}. Next, $(\rho_e^i, u_e^i, v_e^i, T_e^i)$ ($1\leq i\leq N$) satisfy the compressible Euler layer equations \eqref{2.4} (or equivalently \eqref{2.23}-\eqref{2.35}), while $(\rho_p^i, u_p^i, v_p^i, T_p^i)$ ($1\leq i\leq N$) satisfy the linear compressible Prandtl layer \eqref{3.6-1} (or equivalently \eqref{3.6-1A}). The final layer $(\rho^{(\v)}, u^{(\v)}, v^{(\v)}, T^{(\v)})$ are called the remainders which depend on $\v$.

To understand the asymptotic behavior of $(\rho^{NS}, u^{NS}, v^{NS}, T^{NS})$, a key step is to construct the approximate terms, namely, the Prandtl layer $(\rho^0_p, u^0_p, v^0_p, T^0_p)$ and its high-order corrections $(\rho_e^i, u_e^i, v_e^i, T_e^i)$ \& $(\rho_p^i, u_p^i, v_p^i, T_p^i)$ for $1\leq i\leq N$. The goal of this paper is to study the thermal boundary layer theory for compressible fluids. Precisely, we aim to establish local well-posedness of Prandtl layers for steady compressible fluids and construct Prandtl layer expansion in strong Sobolev norms, a foundation for the validity of Prandtl layer theory \cite{Guo-Wang}.


\subsection{The Prandtl equations for compressible flows}
We consider the 2-dimensional steady compressible Prandtl equations  as the leading order in expansion \eqref{3.4-0}
\begin{align}\label{3.6-00}
	\begin{cases}
		\pa_y(\bar{\rho}^0_p \bar{T}^0_p)=0,\\
		(\bar{\rho}^0_p \bar{u}^0_p)_x + (\bar{\rho}^0_p \bar{v}^0_p)_y =0,\\
		\bar{\rho}^0_p[\bar{u}^0_p \, \pa_x + \bar{v}^0_p \, \pa_y] \bar{u}^0_p + (\bar{\rho}^0_p \bar{T}^0_p)_x = \mu \pa_{yy} \bar{u}^0_p, \\
		\bar{\rho}^0_p[\bar{u}^0_p \, \pa_x + \bar{v}^0_p \, \pa_y] \bar{T}^0_p + \bar{\rho}^0_p \bar{T}^0_p [\pa_x \bar{u}^0_p + \pa_y \bar{v}^0_p] = \kappa \pa_{yy}\bar{T}^0_{p} + \mu |\pa_y \bar{u}^0_p|^2,
	\end{cases} 
\end{align}
with $(x,y)\in (0,L)\times \mathbb{R}_+$, see Lemma \ref{lem2.1} in the appendix for  derivation of \eqref{3.6-00}.  The Prandtl equations \eqref{3.6-00} can be regarded as evolution equations in which $x$ acts as a time-like variable and  $y\in\mathbb{R}_+$ as a space-like variable, thus \eqref{3.6-00} constitutes a system of quasilinear parabolic equations. We impose the following initial and boundary conditions for \eqref{3.6-00}:
\begin{align}\label{3.6-01}
	\begin{cases}
		\dis (\bar{u}^0_p,  \bar{T}^0_p)|_{x=0}= (\mathscr{U}_0(y) ,\, \mathscr{T}_0(y)), \quad 
		(\bar{u}^0_p, \bar{v}^0_p)|_{y=0}=(0,0),\\
		\dis  \lim_{y\to \infty}(\bar{u}^0_p,  \bar{T}^0_p)=(u^0_e(0), T_e^0(0)),\\
		\dis \pa_y\bar{T}^0_p|_{y=0}=0\quad\quad  \mbox{for NBC},\\
		\dis \bar{T}^0_p|_{y=0}=\mathscr{T}_b(x)\quad \mbox{for DBC},
	\end{cases}
\end{align}
where $\mathscr{U}_0 ,\, \mathscr{T}_0$ are given  smooth  initial data, and $u^0_e(0), T_e^0(0)$  are two positive constants determined by the outer Euler shear flow, see \eqref{1.4}-\eqref{1.3}. 

\begin{theorem}\label{thm1}
	For any fix $\mathfrak{m}_0\gg1$. Assume  $\mathscr{U}_0 ,\, \mathscr{T}_0\in C^{2\mathfrak{m}_0+6}(\mathbb{R}_+)$, $\mathscr{T}_b\in C^{\mathfrak{m}_0+3}([0,1])$ satisfying 
	\begin{align}\label{0.1}
		\begin{split}
			\mathscr{U}_0'(0)>0,\quad 	\kappa \mathscr{T}_0''(0) + \mu |\mathscr{U}_0'(0)|^2=0 \quad \mbox{and} \quad \mathscr{U}_0''(0)=\mathscr{U}_0'''(0)=0.
		\end{split}
	\end{align}
	Assume also generic compatibility conditions at the corner up to order $2\mathfrak{m}_0+6$, see Section \ref{sec2.2.4} for details. Let $w_0=(1+ y)^{\mathfrak{l}_0}$ with $\mathfrak{l}_0\gg1$.  Assume 
	\begin{equation}\label{0.2}
	 {\bf I}[\mathscr{T}_b, \mathscr{U}_0, \mathscr{T}_0] :=\|\mathscr{T}_b\|_{C^{\mathfrak{m}_0+3}([0,1])}	+	\sum_{i=0}^{2\mathfrak{m}_0+6}\|(\pa_y^i(\mathscr{U}_0-u_e^0(0)), \pa_y^i(\mathscr{T}_0-T_e^0(0))) \, w_0\|_{L^\infty}<\infty.
	\end{equation}
	Then there exists a small $L_0>0$ depending on initial data ${\bf I}[\mathscr{T}_b, \mathscr{U}_0, \mathscr{T}_0]$ so that the smooth solution  $(\bar{\rho}^0_p, \bar{u}^0_p, \bar{v}^0_p, \bar{T}^0_p)$  of \eqref{3.6-00}-\eqref{3.6-01} exists in $x\in[0,L_0]$ satisfying
	\begin{align}\label{0.3}
		\begin{cases}
			\dis \bar{u}^0_p(y)>0\quad \mbox{for}\,\, y>0,\\ 
			\dis \sup_{x\in[0,L_0]}\sup_{y\in[0,y_0]} \bar{u}_{py} \geq \f12 \mathscr{U}_0'(0)>0\quad \mbox{for some}\,\, y_0\in \mathbb{R}_+,\\[2.5mm]
			\dis \pa_y^2\bar{u}^0_p\big|_{y=0}=\pa_y^3\bar{u}^0_p\big|_{y=0}=0,
		\end{cases}
	\end{align}
	and
	\begin{align}\label{0.4}
		&\sum_{i=0}^{m}\|\pa_x^{i}\bar{v}^0_p\|_{L^\infty} + 	\sum_{0< 2i+j\leq 2m}\|\pa_x^i \pa_y^{j}\bar{u}^0_p \, w_0\|_{L^\infty} + \sum_{0\leq 2i+j\leq 2m}\|\pa_x^i \pa_y^{j} \big(\bar{T}_p^0 -T_e^0(0)\big) \, w_0 \|_{L^\infty} \nonumber\\
		&\leq C\big({\bf I}[\mathscr{T}_b, \mathscr{U}_0, \mathscr{T}_0]\big),
	\end{align}
	for any $m\leq \mathfrak{m}_0$, and   $\bar{v}^0_p(0,\cdot)$ is defined as in \eqref{6.191-2}. Here $y_0>0$ is a constant  determined by $\mathscr{U}_0$. 
\end{theorem}

\subsection{Quotient and pseudo entropy estimates} 
Since the pioneering work of Oleinik \cite{Oleinik-1999}, there have been extensive studies on the Prandtl boundary layer theory for incompressible fluid, which is of parabolic type with $x$ being regarded as the time variable. It has been widely recognized that the non-locality of the normal velocity $v$ is out of control (even leading to a potential derivative loss) in the standard energy method due to the convection term $uu_{x}+vu_{y}$, which presents the single most difficulty in the PDE analysis of the Prandtl equations. In the case of steady flows, the classical von-Mises transform was employed in Oleinik \cite{Oleinik-1999} to circumvent such a difficulty, leading to a degenerate scalar parabolic reformulation. In order to obtain higher regularity, which is necessary for the Prandtl layer expansion, a quotient estimate $\dis\frac{v}{u}$ was introduced in \cite{Guo-Iyer-CMP} to circumvent the same difficulty. 

In the current study of compressible boundary layer problems with high
regularity, it is natural to introduce quotient $\dis q=\frac{\bar{v}}{\bar{u}_{\theta} }$,  so that the main part of
PDE is reformulated as%
\begin{align*}
	\bar{u}_{\theta}^{2}q_{xy}-\mu\bar{v}_{yyy}  & =\{\mu-\frac{\kappa}{2}\} \bar{u} \frac{\Theta_{xyy}}{T}\cdots, \\
	2\bar{u}_{\theta}\Theta_{x}+2\bar{v}\Theta_{y}  & =\kappa\partial_{yy} \Theta+\cdots. 
\end{align*}
for which the troublesome convection $uu_{x}+vu_{y}$ is eliminated. 
We perform an energy estimate with the multiplier $\partial_{x}^{k}q_{y} \, (k\leq m)$  on the first quotient equation (see Lemma \ref{lem6.8}), as well as an elliptic square estimate  
(see Lemma \ref{lem6.12}).

To control the strong nonlinear interaction
\begin{align*}
	\Big\la\{\mu-\frac{\kappa}{2}\} \bar{u} \frac{\partial_{x}^{m}\Theta_{xyy}}{T},\, \partial_{x}^{m}q_y \Big\ra
	\lesssim  L^{\f16-}\|\bar{u}\pa_x^{m+1}\Theta_{yy}\|\cdot \|\pa_x^mq_y\|_{\MX},
\end{align*}
one must first establish the energy estimate for the second temperature
equation with the multiplier $\partial_{x}^{k}\Theta\, (k\leq m)$, then complete with an
elliptic square estimate. As expected, in the estimate of the highest order derivative of $\pa_x^{m+1}$,  the normal velocity presents a basic
difficult contribution of
\begin{align*}
	\big\la \partial_{x}^{m+1}\bar{v}\Theta_{y},\, \partial_{x}^{m+1}\Theta\big\ra= \big\la\bar
	{u}_{\theta}\partial_{x}^{m+1}q\Theta_{y},\, \partial_{x}^{m+1}\Theta\big\ra + \cdots,
\end{align*}
for which $\partial_{x}^{m+1}q$ is out of control from either the energy
control of $\sqrt{\bar{u}_{\theta}}\partial_{x}^{m}q_{yy}$ (loss $\partial_{x}$) or the elliptic norm of $\bar{u}_{\theta}^{3/2}\partial
_{x}^{m+1}q_{y}$ (degeneracy at $y=0,$ note $\dis \lim_{y\rightarrow\infty} 
\partial_{x}^{m+1}q\neq0$). We note that no quotient estimate for $\Theta$ is
possible due to boundary condition $\Theta|_{y=0}\neq0.$ To circumvent such
well-known difficulty, we introduce a key new good unknown, the \textit{pseudo entropy }$S_{k}$, as
\begin{align*}
	S_{k}  & := \partial_{x}^{k}\Theta+\partial_{x}^{k-1}q\,\Theta_{y}, 
\end{align*}
which yields the equation of $S_{m+1}$ as
\begin{align*}
	2 \bar{u}_{\theta} \pa_x S_{m+1} +  2 \bar{v} \pa_yS_{m+1} - \kappa \pa_y^2 S_{m+1} 
	= 2 \bar{u}_{\theta} \Theta_{xy} \pa_x^{m}q  + 2 \bar{v} \big(\Theta_y\pa_x^{m}q\big)_y - \kappa \big(\Theta_y\pa_x^{m}q\big)_{yy}  + \cdots,
\end{align*}
for which the out of control part $\bar{u}_{\theta}\partial_{x}^{m+1} q\Theta_{y}$ is included as part of the new pseudo entropy $S_{m+1}.$ In fact, we note 
\begin{align*}
	\|\sqrt{\bar{u}_{\t}}S_{m+1}\|_{x=x_0} &= \|\sqrt{\bar{u}_{\t}}\big[\pa_x^{m+1}\Theta + \pa_x^{m}q\, \Theta_y\big]\|_{x=x_0}, 
\end{align*}
and
\begin{align*}
	\|\bar{u}^2_{\t} \pa_x S_{m+1} \|^2
	+   \|\bar{u}_{\t}^{\f32} \pa_yS_{m+1}\|^2_{x=x_0} 
	&=\|\bar{u}^2_{\t}  \big[\pa_x^{m+2}\Theta + \pa_x^{m+1}q\, \Theta_y + \pa_x^{m}q\, \Theta_{xy}\big]\|^2 \nonumber\\
	&\quad   + \|\bar{u}^{\f32}_{\theta}  \big[\pa_x^{m+1}\Theta_{y} + \pa_x^{m}q_y\, \Theta_y + \pa_x^{m}q\, \Theta_{yy}\big]\|^2_{x=x_0},
\end{align*}
we point out that the terms $\|\bar{u}^2_{\theta}   \pa_x^{m+1}q\, \Theta_y\|$, $\|\bar{u}^{\f32}_{\t}\pa_x^{m}q\, \Theta_y\|_{x=x_0}$, $\|\bar{u}^{\f32}_{\t}\pa_x^{m}q\, \Theta_{yy}\|^2_{x=x_0}$ are out of control. 
On the other hand, fortunately, all the rest contributions with lower derivatives $\partial
_{x}^{j}q$ ($j\leq m$) can be controlled  in the new
reformulation of the problem in term of such a pseudo entropy, see \eqref{6.190}-\eqref{6.238} for details.
Then it is possible to close both the velocity quotient estimates and pseudo entropy
estimates via careful rearrangements.

\smallskip

For subsequent higher boundary layer constructions, we employ the stream
function $\psi$ and its quotient 
$\dis q=\frac{\bar{T}^{0}_{p}\psi}{\bar{u}_{p,\theta}^{0}}$ accordingly, as well as modify the pseudo entropy as
\begin{align*}
S_{k}:=\partial_{x}^{k}\Theta + \partial_{y}\bar{T}_{p}^{0}\partial_{x}^{k-1}q,
\end{align*}
see Section \ref{sec3} for details.

\subsection{Higher-order linear compressible Prandtl layer equations}
To solve the higher-order linear compressible Prandtl layer equations (see \eqref{3.6-1} or equivalently \eqref{3.6-1A} in the appendix for its derivation) in the expansion \eqref{3.4-0}, we consider the following linear system  for the unknowns $(\rho_p, u_p, v_p, T_p)$:
\begin{align}\label{3.6-10}
	\begin{cases}
		\big(\bar{\rho}^0_p T_p + \bar{T}^0_p \rho_p\big)_y = G_{1},	\\[1.5mm]
		\big(\bar{\rho}^0_p u_p + \rho_p \bar{u}^0_p \big)_x + \big(\bar{\rho}^0_pv_p  + \rho_p \bar{v}^0_p \big)_y = G_{2},\\[1.5mm]
		\bar{\rho}^0_p  (\bar{U}^0_p\cdot \bar{\nabla}) u_p + \bar{\rho}^0_p(U_p\cdot \bar\nabla) \bar{u}^0_p + \rho_p \big(\bar{U}^0_p\cdot \bar{\nabla}\big) \bar{u}^{0}_p  +   \big(\bar{\rho}^0_p T_p + \bar{T}^0_p \rho_p \big)_x 
		=\mu \pa_{yy} u_p + G_{3},\\[1.5mm]
		\bar{\rho}^0_p(\bar{U}^0_p \cdot \bar{\nabla})T_p + \bar{\rho}^0_p(U_p\cdot \bar\nabla) \bar{T}^0_p   + \rho_p\, (\bar{U}^0_p\cdot \bar{\nabla}) \bar{T}^{0}_p +  \big\{\bar{\rho}_p^0 T_p + \bar{T}^0_p \rho_p\big\} \mbox{\rm div} \bar{U}^0_p   \\[1mm]
		\qquad\qquad\qquad\qquad\qquad\quad \,\,  +  \bar{\rho}^0_p \bar{T}^0_p  \big\{\pa_x u_p + \pa_y v_p \big\}
		= \kappa \pa_{yy}T_p  + 2 \mu \pa_y u_p \cdot \pa_y \bar{u}^{0}_p  + G_{4},
	\end{cases}
\end{align}
with $(x,y)\in (0,L)\times \mathbb{R}_+$, 
and we have denoted $\bar{U}^0_p:=(\bar{u}^0_p, \bar{v}^0_p)^{t}$,  $U_p:=(u_p, v_p)^{t}$, $\mbox{\rm div} := \pa_x + \pa_y$, $\bar{\nabla}:=(\pa_x,\pa_y)^t$. Here $(\bar{\rho}^0_p, \bar{u}^0_p, \bar{v}^0_p, \bar{T}^0_p)$ are solution of nonlinear Prandtl problem \eqref{3.6-00}-\eqref{3.6-01} established in Theorem \ref{thm1}, while the forcing terms $G:=(G_1, G_2, G_3, G_4)$ are given functions with enough regularity and space decay.  The following initial boundary conditions are imposed on \eqref{3.6-10}:
\begin{align}\label{3.6-11}
	\begin{cases}
		(u_{p}, T_{p})|_{x=0}=\big(\tilde{u}_{0}(y), \tilde{T}_{0}(y)\big), \quad \lim_{y\to \infty} (u_p, T_p)=(0,0), \\
		(u_{p}, v_{p})|_{y=0} = \big(\tilde{u}_b(x), 0\big),\\[1mm]
		\pa_yT_p |_{y=0} =\tilde{T}_b(x) \quad \mbox{for NBC},\\[1mm]
		T_{p} |_{y=0} =\tilde{T}_b(x) \qquad\, \mbox{for DBC}.
	\end{cases}
\end{align}
By using similar quotient and pseudo entropy estimates as in the proof of Theorem \ref{thm1}, we establish:
\begin{theorem}\label{thm2}
Recall $\mathfrak{m}_0, l_0$ in Theorem \ref{thm1}.  Let $w_1=(1+ y)^{\mathfrak{l}_1}$ with $\mathfrak{l}_0\gg \mathfrak{l}_1 \gg1$.   For any fix $\mathfrak{m}_0\gg \mathfrak{n}_0\gg1$. Assume  $\tilde{u}_0 ,\, \tilde{T}_0\in C^{2\mathfrak{n}_0+8}(\mathbb{R}_+)$, $\tilde{u}_b ,\, \tilde{T}_b\in C^{\mathfrak{n}_0+4}([0,1])$ with
\begin{align}\label{0.17}
\tilde{\bf I}[\tilde{u}_b, \tilde{T}_b, \tilde{u}_0, \tilde{T}_0]:=\|(\tilde{u}_b, \tilde{T}_b)\|_{C^{\mathfrak{n}_0+4}([0,1])}	+	\sum_{i=0}^{2\mathfrak{n}_0+8}\|(\pa_y^i\tilde{u}_0, \pa_y^i\tilde{T}_0) \, w_1\|_{L^\infty}<\infty.
\end{align}
Assume the forcing terms $G=(G_1, G_2, G_3, G_4)$ satisfy
\begin{align}\label{0.18}
\hat{\bf I}[G]:=\sum_{0\leq 2i+j\leq 2 \mathfrak{n}_0+8} \|\pa_x^i\pa_y^j G w_1\|_{L^2} <\infty.
\end{align}
Assume also generic compatibility conditions on the initial boundary data at the corner up to order $2\mathfrak{n}_0+8$, see Section \ref{sec3.2} for details.
Then there is a suitably small constant $L_1\in (0, L_0)$ such that the solution of \eqref{3.6-10}-\eqref{3.6-11} exists in $x\in[0,L_1]$ with
\begin{align}\label{1.10}
&\sum_{i=0}^{k}\|\pa_x^{i}v_p\|_{L^\infty} + \sum_{0\leq 2i+j\leq 2k} \|(\pa_x^i \pa_y^{j}u_p, \, \pa_x^i \pa_y^{j}T_p)  w_1\|_{L^\infty}
\lesssim C\big(\tilde{\bf I}[\tilde{u}_b, \tilde{T}_b, \tilde{u}_0, \tilde{T}_0], \, \hat{\bf I}[G]\big),
\end{align}
for any $k\leq \mathfrak{n}_0$. 
\end{theorem}

\smallskip

\subsection{Construction of approximate compressible Navier-Stokes solution}
Combining Theorems \ref{thm1} and \ref{thm2} with Proposition \ref{prop5.7}, we have
\begin{theorem}\label{thm3}
Assume the shear flow $(\rho^0_e, u^0_e, 0, T^0_e)(Y)$ is sufficiently  smooth with \eqref{1.4}-\eqref{1.3} and the uniform subsonic condition \eqref{1.3-6}.
Let $(\bar{\rho}^0_p, \bar{u}^0_p, \bar{v}^0_p, \bar{T}^0_p)$ be the solution of Prandtl problem \eqref{3.6-00}-\eqref{3.6-01} established in Theorem \ref{thm1}. We denote the weights $\mathfrak{w}=\la y\ra^{\mathfrak{b}}$ and $\bar{\mathfrak{w}}=\la Y\ra^{\mathfrak{b}}$. Let $\mathfrak{l}_0\gg \mathfrak{b}\gg N$ and $\mathfrak{m}_0\gg K_0\gg N$. For the initial data in \eqref{2.4-1} and \eqref{3.6-2A}, we assume 
\begin{align}
\mathfrak{C}_{N, 2K_0}:=\sum_{k=1}^{N}\sum_{j=0}^{2K_0}\big\{\|\pa_Y^j(\tilde{v}_{e,0}^k, \, \tilde{v}_{e,L}^k,\, \tilde{\rho}^k_e, \tilde{u}^k_e, \tilde{T}^k_e) \bar{\mathfrak{w}}\|_{L^\infty_Y} +  \|\pa_y^j(\tilde{u}^k_p, \tilde{T}^k_p) \mathfrak{w}\|_{L^\infty_y}\big\}<\infty.
\end{align}
We assume  generic compatibility conditions on the initial data $\tilde{u}^k_p, \tilde{T}^k_p$ at the corner up to order $2K_0$, see Section \ref{sec3.2} for details.
At the corners $(0,0)$ and $(L,0)$, we also assume the standard elliptic compatibility conditions on $\tilde{v}_{e,0}^k, \, \tilde{v}_{e,L}^k$, see Section 4 for details. 

Then there is a small $L\in (0, L_0)$ such that the solutions  $(\rho_e^k, u_e^k, v_e^k, T_e^k)$ of \eqref{2.23}-\eqref{2.4-1}  and the solutions $(\rho_p^k, u_p^k, v_p^k, T_p^k)$ of \eqref{3.6-1A}-\eqref{3.6-2A}  exist and are smooth on $\Omega=[0,L]\times \mathbb{R}_+$ for $k=1,\cdots, N$.  The following estimates hold:
\begin{align}\label{0.26}
&\sum_{k=1}^{N}\sum_{j=0}^{K_0} \Big\{\|\nabla^j(\rho_e^k, u_e^k, v_e^k, T_e^k) \sqrt{\bar{\mathfrak{w}}}\|_{L^\infty} + \|\pa_x^j v^k_p\|_{L^\infty} + \|\bar{\nabla}^{j} (\rho^k_p, u^k_p, T^k_p) \sqrt{\mathfrak{w}} \|_{L^\infty} \Big\} \nonumber\\
&\lesssim C\big({\bf I}[\mathscr{T}_b, \mathscr{U}_0, \mathscr{T}_0], \mathfrak{C}_{N, 2K_0}\big).
\end{align}
\end{theorem}

\smallskip

\begin{remark}
	With the help of Theorem \ref{thm3}, we can define the approximate solution $(\rho_s, u_s, v_s, T_s)$ {\rm(}as given in \eqref{3.4-0}{\rm)} to the steady compressible Navier-Stokes equations such that both the remainder terms and the boundary errors exhibit high order $\sqrt{\v}$-decay. For the use in \cite{Guo-Wang}, we introduce a new approximate solution $(\tilde{\rho}_s, \tilde{u}_s, \tilde{v}_s, \tilde{T}_s)$ with $\sqrt{\v}^N$ corrections, where $(\tilde{u}_s, \tilde{v}_s, \tilde{T}_s)$ satisfy exactly \eqref{1.2}-\eqref{1.3-1} and also possess some other good properties; see Appendix B for details.
\end{remark}

\begin{remark}
Now we would like to derive some smallness property if the initial data also own smallness. Let $\mathscr{W}$ be any polynomial weight on $y$, $\bar{\mathscr{W}}$ be any polynomial weight on $Y$. For instance
\begin{align}
\mathscr{W}= \la y \ra^{\mathfrak{b}_1},\,\, \la \sigma y\ra^{\mathfrak{b}_1} \quad \mbox{and}\quad \bar{\mathscr{W}}(Y)=\la Y \ra^{\mathfrak{b}_1},\,\, \la \sigma Y\ra^{\mathfrak{b}_1},
\end{align}
where $\sigma>0$ is some parameter, and $\mathfrak{b}_1\leq \f{1}{2}\mathfrak{b}$.
Let $(\rho_p^0, u_p^0, v_p^0, T_p^0)$ be the solution of Prandtl established in Theorem \ref{thm1}, while $(\rho_e^1, u_e^1, v_e^1, T_e^1)$  and $(\rho_p^1, u_p^1, v_p^1, T_p^1)$ be the solutions established in  Theorem \ref{thm3}.  

For any $i+j\leq \f12 K_0 $, we have 
\begin{align*}
\begin{split}
\pa_x^i\pa_y^j (\bar{\rho}^0_p, \bar{u}^0_p, \bar{v}^0_p, \bar{T}^0_p)(x,y) \mathscr{W}(y)&= \pa_x^i\pa_y^j (\bar{\rho}^0_p, \bar{u}^0_p, \bar{v}^0_p, \bar{T}^0_p)(0,y) \,\mathscr{W}(y)+ O ({\bf I}[\mathscr{T}_b, \mathscr{U}_0, \mathscr{T}_0]) L,\\
\pa_x^i\pa_y^j (\rho_p^1, u_p^1, v_p^1, T_p^1)(x,y) \mathscr{W}(y) &= \pa_x^i\pa_y^j (\rho_p^1, u_p^1, v_p^1, T_p^1)(0,y) \, \mathscr{W}(y) + C\big({\bf I}[\mathscr{T}_b, \mathscr{U}_0, \mathscr{T}_0]\big) L,\\
\pa_x^i\pa_y^j (\rho_e^1, u_e^1, v_e^1, T_e^1)(x,Y) \, \bar{\mathscr{W}}(Y) &= \pa_x^i\pa_y^j (\rho_e^1, u_e^1, v_e^1, T_e^1)(0,Y) \, \bar{\mathscr{W}}(Y) + O({\bf I}[\mathscr{T}_b, \mathscr{U}_0, \mathscr{T}_0], \mathfrak{C}_{1, 2K_0})L,
\end{split}
\end{align*}
which yield that the behaviors of $(\rho_p^0, u_p^0, v_p^0, T_p^0)$, $(\rho_e^1, u_e^1, v_e^1, T_e^1)$, and $(\rho_p^1, u_p^1, v_p^1, T_p^1)$ are completely determined by the corresponding initial data. Therefore, if the initial data possess smallness and decay properties {\rm(}see Remarks \ref{rem2.1} $\&$ \ref{rem2.2} for instances on the construction of such initial data{\rm)}, the solution will also exhibit such properties.

Noting \eqref{3.4-0}, we emphasize that we need only to derive the smallness properties for $(\bar{\rho}_p^0, \bar{u}_p^0, \bar{v}_p^0, \bar{T}_p^0)$, $(\rho_e^1, u_e^1, v_e^1, T_e^1)$ and $(\rho_p^1, u_p^1, v_p^1, T_p^1)$ because additional $\sqrt{\v}$ coefficients are available for $(\rho_e^k, u_e^k, v_e^k, T_e^k)$ and $(\rho_p^k, u_p^k, v_p^k, T_p^k)$ with $2\leq k\leq N$. The smallness and decay properties will play a key role in for the validity
of Prandtl layer theory \cite{Guo-Wang}. 
\end{remark}

\smallskip

\subsection{Literature review}
For the steady incompressible Prandtl equations, Oleinik \cite{Oleinik1} established a fundamental theorem concerning the existence and uniqueness of strong solutions through the von Mises transformation and the maximum principle, and further demonstrated that the solution is global-in-$x$ under favorable pressure condition. Recently, Guo-Iyer \cite{Guo-Iyer-CMP} derived higher regularity via quotient estimates and energy methods, leading to a Prandtl layer expansion up to any order. Iyer \cite{Iyer4} proved the asymptotic stability near the self-similar Blasius profile, while Iyer-Masmoudi \cite{Iyer-M0} established global-in-$x$ regularity of the Prandtl system without invoking the classical von Mises transformation. Wang-Zhang \cite{WZ} obtained global $C^\infty$ regularity under favorable pressure gradients using the maximum principle method. Furthermore, Guo-Wang-Zhang \cite{GWZ} proved the orbital and asymptotic stability of Blasius-like steady states against Oleinik’s monotone solutions; see also the very recent result \cite{WZ1}.

It is a challenge problem to rigorously justify the Prandtl expansion \eqref{0.11}  since it is indeed a multi-scales problem. Recently, significant progress has been made for the steady incompressible Navier-Stokes equations. For the case of external force and periodic-in-$x$, Gerard-Varet and Maekawa \cite{GM2} provided a  sobolev stability of Prandtl expansions for the steady Navier–Stokes equations. Guo-Iyer \cite{Guo-Iyer-2024}  justified the validity of the steady Prandtl layer expansion for scaled Prandtl layers, including the celebrated Blasius boundary layer. Subsequently, Gao-Zhang \cite{GZ} established the expansion over a large interval in $x$. Iyer-Masmoudi \cite{Iyer-M} proved the global-in-$x$ stability of steady Prandtl expansions for the 2D incompressible Navier-Stokes equations. For other contributions, see \cite{Guo-N,Iyer,Iyer1,Iyer2,FM,FGLT1,FGLT2} for steady incompressible Navier-Stokes flows with moving boundaries.

Recently, a mathematical framework \cite{Guo-Wang} has been developed by the authors to justify the Prandtl expansion \eqref{0.11} for the steady compressible fluids, for which our current study provides a necessary foundation.

On the other hand, for the {\it time-dependent} problem of incompressible flow, extensive and important progress has been made since the seminal work of Oleinik \cite{Oleinik1,Oleinik-1999}, we refer to \cite{AWXY,Constantin,E,E1, GD,GM1, GMM, GreN0,GGN1,GGN,GreN1,GreN2, GN4,IV, Kato,Kelliher,KVW, LMY,Maekawa,MM,MW,PZ,SC1,SC2,Temam,WWZ,XZ} and the references therein. For the steady case of compressible flows, we refer to recent works \cite{LYZ,LWY,YZ}.

\smallskip


\noindent{\it Organization of the present paper:} In Section \ref{sec2}, we establish the local existence and uniform-in-$L$ estimates for the nonlinear Prandtl boundary layer problem \eqref{3.6-00}--\eqref{3.6-01}. Section \ref{sec3} presents the local existence and uniform-in-$L$ estimates for the steady linear Prandtl layer problem \eqref{3.6-10}--\eqref{3.6-11}. The well-posedness of the linear Euler layer problem \eqref{2.23}--\eqref{2.4-1} is established in Section \ref{sec4}. With the above preparations, we prove Theorem \ref{thm3} in Section \ref{sec5}, i.e., we construct all profiles in the approximate solution $(\rho_s, u_s, v_s, T_s)$ of the compressible Navier-Stokes equations \eqref{1.1}--\eqref{1.3-1}. In  Appendix \ref{App-A}, we derive the nonlinear Prandtl layer equations , the compressible Euler layer equations, and the steady linear Prandtl layer equations, along with the boundary conditions for these three systems. In  Appendix \ref{App-B}, we present the approximate solution of compressible Navier-Stokes equations with desired estimates.
 
\section{Existence of the Compressible Prandtl Layer}\label{sec2}
In this section, we aim to establish the existence and uniform estimates of the solution to \eqref{3.6-00}-\eqref{3.6-01}. For simplicity of presentation, we drop the superscript and subscript of $(\bar{\rho}^0_p, \bar{u}^0_p, \bar{v}^0_p, \bar{T}^0_p)$, and use $(\rho, u, v, T)$ as the unknown functions. Then we reduce \eqref{3.6-00} as
\begin{align}\label{3.58}
	\begin{cases}
		\dis  (\f{1}{T}u)_x + (\f{1}{T}v)_y =0,\\[3.5mm]
	\dis \f{1}{T} \big[u\pa_x u + v\pa_y u\big] = \mu \pa_{yy} u, \\[3.5mm]
	\dis \f{2}{T} \big[u\pa_xT + v\pa_y T\big]   = \kappa \pa_{yy}T + \mu |\pa_y u|^2.
	\end{cases}
\end{align}
where we have used the fact $\rho  T =\rho^0_e(0)\, T^0_{e}(0)=1$, i.e., $\dis \rho=T^{-1}$.
Then the initial and boundary conditions of \eqref{3.58} become
\begin{align}\label{6.2}
	\begin{split}
		u(0,y)=\mathscr{U}_0(y),\quad T(0,y)=\mathscr{T}_0(y), \quad \lim_{y\to\infty}u=u^0_e(0),  \quad \lim_{y\to\infty} T =T^0_e(0),\\
		u|_{y=0}= v|_{y=0}=0, \quad \pa_yT|_{y=0}=0 \,\, \mbox{for NBC} \, \big(\mbox{or}\,\,T|_{y=0} = \mathscr{T}_b(x)\,\, \mbox{for DBC}\big),
	\end{split}
\end{align}
with \eqref{0.1} holding.

Define
\begin{align}\label{3.61}
	\bar{U}:=(\bar{u},\bar{v})^t=(\f{1}{T}u, \f{1}{T}v)\quad \mbox{and}\quad \bar{u}_{\theta}=\theta + \bar{u}.
\end{align}
Then the system \eqref{6.1} is equivalent to  the following
\begin{align}\label{5.3-0}
	\begin{cases}
		\dis \bar{u}_x + \bar{v}_y =0,\\ 
		\dis \bar{u} \pa_x\bar{u} + \bar{v}\pa_y  \bar{u}  = \mu \pa_{yy} \bar{u} + \mathcal{N}_1 ,\\
		\dis 2\big[\bar{u} \pa_xT + \bar{v}\pa_y T\big]   = \kappa \pa_{yy}T +  \mathcal{N}_2,
	\end{cases}
\end{align}
where we have used the following notation
\begin{align}
	\begin{split}
		\mathcal{N}_1(\bar{u},T):&=\big[\mu-\f{\kappa}{2}\big]\f{T_{yy}}{T} \bar{u}  + \f{2\mu}{T} T_{y}\bar{u}_y- \f{\mu}{2T}  T_y^2\bar{u}^3 - \mu T_y \bar{u}_y\bar{u}^2 - \f{1}{2}\mu T\bar{u}^2_y\bar{u},\\
		\mathcal{N}_2(\bar{u},T):&=\mu |\pa_y (T\bar{u})|^2.
	\end{split}
\end{align}

\subsection{Compatibility conditions of initial data}\label{sec2.2.4}

We always assume enough smoothness and decay property for the initial data \eqref{6.2}. In order to obtain higher regularity, we must assume some compatibility conditions for the initial data.

Noting $u=T\bar{u}$, 
we have from $\eqref{5.3-0}_{1,2}$ that
\begin{align}\label{6.191-1}
	\big(\f{\bar{v}}{\bar{u}}\big)_y  = - \f{1}{T\bar{u}^2} \big\{\mu \pa_{yy}u -   \f1{2T}u\big[\kappa \pa_{yy}T + \mu |\pa_yu|^2\big]\big\},
\end{align}
which yields that 
\begin{align}\label{6.191-2}
	\begin{split}
		&\f{\bar{v}}{\bar{u}} = - \int_0^y \f{1}{T\bar{u}^2} \big(\mu  u_{yy} -   \f1{2T}u\big[\kappa T_{yy} + \mu |u_y|^2\big]\big)(z) dz\\
		&\Longleftrightarrow\quad \bar{v} =-\bar{u} \int_0^y \f{1}{T\bar{u}^2} \big(\mu  u_{yy} -   \f1{2T}u\big[\kappa T_{yy} + \mu |u_y|^2\big]\big)(z) dz.
	\end{split}
\end{align}
Then we can define $\bar{v}|_{x=0}$ and
\begin{align}\label{6.191-0}
	\begin{split}
		\pa_x\bar{u}|_{x=0}&:=-\bar{v}_y|_{x=0},\\
		\pa_xT|_{x=0}&:= \f{1}{2\bar{u}} \big\{\kappa T_{yy}+ \mu |u_y|^2 - 2\bar{v} T_y \big\}.
	\end{split}
\end{align}
Noting \eqref{0.1}, we remark that  \eqref{6.191-2}-\eqref{6.191-0} are well-defined.

For later use, we fix a smooth cut-off function 
\begin{align}\label{cutoff}
\chi(s):=
\begin{cases}
1,\quad s\in[0,1),\\
\mbox{smooth monotonic},\quad s\in [1, 2],\\
0,\quad s\in (2,\infty),
\end{cases}
\end{align}
and $\bar{\chi}(s):=1-\chi(s)$.

Next, we provide an example of the initial data $(\mathscr{U}_0, \mathscr{T}_0)$ by taking 
\begin{align}
\mathscr{U}_0(y)= y \chi(y) + \bar{\chi}(y),
\end{align}
and $\mathscr{T}_0$ such that 
\begin{align}
	\begin{cases}
		\kappa \mathscr{T}_0''(y) + \mu |\mathscr{U}_0'(y)|^2=0,\quad \mbox{for}\,\, y\in [0,1],\\
		\mathscr{T}_0'(0)=0 \,\,\, \mbox{for NBC} \,\, \big(\mbox{or}\,\, \mathscr{T}_0(0)=\mathscr{T}_b(0) \,\, \mbox{for DBC}\big),\\
		|\mathscr{T}_0(y)-T_e^0(0)|\lesssim e^{-y}\quad \mbox{as}\,\, y\to \infty.
	\end{cases}
\end{align}
It is clear to check that \eqref{6.191-2}-\eqref{6.191-0} are well defined. Inductively, assuming high order compatibility conditions, we can also define 
\begin{align}\label{2.12}
\pa_x^k\big(\f{\bar{v}}{\bar{u}}\big)_y\big|_{x=0},\quad  
\pa_x^k\bar{u}\big|_{x=0}, \quad    \pa_x^k\bar{v}\big|_{x=0} \quad \mbox{and}\quad \pa_x^kT\big|_{x=0} \quad \mbox{for}\,\,\, k\geq 0.
\end{align}
The details are omitted for simplicity of presentation.

\medskip

In \cite{Guo-Wang}, we need some smallness properties for the solution of Prandtl equations \eqref{3.58}. Now we give two remarks to explain how to construct initial data with some smallness.
\begin{remark}[Low Mach number fluid]\label{rem2.1}
Motivated by the Blasius profile, we can also construct a class of data in the following way. For instance, 
let $\hat{u}(y)$ be given function with the following properties
\begin{align}\label{I.14}
\hat{u}(y)=
\big(y+y^{k_0}g_1(y)\big)\, \chi(y) + \bar{\chi}(y), 
\end{align}
where $g_1 \in C^{k_1}([0,\infty)), k_1\gg1, k_0\geq 4$. Let $\hat{T}(y) \in C^{k_1}([0,\infty))$ satisfy
\begin{align}\label{I.15}
	\begin{cases}
		\kappa \hat{T}''(0) + \mu |\hat{u}'(0)|^2=0,\\
		\hat{T}'(0)=0 \,\,\, \mbox{for NBC} \,\, \big(\mbox{or given}\,\, \hat{T}(0) \,\, \mbox{for DBC}\big),\\
		|\hat{T}(y)|\lesssim e^{-y}\quad \mbox{for}\,\, y\geq 1.
	\end{cases}
\end{align}
Let $\s>0$ be the amplitude of Prandtl layer, and denote
\begin{align}\label{I.16}
\begin{cases}
\mathscr{U}^\s_0(y)=\s \hat{u}(\sqrt{\s}y),\\
\mathscr{T}^\s_0(y)=1+ \s^2 \hat{T}(\sqrt{\s}y).
\end{cases}
\end{align}
Applying \eqref{I.14}-\eqref{I.16}, one obtains 
\begin{align}
|\bar{v}^\s(y) |
&\lesssim \sqrt{\s}\hat{u}(\sqrt{\s}y) \int_0^{\sqrt{\s}y} \f{1}{|\hat{u}(z)|^2} \Big(\mu  |\hat{u}''(z)| + \sigma^2 \hat{u}(z) \big|\big[\kappa \hat{T}''(z) + \mu  |\hat{u}'(z)|^2\big]\big|\Big)  dz\nonumber\\
&\lesssim \sqrt{\s}\hat{u}(\sqrt{\s}y)  \min\{1, (\sqrt{\s}y)^{k_0-3}\} = \sqrt{\s} |\hat{u}(\sqrt{\s}y)|^{k_0-2}.
\end{align}
Then we can check that 
\begin{align}
\begin{split}
\begin{cases}
\dis \pa_x\bar{u}^\s|_{x=0}:= \pa_y \bar{v}^\s\cong \s |\hat{u}(\sqrt{\s}y)|^{k_0-3} e^{-\sqrt{\s}y}, \\[1mm]
\dis \pa_xT^\s|_{x=0}:=\f{1}{2\bar{u}^\s} \big\{\kappa (\mathscr{T}^{\s}_0)''+ \mu |(\mathscr{U}^{\s}_0)'|^2 - 2\bar{v}^{\s} ( \mathscr{T}_0)' \big\} \cong \s^2 e^{-\sqrt{\s}y}.
\end{cases}
\end{split}
\end{align}
Thus we can inductively define
\begin{align}
\pa_x^k\big(\f{\bar{v}^{\s}}{\bar{u}^{\s}}\big)_y\big|_{x=0},\quad  
\pa_x^k\bar{u}^{\s}\big|_{x=0},\quad  \pa_x^k\bar{v}^{\s}\big|_{x=0} \quad \mbox{and}\quad \pa_x^kT^{\s}\big|_{x=0} \quad \mbox{for}\,\,\, k\geq 1.
\end{align}
It is clear that we have smallness property for the initial data of \eqref{I.16} provided $\s>0$  is suitably small.
\end{remark}

\begin{remark}[High subsonic flow]\label{rem2.2}
We introduce a class of initial data in the following form 
\begin{align}\label{I.20}
(\mathscr{U}_0^{\s}, \mathscr{T}_0^{\s})(y)=(\hat{u}(\s y), 1+\hat{T}(\s y)),
\end{align}
where $(\hat{u}, \hat{T})$ are the functions in \eqref{I.14}-\eqref{I.15}, and $\sigma>0$ is a small parameter. Then it is clear that \eqref{I.20} is high  subsonic flow, and varies very slowly.

It is direct to check  that $\mathscr{U}_0^{\s}$ and $\mathscr{T}_0^{\s}$ satisfy \eqref{0.1}. Similarly as in \eqref{6.191-1}-\eqref{6.191-2}, we know 
\begin{align}
\bar{v}^{\s} &=-\bar{u}^{\s} \int_0^y \Big(\f{1}{T^{\s}|\bar{u}^{\s}|^2} \Big\{\mu \pa_y^2 u^{\s} -   \f1{2T^{\s}}u^{\s}\big[\kappa \pa_y^2T^{\s}  + \mu |\pa_yu^{\s}|^2\big]\Big\}\Big)(z) dz\nonumber\\
&\cong \s \hat{u}(\s y) \int_0^{\s y} \Big(\f{1}{|\bar{u}|^2} \Big\{\mu |\hat{u}''| +  \hat{u}(z)\big[\kappa \hat{u}''(z) + \mu |\hat{u}'|^2\big]\Big\}\Big)(z) dz \cong \sigma |\hat{u}(\sigma y)|^2.
\end{align}
Then we know that 
\begin{align}\label{6.191-3}
\begin{cases}
\pa_x\bar{u}^{\s} |_{x=0}:=-\pa_y\bar{v}^{\s}\cong \sigma^2 \hat{u}(\sigma y) e^{-\sigma y},\\
\pa_xT^{\s} |_{x=0}:= \f{1}{2\bar{u}^{\s}} \big\{\kappa (\mathscr{T}_0^{\s})''+ \mu |(\mathscr{U}_0^{\s})'|^2 - 2\bar{v}^{\s} (\mathscr{T}_0^{\s})' \big\} \cong \sigma^2 e^{-\sigma y}.
\end{cases}
\end{align}
We can also inductively define
\begin{align} 
\pa_x^k\big(\f{\bar{v}^{\s}}{\bar{u}^{\s}}\big)_y\big|_{x=0},\quad  \pa_x^k\bar{u}^{\s}\big|_{x=0},\quad  \pa_x^k\bar{v}^{\s}\big|_{x=0} \quad \mbox{and}\quad \pa_x^kT^{\s}\big|_{x=0} \quad \mbox{for}\,\,\, k\geq 1.
\end{align}
Therefore, for small $\sigma>0$, we can obtain the smallness on derivatives.
\end{remark}

\medskip

\subsection{The $\theta$-approximate problem}\label{sec2.1}
To solve \eqref{3.58}, we consider the following $\theta$-approximate boundary value problem (BVP)
\begin{align}\label{6.1}
	\begin{cases}
		\dis  (\f{1}{T}u)_x + (\f{1}{T}v)_y =0,\\[3.5mm]
		\dis \big[(\theta + \f{1}{T} u)\pa_x u + \f{1}{T}  v\pa_y u\big] = \mu \pa_{yy} u, \\[3.5mm]
		\dis 2\big[(\theta + \f{1}{T} u)\pa_xT + \f{1}{T}  v\pa_y T\big]   = \kappa \pa_{yy}T + \mu |\pa_y u|^2,
	\end{cases}
	(x,y)\in (0,L)\times \mathbb{R}_+,
\end{align}
or the following equivalent system
\begin{align}\label{5.3}
	\begin{cases}
		\dis \bar{u}_x + \bar{v}_y =0,\\ 
		\dis \bar{u}_{\theta}\pa_x\bar{u} + \bar{v}\pa_y  \bar{u}  = \mu \pa_{yy} \bar{u} + \mathcal{N}_1 ,\\
		\dis 2\big[\bar{u}_{\theta}\pa_xT + \bar{v}\pa_y T\big]   = \kappa \pa_{yy}T +  \mathcal{N}_2,
	\end{cases}
\end{align}
where $\theta>0$ is a parameter.

For fix $\theta\in (0,1]$, we shall show the local existence of solution of \eqref{5.3}. Let $\mathcal{N}_1, \mathcal{N}_2$ be given functions, we can apply similar arguments as in  Oleinik-Samokhin \cite[Lemma 2.1.7]{Oleinik-1999} and also energy estimate to derive the existence of  smooth solution $(\bar{u}^{\t}, \bar{v}^{\t})$ for the $\theta$-approximate quasilinear problem $\eqref{5.3}_{1,2}$ in $x\in [0,x_{\theta}]$ where $x_\theta>0$ is a small constant depending on $\theta$. Then with above $(\bar{u}^{\t}, \bar{v}^{\t})$ as coefficients of $\eqref{5.3}_3$, by applying the standard linear parabolic theory, we can obtain the smooth solution $T^{\t}$ of $\eqref{5.3}_3$ in $x\in [0,x_{\theta}]$ (we may further take $x_{\theta}$  smaller if needed). Finally, we consider an iteration for $\mathcal{N}_1, \mathcal{N}_2$ to establish the local smooth solution $(\bar{u}^{\t}, \bar{v}^{\t}, T^{\t})$ for the original nonlinear problem \eqref{5.3} in $x\in [0,x_{\theta}]$ for some small $x_{\theta}>0$. We also point that the estimates (boundedness and decay in $y\in \mathbb{R}_+$) on $(\bar{u}^{\t},\bar{v}^{\t}, T^{\t})$ may depend on $\theta>0$. We remark that 
\begin{align}
	\bar{v}^{\t}(x,y) = \int_0^y \bar{v}^{\t}_y(x,s) ds = -  \int_0^y \bar{u}^{\t}_x(x,s) ds, 
\end{align}
which yields that $\bar{v}^{\t}(x,0)=0$ and $\bar{v}^{\t}(x,\infty)\neq 0$.

Applying the standard maximum principle to  $\eqref{6.1}_2$, one has that 
	\begin{align}\label{3.70}
		0\leq u^{\t}(x,y)\leq u^0_e(0),\quad \mbox{for}\,\, [0,x_{\theta}]\times \mathbb{R}_+.
	\end{align}
It is clear that
\begin{align*}%
\begin{split}
u^{\t}(x,y)&=\mathscr{U}_0(y) + \int_0^x u^{\t}_x(s,y) ds = \mathscr{U}_0(y) + O(1)x^{\f12} \min\{y,1\} \|u^{\t}_{xyy}\la y\ra^2\|,\\
T^{\t}(x,y)&=\mathscr{T}_0(y) + \int_0^x T^{\t}_x(s,y) ds = \mathscr{T}_0(y) + O(1)x^{\f12} \|T^{\t}_{xy}\la y\ra\|,.
\end{split}
\end{align*}
which, together with the smallness of $x_{\theta}>0$, yields that
\begin{align}\label{6.11}
\begin{split}
\f12 \mathscr{U}_0(y) &\leq u^{\t}(x,y)\leq \f32 \mathscr{U}_0(y),  \\
\f14\leq \f12 \mathscr{T}_0(y) &\leq T^{\t}(x,y)\leq \f32 \mathscr{T}_0(y)\leq 3,  
\end{split}
\quad (x,y)\in [0,x_{\theta}]\times \mathbb{R}_+.
\end{align}

Noting $\mathscr{U}'_0(0)>0$, there exist $y_0>0$ such that $\mathscr{U}'_0(y)\geq \f12 \mathscr{U}'_0(0)>0$, then we have 
\begin{align*}
u^{\t}_y(x,y)&=\mathscr{U}_{0y}(y) + \int_0^x u^{\t}_{xy}(s,y) ds = \mathscr{U}_{0y}(y) + O(1)x^{\f12} \|u^{\t}_{xyy}\la y\ra\|,
\end{align*}
which, with the smallness of $x_\theta$,   implies that
\begin{align}\label{6.11-1}
0<\f12 \mathscr{U}_{0y}(y) &\leq u^{\t}_y(x,y)\leq \f32 \mathscr{U}_{0y}(y),   \,\, \mbox{for}\,\,\, y\in [0,y_0],
\end{align}

\smallskip


\subsection{Reformulation and compatibility conditions}
To further establish the solution of \eqref{3.58} in the limit $\theta \to 0+$, it is crucial to derive uniform-in-$\theta$ estimates for the solutions $(\bar{u}^{\t}, \bar{v}^{\t}, T^{\t})$ of \eqref{5.3}. This constitutes the main contribution of our present paper. For simplicity of presentation, we shall drop the superscript $\theta$ from $(\bar{u}^{\t}, \bar{v}^{\t}, T^{\t})$ in the remainder of this section when no confusion arises.

\subsubsection{Quotient formulation on velocity}\label{Sec4.1}
We define
\begin{align}\label{4.10}
	\mathfrak{u}:=\pa_x \bar{u},  
	\quad \mathfrak{v}:=\pa_x \bar{v}. 
\end{align}
Applying $\pa_x$ to $\eqref{5.3}_{1,2}$, one obtains
\begin{align}\label{6.5-1}
	\begin{cases}
		\dis \mathfrak{u}_x + \mathfrak{v}_y =0,\\
		\dis \bar{u}_{\theta}\pa_x\mathfrak{u} + \bar{u}_y \mathfrak{v}  + \bar{v}\pa_y \mathfrak{u} + \mathfrak{u}\, \bar{u}_x  - \mu \pa_{yy} \mathfrak{u}\\
		\dis = \f{2\mu}{T} T_{y}\, \mathfrak{u}_y + \big[\mu-\f{\kappa}{2}\big]\f{T_{yy}}{T} \mathfrak{u} - \f{3\mu}{2T}  T_y^2\bar{u}^2 \mathfrak{u} - \mu T_y [\mathfrak{u}_y\,\bar{u}^2 + 2\bar{u}\bar{u}_y \mathfrak{u}] \\
		\dis\qquad   - \f{1}{2}\mu T [2\bar{u}_y\bar{u} \mathfrak{u}_y + \bar{u}_y^2\mathfrak{u}] + f_1,\\
		(\mathfrak{u}, \mathfrak{v})|_{y=0}=(0,0).
	\end{cases}
\end{align}
where 
\begin{align}\label{6.6}
f_1:= \big[\mu-\f{1}{2}\kappa\big]\bar{u} \pa_x\big(\f{T_{yy}}{T}\big) + \bar{u}_y\,\pa_x\big(\f{2\mu}{T} T_{y}\big)  - \bar{u}^3 \pa_x\big(\f{\mu}{2T}  T_y^2\big) - \mu  \bar{u}_y\bar{u}^2 \pa_x T_y - \f{1}{2}\mu  \bar{u}^2_y\bar{u}\, \pa_x T.
\end{align}

We have from $\eqref{5.3}_1$ and \eqref{4.10} that 
\begin{align*}
\mathfrak{u} \equiv \bar{u}_x =- \bar{v}_y,
\end{align*}
Motivated by  \cite{Guo-Iyer-CMP}, we define the good unknown, i.e., the {\it quotient}:
\begin{align}\label{6.14-0}
	q:= \f{\bar{v}}{\bar{u}_{\theta}}\,\, \Longleftrightarrow \,\,  \bar{v} = \bar{u}_{\theta} \, q.
\end{align}

\smallskip

A direct calculation shows that 
\begin{align*}
\bar{u}_{\theta}\pa_x\mathfrak{u} + \bar{u}_y \mathfrak{v} 
&=-\bar{u}_{\theta} \big(\bar{u}_{\theta}  q_{xy} + \bar{u}_{x}  q_y + \cancel{\bar{u}_{y}  q_x} + \bar{u}_{xy}  q \big)  +  \bar{u}_{y} \big( \cancel{\bar{u}_{\theta} q_x} + \bar{u}_{x} q \big) \nonumber\\
&=-\bar{u}_{\theta}^2  q_{xy} -\bar{u}_{\theta} \bar{u}_x q_y - \big[\bar{u}_{\theta} \bar{u}_{xy}    - \bar{u}_y \bar{u}_{x}\big] q,
\end{align*}
which, together with $\eqref{6.5-1}_2$, yields that 
\begin{align}\label{6.17}
&\bar{u}_{\theta}^2  q_{xy} - \mu  \pa_{yyy} [\bar{u}_{\theta}q]  
\equiv \bar{u}_{\theta}^2  q_{xy} - \mu \bar{v}_{yyy}
= F,
\end{align}
where $F=-f_1 + f_2 + f_3$ with 
\begin{align}\label{6.17-0}
\begin{cases}
\dis f_2:= 2\bar{u}_{\theta} \bar{v}_y q_y ,\\
\dis f_3:  =  \Big\{\big[\mu-\f{\kappa}{2}\big]\f{T_{yy}}{T}  - \f{3\mu}{2T}  T_y^2\bar{u}^2   - 2 \mu T_y \bar{u}\bar{u}_y  - \f{1}{2}\mu T\bar{u}_y^2  \Big\}\,  \bar{v}_y \\
\dis\qquad\,\, + \Big\{\f{2\mu}{T} T_{y} - \mu T_y \bar{u}^2  - \mu T \bar{u}_y\bar{u} \Big\}\, \bar{v}_{yy}.
\end{cases}
\end{align}
It is clear to know that 
\begin{align}\label{6.17-1}
q\big|_{y=0}=0\quad\mbox{and}\quad  q_y\big|_{y=0} = \big[\f{\bar{v}_y}{\bar{u}_{\theta}}  -  \f{\bar{v}}{\bar{u}_{\theta}^2 }\bar{u}_{y} \big]\Big|_{y=0}=0.
\end{align}



\subsubsection{Formulation on temperature} Recall the cut off function $\chi$ in \eqref{cutoff}. 
Define the auxiliary functions
\begin{align}\label{6.195-1}
	T_{B}:=
	\begin{cases}
		1, \quad & \mbox{NBC},\\
		T_b(x)  \chi(y) +  [1-\chi(y)],\,& \mbox{DBC},
	\end{cases}
\end{align}
and
\begin{align}\label{6.195-2}
	G_{B}:=
	\begin{cases}
		0, \quad & \mbox{NBC},\\
		\kappa \pa_{yy}T_B - 2 [\bar{u}\pa_x T_B + \bar{v}\pa_y T_B],\,& \mbox{DBC}.
	\end{cases}
\end{align}

Define
\begin{align*}
	\Theta&:=T(x,y) - T_B(x,y),
\end{align*}
then the temperature equation $\eqref{5.3}_3$ can be rewritten as
\begin{align}\label{6.195-3}
	2 \bar{u}_{\theta}\pa_x\Theta + 2\bar{v}\pa_y \Theta  = \kappa \pa_{yy}\Theta + \mu |\pa_y ([T_B+\Theta]\bar{u})|^2 + G_B.
\end{align}
The initial and boundary conditions become
\begin{align}\label{6.195-30}
\begin{cases}
\dis \Theta\big|_{x=0}= \mathscr{T}_0 - T_B, \quad \lim_{y\to+\infty}\Theta(x,y)=0,\\
\dis	\Theta_y\big|_{y=0}=0 \,\, \mbox{for NBC}\,\, \big(\mbox{or}\,\, \Theta\big|_{y=0}=0 \,\, \mbox{for DBC} \big).
\end{cases}
\end{align}

\subsubsection{Functional space} 
For later use, we denote the weight function
\begin{align}
w(y):=\la y\ra^{\mathfrak{l}} \quad \mbox{for}\,\,\, \mathfrak{l}\gg1.
\end{align}
We define the following norms for  quotient:
\begin{align}
	\begin{split}
		\|q_y\|^2_{\mathbf{X}_{\theta, w}}:&=  \|\sqrt{\bar{u}_{\theta}} q_{yy}w\|^2 + \sup_{x_0\in[0,L]}\|\bar{u}_{\theta} q_yw\|^2_{x=x_0},\\
		\|q_y \|^2_{\mathbf{X}_{\theta, w}^m}:&= \sum_{i=0}^m  \|\pa_x^iq_y \|^2_{\mathbf{X}_{\theta, w}}
		\equiv \sum_{i=0}^m \Big\{\|\sqrt{\bar{u}_{\theta}} \pa_{x}^i q_{yy}w\|^2 + \sup_{x_0\in[0,L]}\|\bar{u}_{\theta}\pa_{x}^i q_yw\|^2_{x=x_0}\Big\}.
	\end{split}
\end{align}

For temperature, we define norms
\begin{align}
	\begin{split}
		\|h \|^2_{\MY_{\t,w}} &:= \kappa \|h_{y}w\|^2 + \sup_{x_0\in[0,L]}\|\sqrt{\bar{u}_{\t}}h w\|^2_{x=x_0},\\
		\|h\|^2_{\MY^m_{\t, w}}&:=\sum_{i=0}^m \|\pa_x^ih\|^2_{\MY_{\t, m}}
		=\sum_{i=0}^k \Big\{\kappa \|\pa_x^ih_{y}w\|^2 + \sup_{x_0\in[0,L]}\|\sqrt{\bar{u}_{\t}}\pa_x^ih w\|^2_{x=x_0}\Big\}.
	\end{split}
\end{align}
For $\theta=0$, we denote $\mathbf{X}_{w}:=\mathbf{X}_{0,w}$,  $\mathbf{X}^m_w:=\mathbf{X}^m_{0,w}$ and  $\mathbf{Y}_w:=\mathbf{Y}_{0,w}$,  $\mathbf{Y}^m_w:=\mathbf{Y}^m_{0,w}$.

For later use, we define the following notations
	\begin{align}
		\||\bar{v}|\|_{m,w}:&= 1+ \|q_y\|^2_{\MX^m_{\t, w}}  + \sum_{i=0}^{m-1}\|\pa_x^{i}\bar{v}_{yy}w\|^2_{x=x_0} + 	\sum_{i=0}^m\|\pa_x^i\bar{v}_yw\|  
		  + \sum_{i=0}^{m-1} \|\pa^i_xq\|_{L^\infty} \nonumber\\
		&\quad + \sum_{i=1}^{m-1} \|\pa^i_x\bar{u}\la y\ra^{l_0}\|_{L^\infty} + \sum_{i=0}^{m-1} \|\pa^i_x\bar{u}_y\la y\ra^{l_0}\|_{L^\infty} + \sum_{i=0}^{1+\f{m}{2}} \|\pa^i_x\bar{u}_{yy}\la y\ra^{l_0}\|_{L^\infty},\\
		[[[\Theta]]]_{m+1,w}:&=\|\Theta\|^2_{\MY^{m}_{\t, w}} + \|S_{m+1}\|^2_{\MY_{\t,w}} +   \|(\Theta_y,  \Theta_{yy},  \Theta_{yyy}, \Theta_{xy}, \Theta_{xxy}) w\|_{L^\infty_xL^2_y}\nonumber\\
		&\quad  
		+ \sum_{i=0}^{1+\f{m}{2}} \|(\pa_x^i\Theta,\pa_x^i\Theta_y)\la y\ra^{l_0}\|_{L^\infty}  + \sum_{i=0}^{\f{m}{2}}\|\pa_x^i\Theta_{yy}\la y\ra^{l_0}\|_{L^\infty},
	\end{align}
where $0<l_0\leq \f12 \mathfrak{l}$ and $S_{m+1}$ 
is the pseudo entropy defined in \eqref{6.231}.

\smallskip

\subsection{Estimates on quotient} In this subsection, we aim to derive some uniform estimates on quotient. 

\subsubsection{Some useful inequalities}
\begin{lemma}\label{lemA.5-1}
	For $i\geq 0$, it holds that 
	\begin{align}\label{6.54}
			|\pa_x^i \bar{u}(x,y)| 
			&\lesssim 
			\begin{cases}
				\min\{\sqrt{y},1\} \|\pa_x^i\bar{u}_{y} \la y\ra\|_{L^2_y},\\
				\min\{y,1\} \|\pa_x^i\bar{u}_{yy} \la y\ra^2\|_{L^2_y},
			\end{cases}
	\end{align}
	and
	\begin{align}\label{6.55} 
		|\pa_x^i\bar{v}(x,y)| 
		&\lesssim 
		\begin{cases}
			\min\{\sqrt{y},1\} \|\pa^{i}_x\bar{v}_y \la y\ra\|_{L^2_y},\\
			\min\{y^{\f32},1\} \| \pa^{i}_x\bar{v}_{yy} \la y\ra^2\|_{L^2_y},\\
			\min\{y^2, 1\} \|\pa^{i}_x\bar{v}_{yyy}\la y\ra^3\|_{L^2_y}.
		\end{cases}
	\end{align}
\end{lemma}

\noindent{\bf Proof.}
For $0\leq i\leq m$, it holds that  
\begin{align}\label{6.165-2}
	\begin{split}
		|\pa_x^i \bar{u}(x,y)|&\lesssim \int_0^y |\pa_x^{i}\bar{u}_y(x,z)| dz \lesssim \min\{y^{\f12},1\} \|\pa_x^i\bar{u}_{y} \la y\ra\|_{L^2_y},\\
		|\pa_x^i \bar{u}(x,y)|&\lesssim \int_0^y |\pa_x^{i}\bar{u}_y(x,z)| dz \lesssim  \int_0^y \Big|\int_z^\infty \pa_x^{i}\bar{u}_{yy}(x,s) ds\Big| dz \\
		& \lesssim \min\{y,1\} \|\pa_x^i\bar{u}_{yy} \la y\ra^2\|_{L^2_y}.
	\end{split}
\end{align}

For $\pa_x^i\bar{v} $, we note
\begin{align}\label{6.165-0} 
	|\pa_x^i\bar{v}(x,y)| &\lesssim \int_0^y |\pa^{i}_x\bar{v}_y(x,s)| ds \lesssim \min\{\sqrt{y},1\} \|\pa^{i}_x\bar{v}_y \la y\ra\|_{L^2_y},
\end{align}
and
\begin{align}\label{6.165-1}
	|\pa_x^i\bar{v}(x,y)| &\lesssim \int_0^y |\pa^{i}_x\bar{v}_y(x,s)| ds {\bf 1}_{\{y\in[0,1]\}}  +  \|\pa^{i}_x\bar{v}_y\la y\ra\|_{L^2_y} {\bf 1}_{\{y\geq 1\}}\nonumber\\
	&\lesssim  \int_0^y\int_0^s |\pa^{i}_x\bar{v}_{yy}(x,t)| dt ds {\bf 1}_{\{y\in[0,1]\}}  +  \|\pa^{i}_x\bar{v}_y  \la y\ra\|_{L^2_y} {\bf 1}_{\{y\geq 1\}} \nonumber\\
	&\lesssim 
	\begin{cases}
		\min\{y^{\f32},1\} \| \pa^{i}_x\bar{v}_{yy} \la y\ra^2\|_{L^2_y},\\
		\min\{y^2, 1\} \|\pa^{i}_x\bar{v}_{yyy}\la y\ra^3\|_{L^2_y}, 
	\end{cases}
\end{align}
where we have used the Hardy's inequality 
\begin{align*}
	\|\pa^{i}_x\bar{v}_y  \la y\ra\|_{L^2_y}\lesssim \|\pa^{i}_x\bar{v}_{yy}  \la y\ra^2 \|_{L^2_y} \lesssim \|\pa^{i}_x\bar{v}_{yyy}  \la y\ra^3\|_{L^2_y}.
\end{align*}
Therefore the proof of Lemma \ref{lemA.5-1} is completed. $\hfill\Box$

\medskip

\begin{lemma}\label{lemA.5-2}
For any $l\geq 0$, it holds that 
\begin{align}\label{6.29}
\begin{split}
|\bar{u}_y(x,y)\la y\ra^{l}| 
&\lesssim  \|\bar{u}_{yy}\la y\ra^{l+1}\|_{x=0} + \sqrt{L} \|\bar{v}_{yyy}\la y\ra^{l+1}\|,  \\
|\bar{u}_{yy}(x,y)\la y\ra^{l}| 
&\lesssim  \|\bar{u}_{yyy}\la y\ra^{l+1}\|_{x=0} + \sqrt{L} \|\bar{v}_{yyyy}\la y\ra^{l+1}\|,
\end{split}
\end{align}
and
\begin{align}\label{6.30}
\begin{split}
\|\bar{u}_y \la y\ra^{l}\|^2_{x=x_0} 
&\leq \|\bar{u}_y \la y\ra^{l}\|^2_{x=0} + L \|\bar{v}_{yy} \la y\ra^{l}\|^2,\\
\|\bar{u}_{yy} \la y\ra^{l}\|^2_{x=x_0} 
&\leq \|\bar{u}_{yy} \la y\ra^{l}\|^2_{x=0} + L \|\bar{v}_{yyy} \la y\ra^{l}\|^2.
\end{split}
\end{align}
For $k\geq 1$, we have 
\begin{align} \label{A6.109}
\begin{split}
|\pa_x^{k}\bar{u}\la y\ra^{l}| &\lesssim  \|\pa_x^{k-1}\bar{v}_{yy}\la y\ra^{l+1}\|_{L^2_y},\\
|\pa_x^{k}\bar{u}_y\la y\ra^{l}| &\lesssim \|\pa_x^{k-1}\bar{v}_{yyy}\la y\ra^{l+1}\|_{L^2_y},\\
\|\pa_x^{k} \bar{u}_{yy}\la y\ra^{l}\|_{L^2_xL^\infty_y} &\lesssim  \|\pa_x^{k-1}\bar{v}_{yyyy} \la y\ra^{l+1}\|.
\end{split}
\end{align}
\end{lemma}

\noindent{\bf Proof.} Let $h\to 0$ as $y\to\infty$, we note that
\begin{align*}
	|h(y)|&\lesssim \int_y^{\infty} |h_y(z_1)| dz_1 \lesssim \Big(\int_y^\infty |h_y(z_1)|^2 \la z_1\ra^{2l+2}  dz_1\Big)^{\f12} \Big(\int_y^\infty \la z_1\ra^{-2l-2} dz_1\Big)^{\f12} \nonumber\\
	&\lesssim \la y\ra^{-l-\f12} \Big(\int_y^\infty |h_y(z_1)|^2 \la z_1\ra^{2l+2}  dz_1\Big)^{\f12},\quad \mbox{for}\,\,l\geq0,
\end{align*}
which implies that 
\begin{align}\label{6.162}
		|\la y \ra^{l+\f12} h(y)|\lesssim  \|\la y\ra^{l+1} h_y\|_{L^2_y},\quad \mbox{for}\,\,l\geq 0.
\end{align}
Also it follows from Sobolev inequality that 
\begin{align*}
|h(y)\la y\ra^{l}|^2 \lesssim \|h(y)\la y\ra^{l}\|^2 + \|h\la y\ra^{l}\|_{L^2_y}\cdot \|h_y\la y\ra^{l}\|_{L^2_y}
\end{align*}

Using \eqref{6.162}, we have 
\begin{align*}
	|\bar{u}_y(x,y)\la y\ra^{l}|&\lesssim \|\bar{u}_{yy}(x,y)\la y\ra^{l+1}\|_{L^2_y} \lesssim \|\bar{u}_{yy}\la y\ra^{l+1}\|_{x=0} + \sqrt{L}\|\bar{u}_{xyy}\la y\ra^{l+1}\| \nonumber\\ 
	&\lesssim \|\bar{u}_{yy}\la y\ra^{l+1}\|_{x=0} + \sqrt{L} \|\bar{v}_{yyy}\la y\ra^{l+1}\|,  \\
	|\bar{u}_{yy}(x,y)\la y\ra^{l}| 
	&\lesssim  \|\bar{u}_{yyy}(x,y)\la y\ra^{l+1}\|_{L^2_y}  \lesssim   \|\bar{u}_{yyy}\la y\ra^{l+1}\|_{x=0} + \sqrt{L}  \|\bar{u}_{xyyy}(x,y)\la y\ra^{l+1}\|  \nonumber\\
	&\lesssim  \|\bar{u}_{yyy}\la y\ra^{l+1}\|_{x=0} + \sqrt{L} \|\bar{v}_{yyyy}\la y\ra^{l+1}\|.
\end{align*}
Also it holds that 
\begin{align*}
\begin{split}
\|\bar{u}_y \la y\ra^{l}\|^2_{x=x_0}&\lesssim \|\bar{u}_y\la y\ra^{l}\|^2_{x=0} +  L \|\bar{u}_{xy}\la y\ra^{l}\|^2 \lesssim \|\bar{u}_y \la y\ra^{l}\|_{x=0} + L \|\bar{v}_{yy} \la y\ra^{l}\|^2,\\
\|\bar{u}_{yy} \la y\ra^{l}\|^2_{x=x_0}&\lesssim \|\bar{u}_{yy} \la y\ra^{l}\|^2_{x=0} +   L \|\bar{u}_{xyy}\la y\ra^{l}\|^2 \leq \|\bar{u}_{yy} \la y\ra^{l}\|_{x=0} + L \|\bar{v}_{yyy} \la y\ra^{l}\|^2,
\end{split}
\end{align*}

\smallskip

For $k\geq 1$, we have from \eqref{6.162} that 
\begin{align*}
	\begin{split}
		|\pa_x^{k}\bar{u}\la y\ra^{l}| &\lesssim  \|\pa_x^{k} \bar{u}_{y} \la y\ra^{l+1}\|_{L^2_y} \lesssim \|\pa_x^{k-1}\bar{v}_{yy}\la y\ra^{l+1}\|_{L^2_y},\\
		|\pa_x^{k}\bar{u}_y\la y\ra^{l}| &\lesssim  \|\pa_x^{k} \bar{u}_{yy} \la y\ra^{l+1}\|_{L^2_y} \lesssim \|\pa_x^{k-1}\bar{v}_{yyy}\la y\ra^{l+1}\|_{L^2_y},
	\end{split}
\end{align*}
and
\begin{align*}
	\|\pa_x^{k} \bar{u}_{yy}\la y\ra^{l}\|_{L^2_xL^\infty_y}&\lesssim \|\pa_x^{k} \bar{u}_{yyy} \la y\ra^{l+1}\| \lesssim  \|\pa_x^{k-1}\bar{v}_{yyyy} \la y\ra^{l+1}\|.
\end{align*}
Therefore the proof of Lemma \ref{lemA.5-2} is completed. $\hfill\Box$

\medskip

\begin{lemma}[Hardy-type Inequality]\label{lem6.1}
We denote $a:=\mathscr{U}'_0(0)>0$. It holds that 
	\begin{align}\label{6.20}
		\begin{split}
			\|gw\|^2&\leq  \f{C\xi}{a} \|\sqrt{\bar{u}}g_y \chi(\f{y}{\xi})\|^2 + \f{CL}{a^2\xi^2}  \sup_{x\in[0,x_0]} \|\bar{u} g w\|^2_{L^2_y},
		\end{split}
	\end{align}
	where the constant $C>0$ is independent of $\theta$.
\end{lemma}

\noindent{\bf Proof.} It is noted that 
\begin{align}\label{6.20-1}
	\|g w\|^2& \leq  2\|g\chi(\f{y}{\xi})  \|^2 + 2\|g w [1-\chi(\f{y}{\xi})]\|^2.
\end{align}
It is clear that 
\begin{align}\label{6.21}
	\|gw [1-\chi(\f{y}{\xi})]\|^2\lesssim \f{1}{a^2\xi^2} \| \bar{u} gw [1-\chi(\f{y}{\xi})]\|^2 \lesssim \f{L}{a^2\xi^2}  \sup_{x\in[0,x_0]} \|\bar{u}g w\|^2_{L^2_y}.
\end{align}
Integrating by parts in $y$, one has that 
\begin{align}\label{6.22}
	\|g\chi(\f{y}{\xi}) \|^2&=\iint g^2 \chi^2(\f{y}{\xi}) dydx =-\iint \Big[2y g_y g \chi^2(\f{y}{\xi}) + \f{2}{\xi} y g^2 \chi(\f{y}{\xi}) \chi'(\f{y}{\xi})\Big] dydx\nonumber\\
	&\lesssim \|yg_y\chi(\f{y}{\xi})\|\cdot \|g\chi(\f{y}{\xi})\| + \f{L}{a^2\xi^2}  \sup_{x\in[0,x_0]} \|\bar{u} g w\|^2_{L^2_y}\nonumber\\
	&\leq \f{1}{2} \|g\chi(\f{y}{\xi}) \|^2 + C \f{\xi}{a} \|\sqrt{\bar{u}} g_y \chi(\f{y}{\xi})\|^2 + \f{L}{a^2\xi^2}  \sup_{x\in[0,x_0]} \|\bar{u} g w\|^2_{L^2_y}.
\end{align}

Substituting \eqref{6.21}-\eqref{6.22} into \eqref{6.20-1}, one gets that 
\begin{align*}
	\|g w\|^2
	&\lesssim \f{\xi}{a} \|\sqrt{\bar{u}}g_y \chi(\f{y}{\xi})\|^2 + \f{L}{a^2\xi^2}  \sup_{x\in[0,x_0]} \|\bar{u} g w\|^2_{L^2_y},
\end{align*}
which concludes \eqref{6.20}. Therefore the proof of Lemma \ref{lem6.1} is completed. $\hfill\Box$

\medskip

\begin{lemma}\label{lemA.5}
	For $k\geq 0$, the following inequalities hold
	\begin{align}
		\|\pa_x^kq_yw\| + \|\f{\pa_x^kq}{y}\| &\lesssim 
		L^{\f16-} \|\pa_x^kq_y\|_{\mathbf{X}_w},\label{A6.23}\\
		\|\pa_x^k q_yw\|_{L^\infty_x L^2_y}&\lesssim \|\pa_x^k q_y w\|_{x=0} + L^{\f23-} \|\pa_x^{k+1}q_y\|_{\mathbf{X}_w}, \label{A6.62}\\
		\|\pa_x^k q\|_{L^\infty}&\lesssim \|\pa_x^k q_y\la y\ra\|_{x=0} + L^{\f23-} \|\pa_x^{k+1}q_y\|_{\mathbf{X}_{\la y\ra}},\label{A6.61}\\
	\|\pa_x^kq\|_{L^2_xL^\infty_y} 
		&\lesssim  L^{\f16-} \|\pa_x^k q_y\|_{\MX_{\la y\ra}}. \label{6.115}
	\end{align}
We also have 
	\begin{align}
		\|\pa_x^{k}q_{yy}w\|
		&\lesssim P\big(\||\bar{v}|\|_{k,w}\big)\sum_{i=1}^{k}  \|\pa_x^i\bar{v}_{yyy} \la y\ra^2\|  + L^{\f12-} \|q_y\|_{\MX^{k+1}_w}\nonumber\\
		&\quad +  L^{\f12}\|\pa_x^{k}q_{yy}w\bar{\chi}\|_{x=0} , \quad \mbox{for}\,\, k\geq 0,\label{6.38-0}
	\end{align}
	and
	\begin{align}\label{6.38-1}
	\|q_{yy}w\|_{L_x^\infty L^2_y}&\lesssim  \|q_{yy}w\bar{\chi}\|_{x=0} + L^{\f12-} \|\sqrt{\bar{u}}q_{xyy}w\bar{\chi}\| + P\big(\|\bar{u}_{yy}\|_{L^\infty}\big) \|(\bar{v}_{yy}, \bar{v}_{yyy})\|_{x=x_0}.
	\end{align}
Hereafter $P(\cdot)$ is some polynomial function which may vary from line to line.
\end{lemma}

\noindent{\bf Proof.} 1. Taking $\xi=L^{\f13}$ in Lemma \ref{lem6.1}, we obtain
\begin{align}\label{6.48}
	\|\pa_x^kq_yw\|^2&\lesssim L^{\f13-} \Big\{\|\sqrt{\bar{u}}\pa_x^kq_{yy}\|^2 + \sup_{x\in[0,x_0]} \|\bar{u} \pa_x^kq_y w\|^2_{L^2_y}\Big\}
	\lesssim L^{\f13-} \|\pa_x^kq_y\|^2_{\mathbf{X}_w}.
\end{align}
Applying Hardy inequality and \eqref{6.48},  one obtains
\begin{align}\label{6.48-1}
	\|\f{\pa_x^kq}{y}\|^2\lesssim \|\pa_x^kq_y\|^2 \lesssim L^{\f13-} \Big\{\|\sqrt{\bar{u}}\pa_x^kq_{yy}\|^2 +  \sup_{x\in[0,x_0]}\|\bar{u} \pa_x^kq_y\|^2_{L^2_y} \Big\}
	\lesssim L^{\f13-} \|\pa_x^kq_y\|^2_{\mathbf{X}_1}.
\end{align}
Hence we obtain \eqref{A6.23}.


 It is clear that 
\begin{align}\label{A6.77}
	\int_0^\infty |\pa_x^{k}q_y(x,y)|^2w^2 dy 
	&\lesssim \int_0^\infty |\pa_x^{k}q_y(0,y)|^2w^2 dy + L \int_0^{L} \int_0^\infty |\pa_x^{k}q_{xy}(x,y)|^2w^2 dy \nonumber\\
	&\lesssim  \|\pa_x^{k}q_yw\|^2_{x=0} + L\|\pa_x^{k+1}q_{y}w\|^2\nonumber\\
	& \lesssim \|\pa_x^{k}q_y w\|^2_{x=0} + L^{\f43-} \|\pa_x^{k+1}q_y\|^2_{\mathbf{X}_w},
\end{align}
where we have used \eqref{A6.23}. It follows from \eqref{A6.77} that
\begin{align}\label{A6.81}
	|\pa_x^{k}q(x,y)|^2
	&\lesssim \int_0^\infty |\pa_x^{k}q_y(x,y)|^2\la y\ra^2 dy 
	\lesssim \|\pa_x^{k}q_y\la y\ra\|^2_{x=0} + L^{\f43-} \|\pa_{x}^{k+1}q_y \|^2_{\mathbf{X}_{\la y\ra}}.
\end{align}
Also it is clear to have
\begin{align} 
	\|\pa_x^kq\|^2_{L^2_xL^\infty_y} \lesssim \int_0^L\int_0^\infty |\pa_x^{k}q_y(x,y)|^2\la y\ra^2 dy dx \lesssim L^{\f13-} \|\pa_x^k q_y \|^2_{\MX_{\la y\ra}},
\end{align}
where we have used \eqref{6.48} in the last inequality.

2.  Noting
\begin{align*}
	\begin{split}
		&\f{1}{y}\bar{v}(x,y) =\f{1}{y} \int_0^y \bar{v}_y(x,z) dz = \int_0^1 \bar{v}_y(x,ty) dt\\
		&\Longrightarrow \big(\f{1}{y}\bar{v}\big)_{yy} = \int_0^1 t^2  \bar{v}_{yyy}(x,ty) dt,
	\end{split}
\end{align*}
which implies that 
\begin{align}
\begin{split}
\int_0^1 \big|\big(\f{\bar{v}}{y}\big)_{y}\big|^2 dy 
&\leq  \int_0^1 \int_0^1 t^2 |\bar{v}_{yy}(x,ty)|^2 dt dy 
\lesssim \int_0^1 |\bar{v}_{yy}(x,y)|^2 dy,\\
\int_0^1 \big|\big(\f{\bar{v}}{y}\big)_{yy}\big|^2 dy 
&\leq  \int_0^1 \int_0^1 t^4 |\bar{v}_{yyy}(x,ty)|^2 dt dy  
\lesssim  \int_0^1 |\bar{v}_{yyy}(x,y)|^2 dy.
\end{split}
\end{align}

It follows from $\bar{v}_y(x,0)=\bar{v}(x,0)=0$ and Hardy's inequality that 
\begin{align}
	\begin{split}
		\|\f{\bar{v}}{y}\|_{L^2_y}  \lesssim \|\bar{v}_{y}\|_{L^2_y} \quad \mbox{and}\quad 
		\|\f{\bar{v}}{y^2}\|_{L^2_y}  \lesssim \|\f{\bar{v}_{y}}{y}\|_{L^2_y} \lesssim \|\bar{v}_{yy}\|_{L^2_y}.
	\end{split}
\end{align}

For $y\in [0,1]$, we note that 
\begin{align}\label{A6.64}
	q&=\f{\bar{v}}{\bar{u}_{\theta}} = \f{1}{\bar{u}_{y}(x,0)}\f{\bar{v}}{y} + \Big(\f{1}{\bar{u}_{\theta}} - \f{1}{y \bar{u}_{y}(x,0)} \Big) \bar{v} \nonumber\\
	&= \f{1}{\bar{u}_{y}(x,0)}\f{\bar{v}}{y} - \f{\theta +  r_1(x,y)}{ (\theta+ \bar{u}) \bar{u}_{y}(x,0)} \f{\bar{v}}{y}
\end{align}
where we have used the following Taylor expansion 
\begin{align}
	\bar{u}(x,y)
	= y \bar{u}_y(x,0) +  r_1(x,y), 
\end{align}
where $\dis r_1:=\int_0^y (\pa_y^2\bar{u})(x,t) \, (y-t) dt \equiv y^2 \int_0^1 (1-s) (\pa_y^2\bar{u})(x,sy) ds.$

Using \eqref{A6.64}, we have 
\begin{align}\label{A6.65}
	q_{yy}
	&=\f{1}{\bar{u}_{y}(x,0)} \Big[1 - \f{\theta + r_1}{\theta+ \bar{u}} \Big] \big(\f{\bar{v}}{y}\big)_{yy} - \f{2}{\bar{u}_{y}(x,0)} \Big(\f{\theta + r_1}{\theta+ \bar{u}}\Big)_y \big(\f{\bar{v}}{y}\big)_{y} \nonumber\\
	&\quad - \f{1}{\bar{u}_{y}(x,0)}\Big(\f{\theta + r_1}{\theta+ \bar{u}}\Big)_{yy} \f{\bar{v}}{y}=:I_1 + I_2 + I_3.
\end{align}
Then, for $I_1, I_2$, it is direct to get that 
\begin{align}\label{A6.68-0}
	\begin{split}
		\int_0^{x_0} \int_0^1|I_1|^2 dydx 
		&\lesssim \int_0^{x_0} \int_0^1 \big|\big(\f{\bar{v}}{y}\big)_{yy}\big|^2 dxdy 
		\lesssim \|\bar{v}_{yyy}\|^2_{L^2_{loc}},\\
		\int_0^{x_0} \int_0^1|I_2|^2 dydx &\lesssim \int_0^{x_0} \int_0^1 \big|\big(\f{\bar{v}}{y}\big)_{y}\big|^2 dxdy  \lesssim  \|\bar{v}_{yy}\|^2_{L^2_{loc}},\\
		\int_0^{x_0} \int_0^1|I_3|^2 dydx & \lesssim  \int_0^{x_0} \int_0^1 \big|\f{\bar{v}}{y^2}\big|^2 dydx \lesssim  \|\bar{v}_{yy}\|^2,
	\end{split}
\end{align}
where we have used the Hardy's inequality in the last line. Then it follows from \eqref{A6.65} and \eqref{A6.68-0} that 
\begin{align}\label{A6.68}
	\int_0^{x_0}\int_0^1 |q_{yy}|^2 dydx \lesssim  \|\bar{v}_{yyy}\|^2_{L^2_{loc}} + \|\bar{v}_{yy}\|^2.
\end{align}
For $y\geq 1$, it holds that
\begin{align}\label{A6.69}
	\int_0^{x_0}\int_1^\infty |q_{yy}w|^2 dydx 
	&\lesssim \int_0^{x_0}\int_1^\infty |q_{yy}(0,y)w|^2 dydx + \int_0^{x_0}\int_1^\infty \Big|\int_0^{x}q_{xyy}(s,y)ds\Big|^2 w^2 dydx\nonumber\\
	&\lesssim L\|q_{yy}w\mathbf{1}_{\{y\geq 1\}}\|^2_{x=0} + L^{2-} \|\sqrt{\bar{u}}q_{xyy}w\|^2.
\end{align}
Thus it follows from \eqref{A6.68}-\eqref{A6.69} that
\begin{align}
	\|q_{yy}w\|^2 \lesssim L\|q_{yy}w\bar{\chi}\|^2_{x=0} + L \|\sqrt{\bar{u}}q_{xyy}w\|^2
	+ \|\bar{v}_{yyy}\|^2_{L^2_{loc}} + \|\bar{v}_{yy}\|^2.
\end{align}

Next, we consider the estimate on $\|\pa_x^{k}q_{yy}w\|$ for $k\geq 1$. Similar as \eqref{A6.69}, it holds that 
\begin{align}\label{A6.87}
	\int_0^{x_0}\int_1^\infty |\pa_x^k q_{yy}w|^2 dydx 
	&\lesssim L\|\pa_x^kq_{yy}w\mathbf{1}_{\{y\geq 1\}}\|^2_{x=0} +  L^{1-} \|\sqrt{\bar{u}}\pa_x^{k+1}q_{yy}w\|^2.
\end{align}
It follows from  \eqref{A6.65} that 
\begin{align}\label{A6.88}
	\|\pa_x^kq_{yy}\chi\|^2 & \lesssim P\Big(\sum_{i=1}^{\f{k}{2}}\|\pa_x^i\bar{u} \la y\ra^{l_0}\|_{L^\infty}, \sum_{i=0}^{\f{k}{2}}\|\pa_x^i\bar{u}_{yy} \la y\ra^{l_0}\|_{L^\infty}\Big) \nonumber\\
	&\quad \times \sum_{i=\f{k}{2}}^{k} \big\{\|(\pa_x^{i}\bar{v}_{y},\pa_x^{i}\bar{v}_{yy}, \pa_x^i\bar{v}_{yyy})\|^2  + \|\pa_x^i\bar{u}_{yy}\la y\ra\|^2 \big\}.
\end{align}
Hence we have from \eqref{A6.87}-\eqref{A6.88} that 
\begin{align}
\|\pa_x^{k}q_{yy}w\|^2 
& \lesssim P\big(\||\bar{v}|\|_{k,w}\big)  \sum_{i=1}^{k}  \|\pa_x^i\bar{v}_{yyy} \la y\ra^2\|^2 + L^{1-} \|\pa_x^{k+1}q_y\|^2_{\MX_w} + L\|\pa_x^{k}q_{yy}w\bar{\chi}\|^2_{x=0}.
\end{align}

Similarly, it is direct to obtain
\begin{align}
\|q_{yy}w\|_{L_x^\infty L^2_y}&\lesssim  \|q_{yy}w\bar{\chi}\|_{x=0} + L^{\f12-} \|\sqrt{\bar{u}}q_{xyy}w\bar{\chi}\| + P\big(\|\bar{u}_{yy}\|\big) \|(\bar{v}_{yy}, \bar{v}_{yyy})\|_{x=x_0}.
\end{align}
Therefore the proof of Lemma \ref{lemA.5} is completed. $\hfill\Box$

\medskip

\begin{lemma}\label{lemA.6}
For $m\geq 2$, it holds that 
\begin{align}\label{6.140-2}
\begin{split}
\|\pa_x^m\bar{u} \la y\ra^l\|_{L^\infty} + \|(\pa_x^m\bar{u} \la y\ra^{l}, \, \pa_x^m\bar{u}_y \la y\ra^{l+1})\|_{L^\infty_xL^2_y} 
&\lesssim\||\bar{v}|\|_{m,w},\\
\|\pa_x^m\bar{u}_y \la y\ra^{l}\|_{L^2_xL^\infty_y}&\lesssim   \|\pa_x^{m-1}\bar{v}_{yyy}\la y\ra^{l+1}\|,
\end{split}
\end{align}
and 
\begin{align}
	\begin{split}
\sum_{i=1}^{m-1} \|\pa^i_x\bar{u}\la y\ra^{l_0} \|_{L^\infty} + \sum_{i=0}^{m-1} \|\pa^i_x\bar{u}_{y}\la y\ra^{l_0}\|_{L^\infty} &\lesssim \|\bar{u}_{yy}\la y\ra^{l_0+1}\|_{x=0} + \sum_{i=0}^{m-2} \|\pa_x^{i}\bar{v}_{yyy}\la y\ra^{l_0+2}\|_{x=0}\\
& \quad + \sqrt{L}\sum_{i=0}^{m-1} \|\pa_x^{i}\bar{v}_{yyy}\la y\ra^{l_0+2}\|,\\ 
\sum_{i=0}^{m-1} \|\pa_x^{i}q\|_{L^\infty} & \lesssim \sum_{i=0}^{m-1} \|\pa_x^{i}q_y\la y\ra\|_{x=0} + L^{\f23-}\|q_y\|_{\MX^m_{\la y\ra}}.\label{6.140-3}
	\end{split}
\end{align}
\end{lemma}

\noindent{\bf Proof.} It is clear 
\begin{align*}
	\begin{split}
		\|\pa_x^m\bar{u} \la y\ra^l\|_{L^\infty} &\lesssim \| \pa_x^{m-1}\bar{v}_{yy} \la y\ra^{l+1}\|_{L^\infty_x L^2_y}  \lesssim \||\bar{v}|\|_{m,w},\\
		\|\pa_x^m\bar{u} \la y\ra^{l}\|_{L^2_y} 
		&\lesssim\|\pa_x^m\bar{u}_y \la y\ra^{l+1}\|_{L^2_y}  \lesssim  \| \pa_x^{m-1}\bar{v}_{yy} \la y\ra^{l+1}\|_{L^2_y} \lesssim \||\bar{v}|\|_{m,w},\\
		\|\pa_x^m\bar{u}_y \la y\ra^{l}\|_{L^2_xL^\infty_y}&\lesssim \|\pa_x^m\bar{u}_{yy} \la y\ra^{l+1}\|\lesssim \|\pa_x^{m-1}\bar{v}_{yyy}\la y\ra^{l+1}\|,
	\end{split}
\end{align*}
which imply \eqref{6.140-2}.  

We have from Lemma \ref{lemA.5-2} that
\begin{align*}
	&\sum_{i=1}^{m-1} \|\pa^i_x\bar{u}\la y\ra^{l_0}\|_{L^\infty} + \sum_{i=0}^{m-1} \|\pa^i_x\bar{u}_{y}\la y\ra^{l_0}\|_{L^\infty}\nonumber\\
	&\lesssim \sum_{i=0}^{m-2} \|\pa_x^{i}\bar{v}_{yyy} \la y\ra^{l_0+2}\|_{L^\infty_xL^2_y}  + \|\bar{u}_{yy}\la y\ra^{l_0+1}\|_{x=0} + \sqrt{L} \|\bar{v}_{yyy}\la y\ra^{l_0+1}\| \nonumber\\
	&\lesssim  \|\bar{u}_{yy}\la y\ra^{l_0+1}\|_{x=0} + \sum_{i=0}^{m-2} \|\pa_x^{i}\bar{v}_{yyy}\la y\ra^{l_0+2}\|_{x=0} + \sqrt{L}\sum_{i=0}^{m-1} \|\pa_x^{i}\bar{v}_{yyy}\la y\ra^{l_0+2}\|,
\end{align*}
which yields $\eqref{6.140-3}_1$. Using \eqref{A6.61}, we derive $\eqref{6.140-3}_2$. Therefore the proof of Lemma \ref{lemA.6} is completed. $\hfill\Box$

\medskip

\begin{lemma}\label{lem6.5}
For $1\leq k\leq m$,	it holds that 
\begin{align}
\|\pa_x^{k-1}\bar{v}_{yyy}w\| 
&\lesssim L^{\f16-} P\big(\||\bar{v}|\|_{k,w}\big) + \|\pa_x^{k-1}F w\|, \label{6.126-0}\\
\|\pa_x^{k-1}\bar{v}_{yyyy}w\| 
&\lesssim P\big(\||\bar{v}|\|_{k,w}\big) +  \|\pa_x^{k-1}F_y w\|, \label{6.127-0}\\
\sum_{i=0}^{k}   \|\pa^i_x\bar{u}_{yyy}w\| 
&\lesssim \sqrt{L}\|\bar{u}_{yyy}w\|_{x=0} + P\big(\||\bar{v}|\|_{k,w}\big)  + \sum_{i=0}^{k-1} \|\pa_x^{i} F_{y} w\|,
\end{align}
and
\begin{align}\label{6.130-0}
\|\bar{u}_{\theta} \pa_x^{k-1} q_{yyy}w\| 
&\lesssim  P\Big(\sum_{i=0}^{k-1}\|\pa_x^iq_yw\|_{x=0}, \sqrt{L} \sum_{i=0}^{k-1}\| \pa_x^{i}q_{yy}w\bar{\chi}\|_{x=0},  \sqrt{L}\|\bar{u}_{yyy}w\|_{x=0} \Big)\nonumber\\
&\quad + P\Big(\||\bar{v}|\|_{k,w}, \sum_{i=0}^{k-1}\|\pa_x^{i}F w\|, \sum_{i=0}^{k-2} \|\pa_x^{i} F_{y} w\|\Big).
\end{align}
\end{lemma}

\noindent{\bf Proof.} 1. We have from \eqref{6.17} that 
\begin{align}
	\begin{split}
		\|\bar{v}_{yyy}w\|&\lesssim \|\bar{u}_{\theta}^2 q_{xy}w\| + \|Fw\| \lesssim L^{\f12} \|q_{y}\|_{\MX^1_{\t, w}} + \|Fw\|.
	\end{split}
\end{align}

Applying $\pa_y$ to \eqref{6.17}, one gets
\begin{align}\label{6.81}
	\bar{u}_{\theta}^2  q_{xyy} + 2\bar{u}_{\theta} \bar{u}_{y}   q_{xy}  - \mu \bar{v}_{yyyy}=F_y,
\end{align}
which yields that 
\begin{align}\label{6.81-0}
	\|\bar{v}_{yyyy}w\|&\lesssim \|\bar{u}_{\theta}^2  q_{xyy}w\| + \|\bar{u}_{\theta} \bar{u}_{y}   q_{xy}w\| + \|F_yw\| 
 \lesssim P\big(\||\bar{v}|\|_{1,w}\big) +  \|F_y w\| .
 \end{align}

\medskip

 It follows from \eqref{6.17} that 
\begin{align}\label{6.131-1}
	\mu \pa_x^{k-1}\bar{v}_{yyy}= \sum_{i=0}^{k-1} C_n^i \pa_x^{i+1}q_{y}\cdot \pa_x^{k-1-i} (\bar{u}_{\theta}^2) - \pa_x^{k-1}F,
\end{align}
which yields that 
\begin{align}\label{6.77}
	\|\pa_x^{k-1}\bar{v}_{yyy}w\| 
	&\lesssim \sum_{i=0}^{k-1} \big[1 + \|\pa_x^i \bar{u}\|^2_{L^\infty}\big] \cdot 
	\sum_{i=1}^{k}\|\pa_x^iq_{y}w\| + \|\pa_x^{k-1}F w\| \nonumber\\
	& \lesssim L^{\f16-} P\big(\||\bar{v}|\|_{k,w}\big) +  \|\pa_x^{k-1}F w\|.
\end{align}

Applying $\pa_y$ to \eqref{6.131-1}, one has that 
\begin{align}
\|\pa_x^{k-1}\bar{v}_{yyyy}w\| &\lesssim \|\bar{u}^2_{\theta}\pa_x^{k}q_{yy}w\|  + P\Big(\|q_{y}\|_{\MX^k_{\t, w}}, \,  \sum_{i=0}^{k-1}  \|(\pa_x^i \bar{u},\pa_x^i \bar{u}_y)\|_{L^\infty}  \Big) + \|\pa_x^{k-1} F_y w\|\nonumber\\
&\lesssim P\big(\||\bar{v}|\|_{k,w}\big)+  \|\pa_x^{k-1}F_y w\|.
\end{align}

\medskip

2. We have from  \eqref{6.126-0} that 
\begin{align}\label{6.142}
	\|\bar{u}_{yyy}w\| \lesssim  \sqrt{L}\|\bar{u}_{yyy}w\|_{x=0} +  \sqrt{L}\|\bar{u}_{xyyy}w\| 
	\lesssim  \sqrt{L}\|\bar{u}_{yyy}w\|_{x=0} +  \sqrt{L}\|\bar{v}_{yyyy}w\|,
\end{align}
which, together with \eqref{6.127-0}, yields  that 
\begin{align}\label{6.80}
\sum_{i=0}^{k}   \|\pa^i_x\bar{u}_{yyy}w\| &\lesssim  \|\bar{u}_{yyy}w\|_{x=0}  +    \sum_{i=0}^{k-1}   \|\pa^{i}_x\bar{v}_{yyyy}w\| \nonumber\\
&\lesssim \sqrt{L}\|\bar{u}_{yyy}w\|_{x=0} + P\big(\||\bar{v}|\|_{k,w}\big) + \sum_{i=0}^{k-1} \|\pa_x^{i} F_{y} w\|.
\end{align} 

\medskip

3.  A direct calculation shows that 
\begin{align}\label{6.137}
	\|\bar{u}_{\theta}q_{yyy}w\|
	&\lesssim \|\bar{v}_{yyy}w\| + \hat{e}_0\|\bar{u}_{yyy}w\|  + e_0\|q_{yy}w\| + \|\bar{u}_{yy}\|_{L^2_xL^\infty_y} \cdot \|q_yw\|_{L_x^{\infty}L^2_y} \nonumber\\
	&\lesssim P\Big(\||\bar{v}|\|_{1,w}, \|q_yw\|_{L_x^{\infty}L^2_y},  \|q_{yy}w\|, \|\bar{v}_{yyy}w\|, \|\bar{u}_{yyy} w\|\Big).
\end{align}

Noting
\begin{align}
	\bar{u}_{\theta} \pa_x^{k-1} q_{yyy}&= \pa_x^{k-1}\bar{v}_{yyy} - \bar{u}_{yyy} \pa_x^{k-1} q - 3\bar{u}_{yy} \pa_x^{k-1} q_y - 3\bar{u}_{y} \pa_x^{k-1} q_{yy} \nonumber\\
	&\quad - \sum_{j=1}^{k-1} C_{k-1}^j \big[\pa_x^j\bar{u}_{yyy} \, \pa_x^{k-1-j}q + 3\pa_x^j\bar{u}_{yy} \, \pa_x^{k-1-j}q_y  + 3\pa_x^j\bar{u}_{y} \, \pa_x^{k-1-j}q_{yy}\nonumber\\
	&\qquad\qquad + \pa_x^j\bar{u} \, \pa_x^{k-1-j}q_{yyy}\big],
\end{align}
which yields that 
\begin{align}\label{6.140}
\|\bar{u}_{\theta} \pa_x^{k-1} q_{yyy}w\|
&\lesssim P\Big(\||\bar{v}|\|_{k,w}, \|\pa_x^{k-1}\bar{v}_{yyy}w\|, \sum_{j=0}^{k-1} \|\pa_x^{j} q_{yy}w\|, \sum_{j=0}^{k-1} \|\pa_x^{j}q_y w\|_{L^\infty_xL^2_y}\Big)\nonumber\\
&\quad + P\Big(\sum_{j=0}^{k-1} \|\pa_x^j\bar{u}_{yyy} w\|,  \sum_{j=0}^{k-2} \|\bar{u}_{\theta}\pa_x^{j}q_{yyy} w\|\Big).
\end{align}
Applying induction arguments to \eqref{6.140}, and using \eqref{6.137}, one can obtain 
\begin{align*}
\|\bar{u}_{\theta} \pa_x^{k-1} q_{yyy}w\|
&\lesssim P\Big(\||\bar{v}|\|_{k,w},  \sum_{j=0}^{k-1} \big\{\|(\pa_x^{j}\bar{v}_{yyy}, \pa_x^{j} q_{yy}, \pa_x^j\bar{u}_{yyy})w\| + \|\pa_x^{j}q_y w\|_{L^\infty_xL^2_y}\big\}\Big)\nonumber\\
&\lesssim P\Big(\sum_{i=0}^{k-1}\|\pa_x^iq_yw\|_{x=0}, \sqrt{L} \sum_{i=0}^{k-1}\| \pa_x^{i}q_{yy}w\bar{\chi}\|_{x=0},  \sqrt{L}\|\bar{u}_{yyy}w\|_{x=0} \Big)  \nonumber\\
&\quad + P\Big(\||\bar{v}|\|_{k,w},  \sum_{i=0}^{k-1}\|\pa_x^{i}F w\|, \sum_{i=0}^{k-2} \|\pa_x^{i} F_{y} w\|\Big),
\end{align*}
where we have used \eqref{A6.62}, \eqref{6.38-0}, \eqref{6.77} and \eqref{6.80}. 
Therefore the proof of Lemma \ref{lem6.5} is completed. $\hfill\Box$

\subsubsection{Basic energy estimates on quotient}

\medskip

\begin{lemma}\label{lem6.8}
For $0\leq k\leq m$, it holds that 
\begin{align}\label{6.120}
& \|\bar{u}_{\theta}\pa_x^k q_y w\|^2_{x=x_0} + \|\sqrt{\bar{u}_{\t}}\pa^k_xq_{yy}w\|^2 \nonumber\\
&\lesssim \|\bar{u}_{\theta}\pa_x^k q_y w\|^2_{x=0}  + L^{\f16-}P\Big(\||\bar{v}|\|_{m,w},  \|\pa^{k-1}_x\bar{v}_{yyy} \la y\ra\|,  \|\bar{u}_{yy} w\|_{L^\infty_x L^2_y}, \sum_{i=0}^{k}\|\pa^i_x\bar{u}_{yyy}w\|\Big)\nonumber\\
&\quad + L^{\f16-} P\Big(\|q_{yy}w\|_{L^\infty_xL^2_y},  \sum_{i=0}^{k-1}\|\pa^{i}_xq_yw\|_{L^\infty_xL^2_y} ,   \sum_{i=1}^{k-1} \|(\bar{u}_{\theta}\pa_x^{i}q_{yyy}, \pa_x^{i}q_{yy})w\|, \|\pa_x^k F w\|\Big).
\end{align}
\end{lemma}

\noindent{\bf Proof.}  1. Applying $\pa^k_x$ to \eqref{6.17}, one gets
\begin{align}\label{6.121}
\pa^k_x\big(\bar{u}_{\theta}^2  q_{xy}\big)  - \mu \,\pa_x^k\pa_{yyy} [\bar{u}_{\theta}q] = \pa_x^kF,\,\,\,\mbox{for}\,\, 0\leq k\leq m.
\end{align}
Multiplying \eqref{6.121} by $\pa_x^k q_{y}w^2$, one has that 
\begin{align}\label{6.122}
\big\la \pa^k_x\big(\bar{u}_{\theta}^2  q_{xy}\big),\, \pa_x^kq_{y}w^2\big\ra -  \big\la \pa_x^k [\bar{u}_{\theta}q]_{yyy},\,  \pa^k_x q_y w^2\big\ra
\leq L^{\f16-} \|\pa_x^k F w\|^2 + L^{\f16-}\|\pa_x^kq_{y}\|^2_{\MX_{w}},
\end{align}
where we have used 
\begin{align*}
\big\la \pa_x^k F,\, \pa_x^kq_{y}w^2\big\ra \lesssim \|\pa_x^k F w\|\cdot \|\pa_x^kq_{y}w\| \lesssim L^{\f16-} \|\pa_x^k F w\|^2 + L^{\f16-}\|\pa_x^kq_{y}\|^2_{\MX_{w}}.
\end{align*}

2. Integrating by parts in $x$, and using \eqref{6.140-2}, one obtains
\begin{align}\label{6.123}
\big\la \pa^k_x\big(\bar{u}_{\theta}^2  q_{xy}\big),\, \pa_x^kq_{y}w^2\big\ra 
&=\big\la \bar{u}_{\theta}^2  \pa^k_xq_{xy},\, \pa_x^kq_{y}w^2\big\ra + \sum_{i=0}^{k-1} C_{k}^i \big\la \pa_x^i q_{xy} \pa_x^{k-i}(\bar{u}_{\theta}^2)  ,\, \pa_x^kq_{y}w^2\big\ra\nonumber\\
&\geq \f12 \|\bar{u}_{\theta}\pa_x^k q_y w\|^2_{x=x_0} - \f12 \|\bar{u}_{\theta}\pa_x^k q_y w\|^2_{x=0}  - L^{\f13-}P\big( \||\bar{v}|\|_{k,w}\big).
\end{align}

3.   It is clear that 
\begin{align}\label{6.124}
	-  \big\la \pa^k_x [\bar{u}_{\theta}q]_{yyy},\, \pa_x^kq_{y} w^2\big\ra = -  \big\la [\bar{u}_{\theta} \pa^k_xq]_{yyy},\, \pa_x^kq_{y}w^2\big\ra -  \sum_{i=1}^{k} C_k^i  \big\la [\pa^i_x\bar{u} \cdot \pa^{k-i}_xq]_{yyy},\,  \pa_x^kq_{y}w^2\big\ra .
\end{align}

3.1. For the first term on right hand side (RHS) of \eqref{6.124}, noting $\pa^k_xq_{y}|_{y=0}=0$, we have 
\begin{align}\label{6.125}
	-\big\la \pa_{yyy} [\bar{u}_{\theta} \pa^k_xq],\,   \pa^k_xq_{y}w^2\big\ra 
	&=\big\la \pa_{yy} [\bar{u}_{\theta}\pa^k_xq],\,   \pa^k_xq_{yy}w^2\big\ra + \big\la \pa_{yy} [\bar{u}_{\theta} \pa^k_xq],\,   2\pa^k_xq_{y}ww_y\big\ra \nonumber\\
	&=\|\sqrt{\bar{u}_{\t}}\pa^k_xq_{yy}w\|^2 + \big\la 2 \bar{u}_{y} \pa^k_xq_{y} + \bar{u}_{yy}\pa^k_xq,\,   \pa^k_xq_{yy} w^2\big\ra \nonumber\\
	&\quad  + \big\la \bar{u}_{\theta}\pa^k_xq_{yy} + 2 \bar{u}_{y} \pa^k_xq_{y} + \bar{u}_{yy}\pa^k_xq,\, 2\pa^k_xq_{y} ww_y\big\ra.
\end{align}

Integrating by parts in $y$ and using \eqref{A6.23}, one obtains
\begin{align}\label{6.126}
\big\la 2 \bar{u}_{y} \pa^k_xq_{y} ,\,   \pa^k_xq_{yy} w^2\big\ra
&= - \big\la  (\bar{u}_{y}w^2)_y,\,   |\pa^k_xq_{y}|^2 \big\ra 
\lesssim \|(\bar{u}_{y},\bar{u}_{yy})\|_{L^\infty} \|\pa^k_xq_{y}w\|^2  \lesssim L^{\f13-}  P\big(\||\bar{v}|\|_{k,w}\big),
\end{align}
and
\begin{align}\label{6.127}
\big\la \bar{u}_{\theta}\pa^k_xq_{yy} ,\, 2\pa^k_xq_{y} ww_y\big\ra + \big\la 2 \bar{u}_{y} \pa^k_xq_{y},\, 2\pa^k_xq_{y} ww_y \big\ra 
&\lesssim  L^{\f16-} P\big(\||\bar{v}|\|_{k,w}\big).
\end{align}
 
Noting $\pa^k_xq|_{y=0}=0$, we have 
\begin{align}\label{6.128}
	\big\la \bar{u}_{yy}\pa^k_xq,\,   \pa^k_xq_{yy} w^2\big\ra 
	&\lesssim \f{1}{\sqrt{a}}  \int_0^{x_0} \|\pa^k_xq_{y}\la y\ra\|_{L^2_y} \int_0^{\infty}  |\bar{u}_{yy}  \sqrt{\bar{u}} \pa^k_x q_{yy}| w^2   dydx\nonumber\\
	&\lesssim \f{1}{\sqrt{a}}  \|\bar{u}_{yy} w\|_{L^\infty_x L^2_y} \|\sqrt{\bar{u}} \pa^k_xq_{yy}w\|\cdot \|\pa^k_xq_yw\|\nonumber\\
	&\leq  \f1{8}  \|\sqrt{\bar{u}} \pa^k_xq_{yy}w\|^2  + C L^{\f13-} P\big( \|\bar{u}_{yy} w\|_{L^\infty_x L^2_y}, \|q_y\|_{\MX^k_{w}}\big),
\end{align}
and
\begin{align}\label{6.129}
	\big|\big\la \bar{u}_{yy}\pa^k_xq,\, 2\pa^k_xq_{y} ww_y \big\ra\big| 
	&\lesssim  \int_0^{x_0} \|\pa^k_xq_y\la y\ra \|_{L^2_y} \int_0^{\infty}  |\bar{u}_{yy}  \pa^k_xq_y ww_y|  dydx  \nonumber\\
	&\lesssim  \|\bar{u}_{yy}w\|_{L^\infty_x L^2_y}\cdot \|\pa^k_x q_yw\|^2
	\lesssim L^{\f13-} P\Big(\|\bar{u}_{yy}w\|_{L^\infty_x L^2_y}, \|q_y\|_{\MX^k_w}\Big).
\end{align}
Thus, substituting \eqref{6.126}-\eqref{6.129} into \eqref{6.125}, one obtains that
\begin{align}\label{6.130}
-\big\la \pa_{yyy} [\bar{u}_{\theta} \pa^k_xq],\,   \pa^k_xq_{y}w^2\big\ra   
&\geq \f78 \|\sqrt{\bar{u}_{\t}}\pa^k_xq_{yy}w\|^2 - L^{\f16-}P\big(\||\bar{v}|\|_{k,w}, \|\bar{u}_{yy} w\|_{L^\infty_x L^2_y}\big).
\end{align}

3.2. For the second term on RHS of \eqref{6.124}, we note that ($1\leq i\leq k$)
\begin{align}\label{6.131}
\big\la [\pa^i_x\bar{u} \pa^{k-i}_xq]_{yyy},\,  \pa_x^kq_{y}w^2\big\ra 
&=\big\la \pa^i_x\bar{u}_{yyy} \pa^{k-i}_xq,\,  \pa_x^kq_{y}w^2\big\ra  + \big\la 3\pa^i_x\bar{u}_{yy} \pa^{k-i}_xq_y,\,  \pa_x^kq_{y}w^2\big\ra
\nonumber\\
&\,\,\, + \big\la 3\pa^i_x\bar{u}_y  \pa^{k-i}_xq_{yy},\,  \pa_x^kq_{y}w^2\big\ra  + \big\la \pa^i_x\bar{u} \pa^{k-i}_xq_{yyy},\,  \pa_x^kq_{y}w^2\big\ra.
\end{align}
It is clear that 
\begin{align}
\sum_{i=1}^{k}\big\la \pa^i_x\bar{u}_{yyy} \pa^{k-i}_xq,\,  \pa_x^kq_{y}w^2\big\ra
&\lesssim \sum_{i=1}^{k}\|\pa_x^{k-i}q\|_{L^\infty} \|\pa^i_x\bar{u}_{yyy}w\|\cdot \|\pa_x^kq_{y}w\| \nonumber\\
&\lesssim L^{\f16-} P\Big(\||\bar{v}|\|_{k,w},  \sum_{i=1}^{k}\|\pa^i_x\bar{u}_{yyy}w\|\Big),
\end{align}
and
\begin{align}
\sum_{i=1}^{k}\big\la 3\pa^i_x\bar{u}_{yy} \pa^{k-i}_xq_y,\,  \pa_x^kq_{y}w^2\big\ra
&\lesssim \sum_{i=1}^{k}\|\pa^i_x\bar{u}_{yy}\|_{L^2_xL^\infty_y} \|\pa^{k-i}_xq_yw\|_{L^\infty_xL^2_y} \|\pa_x^kq_{y}w\| \nonumber\\
&\lesssim L^{\f16-} P\Big(\|q_y\|_{\MX^k_w}, \sum_{i=0}^{k-1}\|\pa^{i}_xq_yw\|_{L^\infty_xL^2_y} ,\sum_{i=1}^{k}\|\pa^i_x\bar{u}_{yyy}\la y \ra\| \Big).
\end{align}

Also, it is direct to have 
\begin{align} 
\sum_{i=1}^{k-1}\big\la 3\pa^i_x\bar{u}_y \pa^{k-i}_xq_{yy} + \pa^i_x\bar{u}  \pa^{k-i}_xq_{yyy},\,  \pa_x^kq_{y}w^2\big\ra   
&\lesssim L^{\f16-} P\Big(\||\bar{v}|\|_{k,w},  \sum_{i=1}^{k-1} \|(\bar{u}_{\theta}\pa_x^{i}q_{yyy}, \pa_x^{i}q_{yy})w\|\Big),
\end{align}
and
\begin{align}
\big\la 3\pa^k_x\bar{u}_y q_{yy},\,  \pa_x^kq_{y}w^2\big\ra
&\lesssim L^{\f16-}  \|q_y\|_{\MX^k_w} \|q_{yy}w\|_{L_x^\infty L^2_y} \|\pa^{k}_x\bar{u}_{y}\|_{L^2_xL_y^\infty} \nonumber\\
&\lesssim  L^{\f16-}  P\Big(\|q_{yy}w\|_{L_x^\infty L^2_y}, \|q_y\|_{\MX^k_w}, \|\pa^{k-1}_x\bar{v}_{yyy} \la y\ra\|\Big).
\end{align}
Integrating by parts in $y$, one gets
\begin{align}\label{6.134}
\big\la   \pa^k_x\bar{u}  q_{yyy},\,  \pa_x^kq_{y}w^2\big\ra
&=-\big\la\pa^k_x\bar{u} \, q_{yy},\,  \pa_x^kq_{yy}w^2\big\ra - \big\la   \pa^k_x\bar{u}_y \,  q_{yy},\,  \pa_x^kq_{y}w^2\big\ra  - \big\la \pa^k_x\bar{u} \, q_{yy},\,  2\pa_x^kq_{y} ww_y\big\ra \nonumber\\
&\lesssim \f{1}{\sqrt{a}}\|\sqrt{\bar{u}}\pa_x^kq_{yy}w\|\cdot \|\pa_x^{k-1}\bar{v}_{yy} \la y\ra \| \cdot \|q_{yy}w\|_{L_x^\infty L^2_y} \nonumber\\
&\quad + L^{\f16-}  \|q_y\|_{\MX^k_w} \|q_{yy}w\|_{L_x^\infty L^2_y} \|\pa^{k-1}_x\bar{v}_{yy} \la y\ra\| \nonumber\\
&\lesssim L^{\f12-} P\Big(\|q_{yy}w\|_{L^\infty_xL^2_y},\|q_y\|_{\MX^k_w},  \|\pa^{k-1}_x\bar{v}_{yy} \la y\ra\|_{L^\infty_xL^2_y} \Big).
\end{align}

Thus we have from \eqref{6.131}-\eqref{6.134} that
\begin{align}\label{6.135}
&-  \sum_{i=1}^{k} C_k^i  \big\la [\pa^i_x\bar{u} \cdot \pa^{k-i}_xq]_{yyy},\,  \pa_x^kq_{y}w^2\big\ra\nonumber\\
& \lesssim   L^{\f16-} P\Big(\|q_{yy}w\|_{L^\infty_xL^2_y}, \sum_{i=0}^{k-1}\|\pa^{i}_xq_yw\|_{L^\infty_xL^2_y} ,   \sum_{i=1}^{k-1} \|(\bar{u}_{\theta}\pa_x^{i}q_{yyy}, \pa_x^{i}q_{yy})w\|\Big)\nonumber\\
&\quad + L^{\f16-} P\Big(\||\bar{v}|\|_{k,w},  \|\pa^{k-1}_x\bar{v}_{yyy} \la y\ra\|,  \sum_{i=1}^{k}\|\pa^i_x\bar{u}_{yyy}w\| \Big) .
\end{align}

Combining \eqref{6.124}, \eqref{6.130} and \eqref{6.135}, one obtains that 
\begin{align}\label{6.136}
& -  \big\la \pa^k_x [\bar{u}_{\theta}q]_{yyy},\, \pa_x^kq_{y} w^2\big\ra \nonumber\\
&\geq \f78 \|\sqrt{\bar{u}_{\t}}\pa^k_xq_{yy}w\|^2 - L^{\f16-}P\big(\||\bar{v}|\|_{k,w},  \|\pa^{k-1}_x\bar{v}_{yyy} \la y\ra\|, \|\bar{u}_{yy} w\|_{L^\infty_x L^2_y}, \sum_{i=1}^{k}\|\pa^i_x\bar{u}_{yyy}w\|\big) \nonumber\\
&\quad - L^{\f16-} P\Big(\|q_{yy}w\|_{L^\infty_xL^2_y},  \sum_{i=0}^{k-1}\|\pa^{i}_xq_yw\|_{L^\infty_xL^2_y} ,   \sum_{i=1}^{k-1} \|(\bar{u}_{\theta}\pa_x^{i}q_{yyy}, \pa_x^{i}q_{yy})w\|\Big).
\end{align}

4. Substituting \eqref{6.136}  and \eqref{6.123} into \eqref{6.122}, one 
concludes \eqref{6.120}. Therefore the proof of Lemma \ref{lem6.8} is completed. $\hfill\Box$

\medskip


\begin{proposition}\label{prop4.9}
It holds that 
\begin{align}\label{6.139}
&\|q_y\|^2_{\MX^m_{\t, w}} +  \sum_{i=0}^{m-1} \|\pa^i_xq\|_{L^\infty}  + \sum_{i=1}^{m-1} \|\pa^i_x\bar{u}\la y\ra^{l_0}\|_{L^\infty} + \sum_{i=0}^{m-1} \|\pa^i_x\bar{u}_y\la y\ra^{l_0}\|_{L^\infty} + \sum_{i=0}^{1+\f{m}{2}} \|\pa^i_x\bar{u}_{yy}\la y\ra^{l_0}\|_{L^\infty} \nonumber\\
&\lesssim IN  + L^{\f16-}P\Big(\||\bar{v}|\|_{m,w},  \sum_{i=0}^{m}\|\pa_x^{i}F w\|, \sum_{i=0}^{m-1} \|\pa_x^{i} F_{y} w\|\Big),
\end{align}
where 
\begin{align}\label{6.139-1}
IN:&=1  + \sum_{k=0}^m\|\bar{u}_{\theta}\pa_x^k q_y w\|^2_{x=0}  +   L^{\f16-}P\Big(\|q_{yy}w\bar{\chi}\|_{x=0},  \sqrt{L} \sum_{i=0}^{m-1}\| \pa_x^{i}q_{yy}w\bar{\chi}\|_{x=0}\Big) \nonumber\\
&\quad + P\Big( \|\bar{u}_{yy} w\|_{x=0}, \sum_{i=0}^{m-2} \|\pa_x^{i}\bar{v}_{yyy} w\|_{x=0},  \sum_{i=0}^{m-1} \|\pa_x^{i}q_y w\|_{x=0}   , \sum_{i=0}^{1+\f{m}{2}} \|\pa_x^i \bar{u}_{yyy} w\|_{x=0}\Big).
\end{align}
\end{proposition}

\noindent{\bf Proof.} 1. It follows from \eqref{6.162} that 
\begin{align*}
	|\pa_x^i \bar{u}_{yy} \la y\ra^{l}| 
	&\lesssim \|\pa_x^i \bar{u}_{yyy} \la y\ra^{l+1}\|_{L^2_y} \lesssim \|\pa_x^i \bar{u}_{yyy} \la y\ra^{l+1}\|_{x=0} + \sqrt{L} \|\pa_x^{i+1} \bar{u}_{yyy} \la y\ra^{l+1}\|\nonumber\\
	&\lesssim \|\pa_x^i \bar{u}_{yyy} \la y\ra^{l+1}\|_{x=0} + \sqrt{L} \|\pa_x^{i} \bar{v}_{yyyy} \la y\ra^{l+1}\|,
\end{align*}
which implies that 
\begin{align}\label{6.122-0}
\sum_{i=0}^{1+\f{m}{2}} \|\pa^i_x\bar{u}_{yy}\la y\ra^{l_0}\|_{L^\infty} \lesssim \sum_{i=0}^{1+\f{m}{2}} \|\pa_x^i \bar{u}_{yyy} \la y\ra^{l_0+1}\|_{x=0} + \sqrt{L}\sum_{i=0}^{1+\f{m}{2}}\|\pa_x^{i} \bar{v}_{yyyy} \la y\ra^{l_0+1}\|.
\end{align} 
Hence we obtain from \eqref{6.122-0}, \eqref{6.140-3} and \eqref{6.126-0}-\eqref{6.127-0} that 
\begin{align}\label{6.150}
&\sum_{i=0}^{m-1} \|\pa^i_xq\|_{L^\infty}  + \sum_{i=1}^{m-1} \|\pa^i_x\bar{u}\la y\ra^{l_0}\|_{L^\infty} + \sum_{i=0}^{m-1} \|\pa^i_x\bar{u}_y\la y\ra^{l_0}\|_{L^\infty} + \sum_{i=0}^{1+\f{m}{2}} \|\pa^i_x\bar{u}_{yy}\la y\ra^{l_0}\|_{L^\infty} \nonumber\\
&\lesssim  IN   + L^{\f12-} P\Big(\||\bar{v}|\|_{m,w},  \sum_{i=0}^{m-1} \|\pa_x^iF \la y\ra^{l_0+2}\|,  \sum_{i=0}^{1+\f{m}{2}}\|\pa_x^{i}F_y \la y\ra^{l_0+1}\|  \Big).
\end{align}

\medskip

2. 
It follows from \eqref{6.30} and Lemmas \ref{lemA.5} \& \ref{lem6.5} that 
\begin{align}\label{6.123-0}
&\|(q_{yy}, \bar{u}_{yy})w\|_{L^\infty_xL^2_y} +  \sum_{i=0}^{k-1}\|\pa^{i}_xq_yw\|_{L^\infty_xL^2_y} +   \sum_{i=0}^{k-1} \|(\bar{u}_{\theta}\pa_x^{i}q_{yyy}, \pa_x^{i}q_{yy})w\|   + \sum_{i=0}^{k}   \|\pa^i_x\bar{u}_{yyy}w\|  \nonumber\\
&\lesssim P\Big(\|q_{yy}w\bar{\chi}\|_{x=0}, \sum_{i=0}^{k-1}\|\pa_x^iq_yw\|_{x=0}, \sqrt{L} \sum_{i=0}^{k-1}\| \pa_x^{i}q_{yy}w\bar{\chi}\|_{x=0} , \|\bar{u}_{yy}w\|_{x=0},   \sqrt{L}\|\bar{u}_{yyy}w\|_{x=0}\Big) \nonumber\\
&\quad  +  P\Big(\||\bar{v}|\|_{m,w},   \|\bar{v}_{yyy}\|_{x=0},   \sum_{i=0}^{k-1}\|(\pa_x^{i}F,\pa_x^{i}F_y) w\| \Big),
\end{align}
which, together with \eqref{6.120}, yields that 
\begin{align}\label{6.140-1}
\|q_y\|^2_{\MX^m_{\t, w}}
&\lesssim IN + L^{\f16-} P\Big(\||\bar{v}|\|_{m,w}, \sum_{i=0}^{m}\|\pa_x^{i}F w\|, \sum_{i=0}^{m-1} \|\pa_x^{i} F_{y} w\|\Big).
\end{align}
Combining  \eqref{6.150} and \eqref{6.140-1}, one  
concludes \eqref{6.139}. Therefore the proof of Proposition \ref{prop4.9} is completed. $\hfill\Box$


\medskip

\begin{lemma}\label{lem6.10-1}
It holds that 
\begin{align}
\sum_{i=0}^m\|\pa_x^i\bar{v}_yw\| &\lesssim L^{\f16-} P\Big(\|\bar{u}_yw\|_{x=0},  \||\bar{v}|\|_{m,w}\Big),\label{6.283-00}\\
\sum_{i=0}^{m-1}\|\pa_x^{i}\bar{v}_{yy}w\|^2_{x=x_0} 
&\lesssim \sum_{i=0}^{m-1}\|\pa_x^{i}\bar{v}_{yy}w\|^2_{x=0} + L^{\f16-} P\Big(\|\bar{u}_yw\|_{x=0},   \sum_{i=0}^{m-1}\|\pa_x^{i}Fw\|\Big),\label{6.283-0}\\
\|\pa_x^m\bar{v}_{yy}w\|
&\lesssim  P\Big(\|q_yw\|_{x=0}, \|\bar{u}_{yy}w\|_{x=0} ,  \||\bar{v}|\|_{m,w}, \sum_{i=0}^{m-1}\|\pa_x^iFw\|\Big). \label{6.283-000}
\end{align}
\end{lemma}

\noindent{\bf Proof.} 1. 
For $0\leq k\leq m$, we have 
\begin{align}
 \|\pa_x^k\bar{v}_yw\| 
&\lesssim \sum_{i=0}^k \big\{\|\pa_x^i\bar{u}\pa_x^{k-i}q_yw\|+  \|\pa_x^i\bar{u}_y\pa_x^{k-i}qw\|\big\} \nonumber\\
&\lesssim L^{\f16-} P\big(\||\bar{v}|\|_{m,w}, \|\pa_x^m \bar{u}\|_{L^\infty} \big)  + \||\bar{v}|\|_{m,w} \sum_{i=0}^{k-1} \|\pa_x^{i}\bar{v}_{yy}w\| \nonumber\\
&\quad  + \|w\bar{u}_y\|_{L^\infty_xL^2_y}\cdot \|\pa_x^kq\|_{L^2_xL^\infty_y} \nonumber\\
&\lesssim L^{\f16-} P\Big(\|\bar{u}_yw\|_{x=0}, \||\bar{v}|\|_{m,w}\Big),
\end{align}
where we have used  \eqref{6.140-2}, \eqref{6.30} and \eqref{6.115}.

\medskip

2.	For $0\leq k\leq m$, we have from \eqref{6.17} that 
\begin{align}\label{6.283} 
	\bar{u}_{\theta}^2 \pa_x^{k}q_y - \pa_x^{k-1}\bar{v}_{yyy}  + \sum_{i=1}^{k-1} C_{k-1}^i \pa^i_x(\bar{u}_{\theta}^2)\pa_x^{k-i}q_{y} = \pa_x^{k-1}F.
\end{align}
Multiplying \eqref{6.283} by $\pa_x^{k}\bar{v}_yw^2$, we obtain
\begin{align}\label{6.283-1}
\f12 \|\pa_x^{k-1}\bar{v}_{yy}w\|^2_{x=x_0} 
	&\leq \f12 \|\pa_x^{k-1}\bar{v}_{yy}w\|^2_{x=0} + \big\{\|\pa_x^{k-1}Fw\| + \|\pa_x^{k-1}\bar{v}_{yy}w\|\big\}\|\pa_x^k\bar{v}_yw\|\nonumber\\
	&\quad + L^{\f16-} P\big(\||\bar{v}|\|_{m,w}\big), 
\end{align}
where we have used 
\begin{align*}
- \big\la \pa_x^{k-1}\bar{v}_{yyy}, \pa_x^{k}\bar{v}_yw^2\big\ra
&=\big\la \pa_x^{k-1}\bar{v}_{yy}, \pa_x^{k}\bar{v}_{yy}w^2\big\ra + \big\la \pa_x^{k-1}\bar{v}_{yy}, \pa_x^{k}\bar{v}_y 2ww_y\big\ra\nonumber\\
&\geq \f12 \|\pa_x^{k-1}\bar{v}_{yy}w\|^2 \Big|_{x=0}^{x=x_0} -  2\|\pa_x^{k-1}\bar{v}_{yy}w\|\cdot \|\pa_x^k\bar{v}_yw\|,
\end{align*}
and
\begin{align*}
	&\big\la \bar{u}_{\theta}^2 \pa_x^{k}q_y + \sum_{i=1}^{k-1} C_{k-1}^i \pa^i_x(\bar{u}_{\theta}^2)\pa_x^{k-i}q_{y},\, \pa_x^{k}\bar{v}_y w^2\big\ra 
	\lesssim L^{\f16-} P\big(\||\bar{v}|\|_{m,w}\big).
\end{align*}
%
Substituting \eqref{6.283-00} into \eqref{6.283-1}, one obtains
\begin{align*}
\sum_{i=0}^{m-1}\|\pa_x^{i}\bar{v}_{yy}w\|^2_{x=x_0} 
&\lesssim \sum_{i=0}^{m-1}\|\pa_x^{i}\bar{v}_{yy}w\|^2_{x=0} + L^{\f16-} P\Big(\|\bar{u}_yw\|_{x=0}, \||\bar{v}|\|_{m,w}, \sum_{i=0}^{m-1}\|\pa_x^{i}Fw\|\Big),
\end{align*}
which concludes \eqref{6.283-0}.

\medskip

3. Noting
\begin{align*}
\pa_x^m\bar{v}_{yy}
&= \bar{u}_{\theta} \pa_x^mq_{yy}+ \sum_{i=1}^m C_m^i \pa_x^i\bar{u} \pa_x^{m-i}q_{yy}
	+ \sum_{i=0}^m 2 C_m^i \pa_x^i\bar{u}_{y}\pa_x^{m-i}q_{y} + \sum_{i=0}^m C_m^i \pa_x^i\bar{u}_{yy}\pa_x^{m-i}q,
\end{align*}
we have 
\begin{align*}
\|\pa_x^m\bar{v}_{yy}w\|&\lesssim \|\bar{u}_{\theta} \pa_x^mq_{yy}w\| +  \f{1}{\sqrt{a}} \sum_{i=0}^{m-1} \|\pa_x^{i}\bar{v}_{yy} \la y\ra\|_{L^\infty_xL^2_y}\cdot \sum_{i=0}^{m-1}\|\sqrt{\bar{u}}\pa_x^{i}q_{yy}\|\nonumber\\
&\quad + \|q_yw\|_{L^\infty L^2_y} \|\pa_x^{m-1}\bar{v}_{yyy}\la y\ra\| + L^{\f16-}  P\big(\||\bar{v}|\|_{m,w}\big)\nonumber\\
&\quad + \|\bar{u}_{yy}w\|_{L^\infty_xL^2_y}\cdot \|\pa_x^mq\|_{L^2_xL^\infty_y} + \||\bar{v}|\|_{m,w}\sum_{i=0}^{m-1}\|\pa_x^i\bar{v}_{yyy}w\|\nonumber\\
&\lesssim  P\Big(\|q_yw\|_{x=0}, \|\bar{u}_{yy}w\|_{x=0},  \||\bar{v}|\|_{m,w}, \sum_{i=0}^{m-1}\|\pa_x^iFw\|\Big),
\end{align*}
where we have used \eqref{6.115}, \eqref{A6.62}, \eqref{6.30} and \eqref{6.126-0}. 
Therefore the proof Lemma \ref{lem6.10-1}  is completed. $\hfill\Box$

\medskip

\begin{lemma}[Elliptic estimate]\label{lem6.12}
It holds that 
\begin{align}\label{6.266}
& \sum_{k=0}^{m} \Big\{\|\bar{u}_{\theta}^{\f52}\pa_x^k q_{yy}w\|^2_{x=x_0} +  \|\bar{u}_{\theta}^{3}\pa_x^{k+1} q_yw\|^2 + \|\bar{u}_{\theta} \pa_x^k\bar{v}_{yyy}w\|^2 \Big\}  \nonumber\\
&\lesssim \sum_{k=0}^{m}  \|\bar{u}_{\theta}^{\f52}\pa_x^k q_{yy}w\|^2_{x=0} + P\Big(IN, \||\bar{v}|\|_{m,w}, \sum_{k=0}^m\|\pa_x^k Fw\|, \sum_{k=0}^{m-1}\|\pa_x^iF_yw\|\Big).
\end{align}
We remark that the elliptic estimation \eqref{6.266} play no role for us to close the uniform-in-$\theta$ estimate. It only help us to know that $\pa_x^{m+1}q_y$ is well-defined in some sense.
\end{lemma}

\noindent{\bf Proof.} Let $0\leq k\leq m$. We have from \eqref{6.283} that 
\begin{align}\label{6.268}
&\|\bar{u}_{\theta}^{3}\pa_x^{k+1} q_yw\|^2 + \|\bar{u}_{\theta}\pa_x^k\bar{v}_{yyy}w\|^2 - 2\mu \big\la \bar{u}_{\theta}^{4} \pa_x^{k+1} q_y,\, \pa_x^k\bar{v}_{yyy}w^2\big\ra \nonumber\\
&\lesssim \|\pa_x^k Fw\|^2 + L^{\f13-}P\big(\||\bar{v}|\|_{m,w}\big),
\end{align}
where   we have used 
\begin{align*}
\sum_{i=1}^{k} \|\pa_x^i(\bar{u}_{\theta}^2)\, \pa_x^{k+1-i}q_{y} w\|^2 
&\lesssim L^{\f13-}P\big(\|\pa_x^m\bar{u}\|_{L^\infty}\big)\|q_y\|^2_{\MX_{\t, w}^{m}} \lesssim L^{\f13-}P\big(\||\bar{v}|\|_{m,w}\big).
\end{align*}

Integrating by parts, one obtains 
\begin{align}\label{6.269}
&- \big\la \bar{u}_{\theta}^{4} \pa_x^{k+1} q_y,\, \pa_x^k\bar{v}_{yyy}w^2\big\ra \nonumber\\
&=\big\la \bar{u}_{\theta}^{4} \pa_x^{k+1} q_{yy},\, \pa_x^k(\bar{u}_{\theta}q_{yy}) w^2\big\ra + \big\la \bar{u}_{\theta}^{4} \pa_x^{k+1} q_{yy},\, \pa_x^k(2\bar{u}_y q_{y}) w^2\big\ra  + \big\la \bar{u}_{\theta}^{4} \pa_x^{k+1} q_{yy},\, \pa_x^k(\bar{u}_{yy}q)w^2\big\ra\nonumber\\
&\quad + 4\big\la \bar{u}_{\theta}^{3}\bar{u}_y \pa_x^{k+1} q_y,\, \pa_x^k\bar{v}_{yy}w^2\big\ra + \big\la \bar{u}_{\theta}^{4} \pa_x^{k+1} q_y,\, \pa_x^k\bar{v}_{yy} 2ww_y\big\ra.
\end{align}

For the first term on RHS of \eqref{6.269}, 
a direct calculation shows that 
\begin{align}\label{6.277}
&\big\la \bar{u}_{\theta}^{4} \pa_x^{k+1} q_{yy},\, \pa_x^k(\bar{u}_{\theta}q_{yy}) w^2\big\ra\nonumber\\
&=\big\la \bar{u}_{\theta}^{4} \pa_x^{k+1} q_{yy},\, \bar{u}_{\theta} \pa_x^{k}q_{yy} w^2\big\ra 
+ \sum_{i=1}^k C_k^i \big\la \bar{u}_{\theta}^{4} \pa_x^{k+1} q_{yy},\, \pa_x^i\bar{u}_{\theta} \pa_x^{k-i}q_{yy}w^2\big\ra\nonumber\\
&=\f12 \|\bar{u}_{\theta}^{\f52}\pa_x^k q_{yy}w\|^2_{x=x_0} - \f12 \|\bar{u}_{\theta}^{\f52}\pa_x^k q_{yy}w\|^2_{x=0} - \f52 \big\la \bar{u}_{\theta}^{4}\bar{u}_x ,\,  |\pa_x^{k}q_{yy}|^2w^2\big\ra \nonumber\\
&\quad + \sum_{i=1}^k C_k^i \big\la \bar{u}_{\theta}^{4} \pa_x^{k+1} q_{yy},\, \pa_x^i\bar{u} \pa_x^{k-i}q_{yy}w^2\big\ra\nonumber\\
&\geq \f12 \|\bar{u}_{\theta}^{\f52}\pa_x^k q_{yy}w\|^2_{x=x_0} - \f12 \|\bar{u}_{\theta}^{\f52}\pa_x^k q_{yy}w\|^2_{x=0} - \f18 \|\bar{u}_{\theta}^{3}\pa_x^{k+1} q_yw\|^2  \nonumber\\
&\quad   - C P\Big(\||\bar{v}|\|_{m,w}, \|\pa_x^{m-1}\bar{v}_{yyy}\la y\ra\|,  \|q_{yy}w\|_{L^\infty_xL^2_y},   \sum_{i=0}^{m-1}\|(\bar{u}_{\theta}  \pa_x^i q_{yyy},  \pa_x^i q_{yy})w\|\Big),
\end{align}
where we have used 
\begin{align}\label{6.277-1}
&\sum_{i=1}^k C_k^i \big\la \bar{u}_{\theta}^{4} \pa_x^{k+1} q_{yy},\, \pa_x^i\bar{u} \pa_x^{k-i}q_{yy}w^2\big\ra\nonumber\\
&=-\sum_{i=1}^k C_k^i \big\la \bar{u}_{\theta}^{4} \pa_x^{k+1} q_{y},\, \pa_x^i\bar{u}  \pa_x^{k-i}q_{yyy}w^2\big\ra - (4)\sum_{i=1}^k C_k^i \big\la \bar{u}_{\theta}^{3} \bar{u}_y\pa_x^{k+1} q_{y},\, \pa_x^i\bar{u} \pa_x^{k-i}q_{yy}w^2\big\ra  \nonumber\\
&\quad - \sum_{i=1}^k C_k^i \big\la \bar{u}_{\theta}^{4} \pa_x^{k+1} q_{y},\, \pa_x^i\bar{u}_y  \pa_x^{k-i}q_{yy}w^2\big\ra - \sum_{i=1}^k C_k^i \big\la \bar{u}_{\theta}^{4} \pa_x^{k+1} q_{y},\, \pa_x^i\bar{u} \pa_x^{k-i}q_{yy} 2ww_y\big\ra\nonumber\\
&\geq - \f18 \|\bar{u}_{\theta}^{3}\pa_x^{k+1} q_yw\|^2 - C P\Big( \|\pa_x^{m-1}\bar{v}_{yyy}\la y\ra\|,  \|q_{yy}w\|_{L^\infty_xL^2_y},   \sum_{i=0}^{m-1}\|(\bar{u}_{\theta} \pa_x^i q_{yyy}, \pa_x^i q_{yy})w\|\Big) \nonumber\\
&\quad - C P\big(\||\bar{v}|\|_{m,w}\big).
\end{align}
Here we have used $\eqref{6.140-2}_2$ in \eqref{6.277}.

For the second term on RHS of \eqref{6.269}, it is clear that
\begin{align}\label{6.274-0}
&|\big\la \bar{u}_{\theta}^{4} \pa_x^{k+1} q_{yy},\, \pa_x^k(2\bar{u}_y q_{y}) w^2\big\ra|\nonumber\\
&\leq 2 |\big\la \bar{u}_{\theta}^{4} \pa_x^{k+1} q_{y},\, \pa_x^k(\bar{u}_y q_{yy} + \bar{u}_{yy} q_{y}) w^2\big\ra| + 4 |\big\la \bar{u}_{\theta}^{4} \pa_x^{k+1} q_{y},\, \pa_x^k(\bar{u}_y q_{y}) ww_y\big\ra| \nonumber\\
&\quad + 4 |\big\la \bar{u}_{\theta}^{3} \bar{u}_y \pa_x^{k+1} q_{y},\, \pa_x^k(2\bar{u}_y q_{y}) w^2\big\ra| \nonumber\\
&\leq  \f18 \|\bar{u}_{\theta}^{3}\pa_x^{k+1} q_yw\|^2 + C P\Big(\||\bar{v}|\|_{m,w},   \|q_{yy}w\|_{L^\infty_xL^2_y},   \sum_{i=0}^{m}\|\bar{u}_{\theta}  \pa_x^i q_{yy}w\|, \sum_{i=0}^{m-1}\|\pa_x^{i}\bar{v}_{yyy}\la y\ra\|\Big) \nonumber\\
&\quad + C P\Big(\sum_{i=1}^m \|\pa_x^i \bar{u}_{yyy}\la y\ra\|, \sum_{i=0}^{\f{m}{2}} \| \pa_x^i q_yw\|_{L^\infty_xL^2_y}\Big).
\end{align}
For 3th term on RHS of \eqref{6.269}, one gets
\begin{align}\label{6.275-0}
|\big\la \bar{u}_{\theta}^{4} \pa_x^{k+1} q_{yy},\, \pa_x^k(\bar{u}_{yy}q)w^2\big\ra| 
&\leq |\big\la \bar{u}_{\theta}^{4} \pa_x^{k+1} q_{y},\, \pa_x^k(\bar{u}_{yy}q)_y w^2\big\ra| + |\big\la \bar{u}_{\theta}^{4} \pa_x^{k+1} q_{y},\, \pa_x^k(\bar{u}_{yy}q)2w w_y\big\ra| \nonumber\\
&\quad + 4 |\big\la \bar{u}_{\theta}^{3} \bar{u}_y \pa_x^{k+1} q_{y},\, \pa_x^k(\bar{u}_{yy}q)w^2\big\ra| \nonumber\\
&\leq  \f18 \|\bar{u}_{\theta}^{3}\pa_x^{k+1} q_yw\|^2 + P\Big(\|(\bar{u}_{yy}, \bar{u}_{yyy})w\|_{x=0}, \||\bar{v}|\|_{m,w},  \sum_{i=1}^m \|\pa_x^i \bar{u}_{yyy} w\|\Big)\nonumber\\
&\quad + P\Big(\sum_{i=0}^{m-1}\|\pa_x^{i}\bar{v}_{yyy}\la y\ra\|\Big),
\end{align}
where we have  used the following fact
\begin{align*}
\|\bar{u}_{\theta} \pa_x^k(\bar{u}_{yy}q)_yw\|^2 
&\lesssim  P\Big(\|(\bar{u}_{yy}, \bar{u}_{yyy})w\|_{x=0}, \||\bar{v}|\|_{m,w},  \sum_{i=1}^m \|\pa_x^i \bar{u}_{yyy} w\|, \sum_{i=0}^{\f{m}{2}} \|\bar{u}_{\theta} \pa_x^i q_yw\|_{L^\infty_xL^2_y}\Big)\nonumber\\
&\quad + P\Big(\sum_{i=0}^{m-1}\|\pa_x^{i}\bar{v}_{yyy}\la y\ra\|\Big).
\end{align*}
For the last two terms on RHS of \eqref{6.269}, we have 
\begin{align}\label{6.282}
&4|\big\la \bar{u}_{\theta}^{3}\bar{u}_y \pa_x^{k+1} q_y,\, \pa_x^k\bar{v}_{yy}w^2\big\ra| + |\big\la \bar{u}_{\theta}^{4} \pa_x^{k+1} q_y,\, \pa_x^k\bar{v}_{yy} 2ww_y\big\ra|\nonumber\\
&\leq  \f18 \|\bar{u}_{\theta}^{3}\pa_x^{k+1} q_yw\|^2  +  C P\big(\||\bar{v}|\|_{m,w}.  \|\pa_x^m \bar{v}_{yy}w\| \big),
\end{align}

Substituting \eqref{6.269}, \eqref{6.277}-\eqref{6.282} into \eqref{6.268},  one obtains
\begin{align*} 
& \|\bar{u}_{\theta}^{\f52}\pa_x^k q_{yy}w\|^2_{x=x_0} +  \|\bar{u}_{\theta}^{3}\pa_x^{k+1} q_yw\|^2 + \|\bar{u}_{\theta}\pa_x^k\bar{v}_{yyy}w\|^2   \nonumber\\
&\lesssim  \|\bar{u}_{\theta}^{\f52}\pa_x^k q_{yy}w\|^2_{x=0} + P\Big(\|(\bar{u}_{yy}, \bar{u}_{yyy})w\|_{x=0}, \||\bar{v}|\|_{m,w},  \sum_{i=1}^m \|\pa_x^i \bar{u}_{yyy} w\|, \|\pa_x^m \bar{v}_{yy}w\| , \sum_{i=0}^{m-1}\|\pa_x^{i}\bar{v}_{yyy}w\|\Big)\nonumber\\
&\quad + P\Big( \|q_{yy}w\|_{L^\infty_xL^2_y},  \sum_{i=0}^{\f{m}{2}} \|\pa_x^i q_yw\|_{L^\infty_xL^2_y}, \sum_{i=0}^{m-1}\|(\bar{u}_{\theta} \pa_x^i q_{yyy},  \pa_x^i q_{yy})w\|, \|\pa_x^k Fw\|\Big),
\end{align*}
which, together with \eqref{6.123-0}, \eqref{6.283-000} and \eqref{6.126-0}, concludes \eqref{6.266}. Therefore the proof Lemma \ref{lem6.12} is completed. $\hfill\Box$

\medskip
 
  
\subsection{Estimates on temperature} \label{sec2.5}
For later use, we have from Lemma \ref{lem6.1} that 
\begin{align}\label{6.161}
	\|\pa_x^ih w\| \lesssim L^{\f16-} \big\{\|\sqrt{\bar{u}} \pa_x^ih_y w\| + \sup_{x\in[0,x_0]} \|\bar{u}\pa_x^ih w\|_{L^2_y}\big\} \lesssim L^{\f16-} \|\pa_x^ih\|_{\MY_w},\quad \mbox{for}\,\,\,i\geq 0.
\end{align}


 
\subsubsection{Basic energy estimate} Applying $\pa_x^k$ to \eqref{6.195-3}, one has
\begin{align}\label{6.173}
	&2 \bar{u}_{\theta}\pa_x^{k}\Theta_x +  2 \bar{v} \pa_x^k\Theta_y + \sum_{i=1}^k 2 C_k^i \pa_x^i \bar{u}_{\theta}\cdot \pa_x^{k-i}\Theta_x + \sum_{i=1}^k 2C_k^i \pa_x^i\bar{v}\cdot\pa_x^{k-i}\Theta_y \nonumber\\
	& = \kappa  \pa_x^k\Theta_{yy} + \mu \pa_x^k \{|([T_B+\Theta]\bar{u})_y|^2\}  + \pa_x^k G_B.
\end{align}
\begin{lemma}\label{lem6.10}
It holds that 
\begin{align}\label{6.169-0}
\|\Theta \|_{\MY^{m}_{\t, w}}^2
&\lesssim \sum_{i=0}^{m} \|\sqrt{\bar{u}_{\theta}}\pa_x^i\Theta w\|^2_{x=0} +  L^{\f16-} P\Big(\|\bar{u}_y w\|_{x=0}, \||\bar{v}|\|_{m,w},  [[[\Theta]]]_{m+1,w}\Big).
\end{align}
\end{lemma}

\noindent{\bf Proof.} Let $0\leq k\leq m$.
Multiplying \eqref{6.173} by $\pa_x^k\Theta w^2$, one obtains that 
\begin{align}\label{6.179}
&\|\sqrt{\bar{u}_{\theta}}\pa_x^k\Theta  w\|^2_{x=x_0}+ \kappa \|\pa_x^k\Theta_y w\|^2 \nonumber\\
&=\|\sqrt{\bar{u}_{\theta}}\pa_x^k\Theta w\|^2_{x=0} + \big\la 2\bar{v} ww_y,\, |\pa_x^k\Theta|^2\big\ra - \kappa \big\la \pa_x^k\Theta_y,\,  \pa_x^k\Theta 2 ww_y \big\ra\nonumber\\
&\quad -\sum_{i=1}^k\big\la 2 C_k^i \pa_x^i \bar{u}_{\theta}\cdot \pa_x^{k-i}\Theta_x,\, \pa_x^k\Theta w^2\big\ra - \sum_{i=1}^k\big\la 2C_k^i \pa_x^i\bar{v}\cdot\pa_x^{k-i}\Theta_y,\, \pa_x^k\Theta w^2\big\ra\nonumber\\
&\quad + \big\la \mu \pa_x^k \{|([T_B+\Theta]\bar{u})_y|^2\} ,\, \pa_x^k\Theta w^2\big\ra  + \big\la \pa_x^m G_B,\, \pa_x^k\Theta w^2\big\ra.
\end{align}
Using \eqref{6.161}, one obtains 
\begin{align}\label{6.180}
\begin{split}
&\big\la 2\bar{v} ww_y,\, |\pa_x^k\Theta|^2\big\ra - \kappa \big\la \pa_x^k\Theta_y,\,  \pa_x^k\Theta 2 ww_y \big\ra \lesssim L^{\f16-} \|\Theta\|^2_{\MY^{k}_{\t, w}},\\
&\sum_{i=1}^k C_{k}^i \big\la \pa_x^i \bar{u}_{\theta}\cdot \pa_x^{k-i}\Theta_x,\, \pa_x^k\Theta \, w^2\big\ra \lesssim L^{\f13-}  \sum_{i=1}^{k} \|\pa_x^i \bar{u}\|_{L^\infty} \|\Theta \|^2_{\MY^k_{\t, w}} \lesssim L^{\f13-}  \||\bar{v}|\|_{m,w} \|\Theta \|^2_{\MY^k_{\t, w}}.
\end{split}
\end{align}
It follows from \eqref{6.161} and \eqref{6.55} that 
\begin{align}\label{6.181}
&\sum_{i=1}^{k} C_k^i\big\la  \pa_x^i\bar{v}\cdot\pa_x^{k-i}\Theta_y,\, \pa_x^k\Theta w^2\big\ra \nonumber\\
&\lesssim  \|\pa_x^k\Theta w\| \Big\{\sum_{i=1}^{k-1} \|\pa_x^{i}\bar{v}\|_{L^\infty} \cdot \sum_{i=1}^{k-1}\|\pa_x^{i}\Theta_y w\| + \|\pa^k_x\bar{v}\|_{L^2_xL^\infty_y}\|\Theta_yw\|_{L^\infty_xL^2_y}\Big\} \nonumber\\
&\lesssim   L^{\f16-} P\big([[[\Theta]]]_{m+1,w}, \||\bar{v}|\|_{m,w}\big),
\end{align}
where we have used the following fact
\begin{align}\label{6.147}
	\|\pa_x^k\bar{v}\|_{L^2_xL^\infty_y} \lesssim \|\pa_x^k\bar{v}_y\la y\ra\| = \|\pa_x^k\bar{v}_y\la y\ra\|\lesssim \||\bar{v}|\|_{m,w}\quad \mbox{for}\,\,\,0\leq k\leq m.
\end{align}

For the last two terms, we   obtain 
\begin{align}\label{6.18-0}
&\big\la \pa_x^k \{|([T_B+\Theta]\bar{u})_y|^2\},\, \pa_x^k\Theta w^2\big\ra \nonumber\\
&\lesssim L^{\f16-} \|\Theta\|_{\MY^k_{\theta,w}} \sum_{i=0}^k \|\pa_x^i ([T_B+\Theta]\bar{u})_y  \pa_x^{k-i} ([T_B+\Theta]\bar{u})_y  \, w \|\nonumber\\
 &\lesssim L^{\f16-}  P\Big(\||\bar{v}|\|_{m,w}, \sum_{i=0}^{\f{k}{2}} \|(\pa_x^i\Theta,\pa_x^i\Theta_y)\la y\ra^{l_0}\|_{L^\infty},  \|\Theta\|_{\MY^k_{\t, w}},  \sum_{i=1}^{k}\|(\pa_x^i\bar{u},\pa_x^i\bar{u}_y) w^{\f34}\|,  \|\bar{u}_yw\|\Big) \nonumber\\
 &\lesssim L^{\f16-} P\Big(\|\bar{u}_y w\|_{x=0}, \||\bar{v}|\|_{m,w},  [[[\Theta]]]_{m+1,w}\Big),
\end{align}
and 
\begin{align}\label{6.182}
\big\la \pa_x^m G_B,\, \pa_x^k\Theta w^2\big\ra \lesssim L^{\f16-}\|\Theta\|_{\MY^k_{\t, w}}.
\end{align}

Substituting  \eqref{6.180}-\eqref{6.182} into \eqref{6.179}, one gets
\begin{align*}
&\|\sqrt{\bar{u}_{\theta}}\pa_x^k\Theta w\|^2_{x=x_0}+ \kappa \|\pa_x^k\Theta_y w\|^2 \nonumber\\
&\lesssim \|\sqrt{\bar{u}_{\theta}}\pa_x^k\Theta w\|^2_{x=0} + L^{\f16-} P\Big(\|\bar{u}_yw\|_{x=0}, \||\bar{v}|\|_{m,w},  [[[\Theta]]]_{m+1,w}\Big),
\end{align*}
which immediately yields \eqref{6.169-0}.
Therefore the proof of Lemma \ref{lem6.10} is completed. $\hfill\Box$


\subsubsection{Estimates for higher order derivatives}\label{sec2.4.2}
Noting the definition of $F$, to close \eqref{6.139}, we have to control $\|\bar{u}\pa_x^{m+1}T_{yy}w\|$.  Due to the appearance of $\pa_x^{m+1} \bar{v}=\bar{u}_{\theta}\pa_x^{m+1}q + \cdots$, it is hard to use the arguments in Lemma  \ref{lem6.10} to derive $\|\pa_x^{m+1}\Theta\|_{\MY_{\t, w}}$. To overcome the difficulty, we introduce the {\it Pseudo entropy} $S_k$
\begin{align}\label{6.231}
S_k:=\pa_x^k\Theta + \pa_x^{k-1}q\, \Theta_y\quad \Longleftrightarrow\quad \pa_x^k\Theta = S_k - \Theta_y\pa_x^{k-1}q.
\end{align}
Noting
\begin{align*}
	\pa_x^k\bar{v} = \pa_x^k(\bar{u}_{\theta}q) = \bar{u}_{\theta} \pa_x^k q + \sum_{i=1}^{k} C_k^i \pa_x^i\bar{u} \pa_x^{k-i} q,
\end{align*}
we can rewrite  \eqref{6.173} as
\begin{align}\label{6.234}
&2 \bar{u}_{\theta} \pa_x S_k +  2 \bar{v} \pa_yS_k - \kappa \pa_y^2 S_{k}\nonumber\\
&= 2 \bar{u}_{\theta} \Theta_{xy} \pa_x^{k-1}q  + 2 \bar{v} \big(\Theta_y\pa_x^{k-1}q\big)_y - \kappa \big(\Theta_y\pa_x^{k-1}q\big)_{yy} - \sum_{i=1}^{k} 2  C_{k}^i  \Theta_y\pa_x^i\bar{u} \pa_x^{k-i} q  \nonumber\\
&\quad-\sum_{i=1}^k 2 C_k^i \pa_x^i \bar{u}_{\theta}\cdot \pa_x^{k-i}\Theta_x- \sum_{i=1}^{k-1} 2C_k^i \pa_x^i\bar{v}\cdot\pa_x^{k-i}\Theta_y  + \mu \pa_x^k \{|([T_B+\Theta]\bar{u})_y|^2\}  + \pa_x^k G_B.
\end{align}
It is easy to know that 
\begin{align}
\begin{split}
\pa_yS_k\big|_{y=0}&=0,\quad \mbox{for NBC},\\
S_k\big|_{y=0} &=0,\quad  \mbox{for DBC}.
\end{split}
\end{align}

For later use, we have from \eqref{6.231} and \eqref{6.115} that 
\begin{align}
\|\pa_x^k\Theta w\| &\leq  \|S_k w\| + \|\Theta_yw\|_{L^\infty_xL^2_y} \|\pa_x^{k-1}q\|_{L^2_xL^\infty_y}\nonumber\\
&\leq \|S_k w\| + L^{\f16-} \|\Theta_yw\|_{L^\infty_xL^2_y} \|q_y \|_{\MX^{k-1}_{\la y\ra}},\label{6.188}\\
\|\pa_x^k\Theta_y w\| &\leq  \|\pa_yS_k w\|  +   \|\Theta_{yy} w\|_{L^\infty_xL^2_y} \|\pa_x^{k-1}q\|_{L^2_xL^\infty_y} + \|\Theta_y\|_{L^\infty} \|\pa_x^{k-1}q_y w\| \nonumber\\
&\leq  \|\pa_yS_k w\|  + L^{\f16-}   \|\Theta_{yy}w\|_{L^\infty_xL^2_y} \|q_y \|_{\MX^{k-1}_{\la y\ra}}.\label{6.188-1}
\end{align}

\smallskip
 
\begin{lemma}[Energy estimate]\label{lem6.14}
It holds that 
\begin{align}\label{6.173-0}
& \|\sqrt{\bar{u}_{\theta}}S_{m+1}w\|^2_{x=x_0}  + \kappa \|\pa_yS_{m+1}w\|^2  \nonumber\\
&\lesssim \|\sqrt{\bar{u}_{\theta}}S_{m+1} w\|^2_{x=0} + L^{\f16-} P\Big(\|\bar{u}_yw\|_{x=0},  \||\bar{v}|\|_{m,w},  [[[\Theta]]]_{m+1,w}, \|\pa_x^m\bar{v}_{yy}w\|\Big).
\end{align}
\end{lemma}

\noindent{\bf Proof.} 
Multiplying \eqref{6.234} by $S_k w^2$ ($k=m+1$), one obtains 
\begin{align}\label{6.235}
& \|\sqrt{\bar{u}_{\theta}}S_kw\|^2_{x=x_0}  + \kappa \|\pa_yS_kw\|^2  \nonumber\\
&\leq \|\sqrt{\bar{u}_{\theta}}S_k w\|^2_{x=0}+ \big\la 2 \bar{u}_{\theta} \Theta_{xy} \pa_x^{k-1}q ,\, S_k w^2\big\ra   + \big\la 2 \bar{v} \big(\Theta_y\pa_x^{k-1}q\big)_y,\, S_k w^2\big\ra\nonumber\\
&\quad - \kappa \big\la \big(\Theta_y\pa_x^{k-1}q\big)_{yy},\, S_k w^2\big\ra - \sum_{i=1}^{k} 2  C_{k}^i  \big\la  \Theta_y\pa_x^i\bar{u} \pa_x^{k-i} q,\, S_k w^2\big\ra   - \sum_{i=1}^k 2 C_k^i \big\la \pa_x^i \bar{u}  \pa_x^{k-i}\Theta_x,\, S_k w^2\big\ra  \nonumber\\
&\quad - \sum_{i=1}^{k-1} 2C_k^i  \big\la \pa_x^i\bar{v}\cdot\pa_x^{k-i}\Theta_y,\, S_k w^2\big\ra +  \big\la \mu \pa_x^k \{|([T_B+\Theta]\bar{u})_y|^2\} +  \pa_x^k G_B ,\, S_kw^2\big\ra\nonumber\\
&\quad  + L^{\f16-} [1+\|\bar{v}_{yy}\la y\ra^3\|_{x=x_0}]\|S_k\|^2_{\MY_{\t, w}} 
\end{align}
 
\smallskip

 A direct calculation shows that 
\begin{align}
\begin{split}
\|\bar{u}_{\theta} \Theta_{xy} \pa_x^{k-1}q w\| 
&\lesssim \|\Theta_{xy}w\|_{L^\infty_xL^2_y} \|\pa_x^{k-1}q\|_{L^2_xL^\infty_y}\lesssim L^{\f16-} \|\Theta_{xy}w\|_{L^\infty_xL^2_y} \|q_y\|_{\MX^{k-1}_{\t, \la y\ra}},\label{6.190}\\
\|\bar{v} \big(\Theta_y\pa_x^{k-1}q\big)_yw\| 
&\lesssim L^{\f16-} \|\bar{v}_{yy}\la y\ra^2 \|_{L^\infty_xL^2_y}\|(\Theta_y,  \Theta_{yy})w\|_{L^\infty_xL^2_y} \|q_y\|_{\MX^{k-1}_{\t, \la y\ra}}, 
\end{split}
\end{align}
Integrating by parts in $y$, one has 
\begin{align}\label{6.238}
& - \big\la \big(\Theta_y\pa_x^{k-1}q\big)_{yy},\, S_k w^2\big\ra \nonumber\\
& =  \big\la \big(\Theta_y\pa_x^{k-1}q\big)_{y},\, \pa_y S_k w^2\big\ra + \big\la \big(\Theta_y\pa_x^{k-1}q\big)_{y},\, 2ww_y S_k \big\ra\nonumber\\
&\lesssim \|(\pa_y S_k,S_k) w\| \Big\{\|\Theta_y\|_{L^\infty} \|\pa_x^{k-1}q_y w\| + \|\Theta_{yy}w\|_{L^\infty_xL^2_y}\|\pa_x^{k-1}q\|_{L^2_xL^\infty_y}\Big\} \nonumber\\
&\lesssim  L^{\f16-}P\big(\||\bar{v}|\|_{m,w},  [[[\Theta]]]_{m+1,w}\big).
\end{align}
where we have used \eqref{6.115}.


Then,  for the 4th, 5th, 6th, 7th and 8th terms on RHS of \eqref{6.234}, one has
\begin{align}\label{6.149}
	\sum_{i=1}^{k} \|\Theta_y\pa_x^i\bar{u} \pa_x^{k-i} q  w\|  
	&\lesssim \|\bar{u}_x\|_{L^\infty} \|\Theta_y w\|_{L^\infty_xL^2_y} \|\pa_x^{k-1}q\|_{L^2_xL^\infty_y} + \|\Theta_y\|_{L^\infty} \sum_{i=0}^{k-2}\|\pa_x^iq\|_{L^\infty} \sum_{i=1}^{k-1} \|\pa_x^i\bar{v}_yw\| \nonumber\\
	&\lesssim P\big(\||\bar{v}|\|_{m,w},  [[[\Theta]]]_{m+1,w}\big),
\end{align}
\begin{align}\label{6.150-1}
	\sum_{i=1}^k \|\pa_x^i \bar{u}_{\theta}\cdot \pa_x^{k-i}\Theta_x w\| 
	&\lesssim \||\bar{v}|\|_{m,w} \sum_{i=3}^{k-1}\|\pa_x^i\Theta w\|   + [[[\Theta]]]_{m+1,w} \|(\pa_x^{k-1}\bar{v}_y,\pa_x^{k-2}\bar{v}_y)w\|  \nonumber\\
	&\lesssim    P\Big(\||\bar{v}|\|_{m,w}, [[[\Theta]]]_{m+1,w}\Big).
\end{align}

%
\begin{align}\label{6.151-0}
	\sum_{i=1}^{k-1} \|\pa_x^i\bar{v}\cdot\pa_x^{k-i}\Theta_y w\| 
	&\lesssim  \|(\Theta_{xy},\Theta_{xxy}) w\|_{L^\infty_xL^2_y}  \|(\pa_x^{k-1}\bar{v}, \pa_x^{k-2} \bar{v})\|_{L^2_xL^\infty_y} \nonumber\\
	&\quad  +  \sum_{i=0}^{m-2} \|\pa_x^{i}\bar{v}_{yy} \la y\ra^2 \|_{L^\infty_x L^2_y}  \Big\{ \|\Theta \|_{\MY^{m}_w} + \|\pa_x^m\Theta_y w\| \Big\}\nonumber\\
	&\lesssim 
  P\Big(\||\bar{v}|\|_{m,w}, [[[\Theta]]]_{m+1,w}\Big),
\end{align}
and
\begin{align}\label{6.197}
	&\|\pa_x^k \{|([T_B+\Theta]\bar{u})_y|^2\} w\|  + \|\pa_x^k G_B w\| \nonumber\\
	&\lesssim P\Big(\||\bar{v}|\|_{m,w}, \sum_{i=0}^{\f{m}{2}} \|(\pa_x^i\Theta,\pa_x^i\Theta_y)\la y\ra^{l_0}\|_{L^\infty},  \sum_{i=0}^k\|(\pa_x^i\Theta, \pa_x^i\Theta_y,  \pa_x^i\bar{u}_y) w\|, \sum_{i=1}^k\|\pa_x^i\bar{u} w\|\Big) \nonumber\\
	&\lesssim P\Big(\|\bar{u}_yw\|_{x=0}, \||\bar{v}|\|_{m,w}, [[[\Theta]]]_{m+1,w}, \|\pa_x^m\bar{v}_{yy}w\| \Big),
\end{align}
where we have used \eqref{6.165-1},  \eqref{6.115}, \eqref{6.188}-\eqref{6.188-1} and \eqref{6.147}.

\smallskip 

Finally, substituting \eqref{6.190}-\eqref{6.197}   into \eqref{6.235}, one  concludes \eqref{6.173-0}. 
 Therefore the proof of Lemma \ref{lem6.14} is completed. $\hfill\Box$

\medskip
\begin{lemma}[Elliptic estimate]\label{lem6.15}
It holds that 
\begin{align}\label{6.189}
& \|\bar{u}^2_{\theta} \pa_x S_{m+1} w\|^2 +  \|\bar{u}_{\theta}\pa_y^2 S_{m+1}w\|^2 +  \|\bar{u}_{\theta}^{\f32} \pa_yS_{m+1}\|^2_{x=x_0}  \nonumber\\
&\lesssim \|\bar{u}_{\theta}^{\f32} \pa_yS_{m+1}\|^2_{x=0} +  P\big(\|\bar{u}_yw\|_{x=0},  \|\pa_x^m\bar{v}_{yy}w\|,  \||\bar{v}|\|_{m,w}, [[[\Theta]]]_{m+1,w}\big).
\end{align}
\end{lemma}

\noindent{\bf Proof.}  We square on both side of \eqref{6.234} with $k=m+1$  to obtain
\begin{align}\label{6.253}
	&3\|\bar{u}^2_{\theta} \pa_x S_k w\|^2 + \f78 \kappa^2 \|\bar{u}_{\theta}\pa_y^2 S_{k}w\|^2  + 2\kappa \|\bar{u}_{\theta}^{\f32} \pa_yS_k\|^2_{x=x_0} \nonumber\\
	&\lesssim  \|\bar{u}_{\theta}^{\f32} \pa_yS_k\|^2_{x=0} +  P\big(\|\bar{u}_yw\|_{x=0},  \|\pa_x^m\bar{v}_{yy}w\|,  \||\bar{v}|\|_{m,w}, [[[\Theta]]]_{m+1,w}\big).
\end{align}
where we have used \eqref{6.190}, \eqref{6.149}-\eqref{6.197}, and the following 
\begin{align}\label{6.154}
	\|\bar{u}_{\theta}\big(\Theta_y\pa_x^{k-1}q\big)_{yy}w\| 
	&\lesssim \|(\Theta_y,\Theta_{yy})\|_{L^\infty} \|q_y\|_{\MX^{k-1}_{\t, w}}   + \|\Theta_{yyy}w\|_{L^\infty_xL^2_y} \|\pa_x^{k-1}q\|_{L^2_xL^\infty_y}\nonumber\\
	&\lesssim P\Big(\|\Theta_{yyy}w\|_{L^\infty_xL^2_y} ,  \||\bar{v}|\|_{m,w}\Big),
\end{align}
and 
\begin{align}\label{6.199}
	-\big\la \bar{u}_{\theta}^3 \pa_x S_k,\, \pa_y^2 S_{k} w^2\big\ra 
	&=\big\la \bar{u}_{\theta}^3 \pa_{xy} S_k,\, \pa_y S_{k} w^2\big\ra + \big\la  (\bar{u}_{\theta}^3w^2)_y \pa_x S_k,\, \pa_y S_{k} \big\ra\nonumber\\
	&=\f12 \|\bar{u}_{\theta}^{\f32} \pa_yS_k w\|^2_{x=x_0} - \f12 \|\bar{u}_{\theta}^{\f32}\pa_yS_kw\|^2_{x=0} - \f32 \big\la \bar{u}_{\theta}^2 \bar{u}_x  ,\, |\pa_y S_{k}|^2 w^2\big\ra\nonumber\\
	&\quad  +  \big\la [3\bar{u}_{\theta}^2 \bar{u}_y w^2 + 2\bar{u}_{\theta}^3 w w_y] \pa_x S_k,\, \pa_y S_{k}\big\ra\nonumber\\
	&\geq \f12 \|\bar{u}_{\theta}^{\f32} \pa_yS_k\|^2 \Big|_{x=0}^{x=x_0}  -\f1{8\kappa} \|\bar{u}_{\theta}^2\pa_xS_kw\|^2 -  P\big(\|(\bar{u}_x, \bar{u}_y)\|_{L^\infty}\big) \|\pa_y S_{k} w\|^2 .
\end{align}
Therefore the proof of Lemma \ref{lem6.15} is completed. $\hfill\Box$

\begin{lemma}\label{lem6.15-0}
It holds that 
\begin{align}\label{6.194}
&\sum_{i=0}^{1+\f{m}{2}} \|(\pa_x^i\Theta,\pa_x^i\Theta_y)\la y\ra^{l_0}\|_{L^\infty} + \sum_{i=0}^{\f{m}{2}}\|\pa_x^i\Theta_{yy}\la y\ra^{l_0}\|_{L^\infty} + \|(\Theta_y,  \Theta_{yy},  \Theta_{yyy}, \Theta_{xy}, \Theta_{xxy}) w\|_{L^\infty_xL^2_y}\nonumber\\
&\lesssim  \|(\Theta_y, \Theta_{yy}, \Theta_{xy}, \Theta_{xxy}) w\|_{x=0}   + \sum_{i=0}^{1+\f{m}{2}} \|\pa_x^i\Theta_{yyy}w\|_{x=0} + \sqrt{L} \|\Theta\|_{\MY^{m}_w}  + \sqrt{L} \|\Theta_{xyy} w\|  \nonumber\\
&\quad + \sqrt{L} \sum_{i=0}^{2+\f{m}{2}} \|\pa_x^i\Theta_{yyy} w\|.
\end{align}
\end{lemma}

\noindent{\bf Proof.} A direct calculation shows that 
\begin{align*}
\sum_{i=0}^{k} \|(\pa_x^i\Theta,\pa_x^i\Theta_y)\la y\ra^{l_0}\|_{L^\infty}&\lesssim \sum_{i=0}^{k} \|\pa_x^i\Theta_{y}\la y\ra^{l_0+1}\|_{L^2_y}
\lesssim \sum_{i=0}^{k} \|\pa_x^i\Theta_{y} \la y\ra^{l_0+2}\|_{x=0} + \sqrt{L}\|\Theta \|_{\MY^{k+1}_w},
\end{align*}
and
\begin{align*}
\sum_{i=0}^{k} \|\pa_x^i \Theta_{yy} \la y\ra^{l_0} \|_{L^\infty}& \lesssim \sum_{i=0}^{k} \|\pa_x^i \Theta_{yyy}\la y\ra^{l_0+1}\|_{x=0} + \sqrt{L} \sum_{i=0}^{k+1} \|\pa_x^i \Theta_{yyy}\la y\ra^{l_0+1}\|. 
\end{align*}
Also it holds that 
\begin{align*}
&\sum_{i=1}^3\|\pa_y^i\Theta \, w\|_{L^\infty_xL^2_y} + \sum_{i=1}^2\|\pa_x^i\Theta_{y}\, w\|_{L^\infty_xL^2_y} \nonumber\\
&\lesssim \sum_{i=1}^3\|\pa_y^i\Theta \, w\|_{x=0} + \sum_{i=1}^2\|\pa_x^i\Theta_{y}\, w\|_{x=0}  + \sqrt{L}\sum_{i=1}^3\| \pa_y^i\Theta_{x} w\|  + \sqrt{L} \sum_{i=2}^3\|\pa_x^i\Theta_{y}\, w\| \nonumber\\
&\lesssim  \sum_{i=1}^3\|\pa_y^i\Theta \, w\|_{x=0} + \sum_{i=1}^2\|\pa_x^i\Theta_{y}\, w\|_{x=0}  + \sqrt{L} \|\Theta\|_{\MY^{m}_w}  + \sqrt{L}\sum_{i=2}^3\| \pa_y^i\Theta_{x} w\|.
\end{align*}
Combining above three estimates, we obtain
\begin{align*} 
&\sum_{i=0}^{1+\f{m}{2}} \|(\pa_x^i\Theta,\pa_x^i\Theta_y)\la y\ra^{l_0}\|_{L^\infty}  + \sum_{i=0}^{\f{m}{2}}\|\pa_x^i\Theta_{yy}\la y\ra^{l_0}\|_{L^\infty} + \|(\Theta_y,  \Theta_{yy},  \Theta_{yyy}, \Theta_{xy}, \Theta_{xxy}) w\|_{L^\infty_xL^2_y}\nonumber\\
&\lesssim   \|(\Theta_y, \Theta_{yy}, \Theta_{xy}, \Theta_{xxy}) w\|_{x=0}   + \sum_{i=0}^{1+\f{m}{2}} \|\pa_x^i\Theta_{yyy}w\|_{x=0} + \sqrt{L} \|\Theta\|_{\MY^{m}_w}   + \sqrt{L} \|\Theta_{xyy} w\| \nonumber\\
&\quad + \sqrt{L} \sum_{i=0}^{2+\f{m}{2}} \|\pa_x^i\Theta_{yyy} w\|,
\end{align*}
where we have used Hardy's inequality. Therefore the proof of Lemma \ref{lem6.15-0} is completed. $\hfill\Box$

\medskip

\subsection{Estimates on $F$} 
It is clear that 
\begin{align}\label{6.196}
	\begin{split}
		\sum_{i=0}^{m} \|\pa_x^i\Theta w\| &\lesssim L^{\f16-} \|\Theta \|_{\MY^{m}_w}\quad\mbox{and} \quad 
		\sum_{i=0}^{m} \|\pa_x^i\Theta_yw\|  \lesssim \|\Theta\|_{\MY^{m}_w}.
	\end{split}
\end{align}
which, together with \eqref{6.188}-\eqref{6.188-1}, yields  that
\begin{align}\label{6.197-0}
\sum_{k=0}^{m+1} \|(\pa_x^k \Theta, \pa_x^k\Theta_y)w\| &\lesssim P\Big(\|(\Theta_y, \Theta_{yy}) w\|_{L^\infty_xL^2_y}, \|\Theta\|_{\MY^{m}_w},   \|S_{m+1}\|_{\MY_{\t, w}},  \|q_y\|_{\MX^m_w}\Big)\nonumber\\
&\lesssim P\big([[[\Theta]]]_{m+1,w},  \||\bar{v}|\|_{m,w}\big).
\end{align}
Noting \eqref{6.231} and \eqref{6.115}, one can obtain
\begin{align}\label{6.198}
\|\bar{u}_{\theta}\pa_x^{m+1}\Theta_{yy}w\| &\lesssim \|\bar{u}_{\theta}\pa_y^2S_{m+1}w\| + \|\Theta_{yyy}w\|_{L^\infty_xL^2_y}\|q_y\|_{\MX^m_w}\nonumber\\
&\lesssim \|\bar{u}_{\theta}\pa_y^2S_{m+1}w\|+  P\big([[[\Theta]]]_{m+1,w},  \||\bar{v}|\|_{m,w}\big).
\end{align}

\medskip

Using \eqref{6.173}, we have 
\begin{align}\label{6.199-0}
\sum_{k=0}^m\|\pa_x^k\Theta_{yy}w\|&\lesssim  P\Big(\|(\Theta_y,\Theta_{yy}) w\|_{L^\infty_xL^2_y}, \sum_{i=0}^{\f{m}{2}} \|(\pa_x^i\Theta,\pa_x^i\Theta_y)\la y\ra^{l_0}\|_{L^\infty} , \sum_{k=0}^{m+1} \|\pa_x^{k}\Theta w\|, \sum_{k=0}^{m} \|\pa_x^{k}\Theta_yw\|\Big)\nonumber\\
&\quad + P\big(\|\bar{u}_yw\|_{x=0}, \||\bar{v}|\|_{m,w}\big) \nonumber\\
&\lesssim P\big(\|\bar{u}_yw\|_{x=0}, [[[\Theta]]]_{m+1,w},  \||\bar{v}|\|_{m,w}\big),
\end{align}
where we have used  \eqref{6.197} and \eqref{6.197-0}.

Also we have from \eqref{6.173} that
\begin{align}\label{6.248}
	\sum_{k=0}^m\|\pa_x^{k}\Theta_{yyy}w\| 
	&\lesssim P\Big(\||\bar{v}|\|_{m, w},  \sum_{i=0}^{m+1}\| (\pa_x^{i}\Theta, \pa_x^{i}\Theta_y) w\|, \sum_{i=0}^{m}\|\pa_x^{i}\Theta_{yy} w\|, \sum_{i=0}^{\f{m}{2}} \|(\pa_x^i\Theta,\pa_x^i\Theta_y)\la y\ra^{l_0}\|_{L^\infty}\Big) \nonumber\\
	&\quad +  P\Big(\|(\Theta_y, \Theta_{yy},  \Theta_{xy})w\|_{L^\infty_xL^2_y},   \sum_{i=0}^{\f{m}{2}}\|\pa_x^i\Theta_{yy}\|_{L^\infty},    \sum_{k=0}^{m-1}\|\pa_x^{i}\bar{v}_{yyy}w \|\Big) \nonumber\\
	&\lesssim P\Big(\|\bar{u}_yw\|_{x=0},  [[[\Theta]]]_{m+1,w},  \||\bar{v}|\|_{m, w} ,  \sum_{i=0}^{m-1}\|\pa_x^{i}F w \|\Big),
\end{align}
where we have used  \eqref{6.197-0}, \eqref{6.199-0} and \eqref{6.126-0}.

\


\subsubsection{Estimates on $\|\pa_x^kFw\|$}
Recall $F=-f_1 + f_2 + f_3$. It follows from \eqref{6.6} that 
\begin{align}\label{6.294}
	\pa_x^kf_1 &\cong  \pa_x^k \big\{\bar{u}_y\pa_x(\f{T_y}{T})\big\} + \pa_x^k\big\{\bar{u} \pa_x\big(\f{T_{yy}}{T}\big)\big\} + \pa_x^k \big\{\bar{u}^3 \pa_x\big(\f{T_y^2}{T}\big)\big\} + \pa_x^k\big\{\bar{u}^2 \bar{u}_y \pa_x T_y\big\} \nonumber\\
	&\quad + \pa_x^k\big\{\bar{u}\bar{u}^2_y\, \pa_x T\big\},
\end{align}
which, together with \eqref{6.197-0} and \eqref{6.199-0},  yields that 
\begin{align}\label{6.309}
\sum_{k=0}^m\|\pa_x^{k}f_1w\| &\lesssim  P\Big( \|\bar{u}_{\theta}\pa_x^{m+1}\Theta_{yy}w\|,  \sum_{i=0}^{m+1}\| (\pa_x^{i}\Theta, \pa_x^{i}\Theta_y) w\|, \sum_{i=0}^{m}\|\pa_x^{i}\Theta_{yy} w\| \Big)\nonumber\\
&\quad +  P\Big(\sum_{i=0}^{\f{m}{2}} \|(\pa_x^i\Theta,\pa_x^i\Theta_y)\la y\ra^{l_0}\|_{L^\infty}, \sum_{i=0}^{\f{m}{2}}\|\pa_x^i\Theta_{yy}\|_{L^\infty}, \||\bar{v}|\|_{m,w}\Big) \nonumber\\
&\lesssim P\big(\|\bar{u}_yw\|_{x=0}, [[[\Theta]]]_{m+1,w},  \||\bar{v}|\|_{m,w}, \|\bar{u}_{\theta}\pa_x^{m+1}\Theta_{yy}w\|\big).
\end{align}

From $\eqref{6.17-0}_1$, it is clear that
\begin{align*}
\pa_x^k f_2= 2\sum_{i=0}^k \pa_x^i(\bar{u}_{\theta} \bar{v}_y)\, \pa_x^{k-i} q_y,
\end{align*}
which yields that 
\begin{align}\label{6.231-0}
\sum_{k=0}^m\|\pa_x^kf_2w\|
&\lesssim P\big(\||\bar{v}|\|_{m,w}, \|\pa_x^{m}\bar{v}_{yy} w\|\big).
\end{align}

\medskip

Using \eqref{6.197-0} and \eqref{6.199-0}, we have from $\eqref{6.17-0}_2$ that
\begin{align}\label{6.233-0}
\sum_{k=0}^m\|\pa_x^k f_3w\|&\lesssim P\Big(\sum_{i=0}^{\f{m}{2}} \|(\pa_x^i\Theta,\pa_x^i\Theta_y)\la y\ra^{l_0}\|_{L^\infty},  \sum_{i=0}^{\f{m}{2}}\|\pa_x^i\Theta_{yy}\|_{L^\infty},  \sum_{i=0}^{m} \|(\pa_x^i\Theta,  \pa_x^i\Theta_y,\pa_x^{i}\Theta_{yy})w\|\Big) \nonumber\\
&\quad + P\Big(\||\bar{v}|\|_{m,w},\|\pa_x^{m}\bar{v}_{yy} w\| \Big) \nonumber\\
&\lesssim  P\Big(\|\bar{u}_yw\|_{x=0},   [[[\Theta]]]_{m+1,w},  \||\bar{v}|\|_{m,w}, \|\pa_x^{m}\bar{v}_{yy} w\|\Big).
\end{align}


\smallskip

Hence it is follows from \eqref{6.309}-\eqref{6.233-0} that
\begin{align}\label{6.210}
	\sum_{i=0}^{m-1}\|\pa_x^i F w\| 
	&\lesssim P\big(\|\bar{u}_yw\|_{x=0}, [[[\Theta]]]_{m+1,w},  \||\bar{v}|\|_{m,w},\,  \|\bar{u}_{\theta}\pa_x^{m}\Theta_{yy}w\|\big)\nonumber\\
	&\lesssim P\big(\|\bar{u}_yw\|_{x=0}, [[[\Theta]]]_{m+1,w},  \||\bar{v}|\|_{m,w}\big).
\end{align}

Also we have from  \eqref{6.309}, \eqref{6.231-0}  and \eqref{6.233-0} that
\begin{align}\label{6.209}
\sum_{i=0}^{m}\|\pa_x^i F w\| 
&\lesssim P\big(\|\bar{u}_yw\|_{x=0}, [[[\Theta]]]_{m+1,w},  \||\bar{v}|\|_{m,w},\,  \|\bar{u}_{\theta}\pa_x^{m+1}\Theta_{yy}w\|, \|\pa_x^{m}\bar{v}_{yy} w\|\big) \nonumber\\
&\lesssim P\big(\|(q_y, \bar{u}_{y}, \bar{u}_{yy})w\|_{x=0}, [[[\Theta]]]_{m+1,w},  \||\bar{v}|\|_{m,w},\,  \|\bar{u}_{\theta}\pa_y^2S_{m+1}w\|\big),
\end{align}
where we have used \eqref{6.283-000}, \eqref{6.210} and \eqref{6.198}.

\smallskip

\subsubsection{Estimates on $\|\pa_x^k\pa_yFw\|$}  A direct calculation shows that 
\begin{align}\label{6.211}
\sum_{i=0}^{m-1}\|\pa_x^{i}\pa_y f_2w\|
&\lesssim P\Big(\||\bar{v}|\|_{m, w},  \sum_{i=0}^{m-1}\|\pa_x^i \bar{v}_{yyy}w\|\Big)\lesssim P\Big(\||\bar{v}|\|_{m, w},  \sum_{i=0}^{m-1}\|\pa_x^{i} F w\|\Big),
\end{align}
\begin{align}\label{6.212}
\sum_{i=0}^{m-1}\|\pa_x^{i}\pa_y f_1w\| &\lesssim   P\Big(\sum_{i=0}^{\f{m}{2}} \|(\pa_x^i\Theta,\pa_x^i\Theta_y)\la y\ra^{l_0}\|_{L^\infty}, \sum_{i=0}^{m}\|(\pa_x^i\Theta,  \pa_x^i\Theta_y,\pa_x^{i}\Theta_{yy})w\|, \sum_{i=0}^{m}\|\pa_x^i \Theta_{yyy}w\|\Big) \nonumber\\
&\quad + P\Big(  \sum_{i=0}^{\f{m}{2}}\|\pa_x^i\Theta_{yy}\|_{L^\infty}, \||\bar{v}|\|_{m, w},  \sum_{i=1}^{m-2} \|\pa_x^{i}\bar{v}_{yyy} w\| \Big)\nonumber\\
&\lesssim P\Big(\|\bar{u}_yw\|_{x=0},  [[[\Theta]]]_{m+1,w},  \||\bar{v}|\|_{m, w} ,  \sum_{i=0}^{m-1}\|\pa_x^{i}F w \|\Big),
\end{align}
and
\begin{align}\label{6.213}
\sum_{i=0}^{m-1}\|\pa_x^{i}\pa_y f_3w\|
&\lesssim  P\Big(\sum_{i=0}^{m-1}\|(\pa_x^i\Theta,  \pa_x^i\Theta_y,\pa_x^{i}\Theta_{yy}, \pa_x^i \Theta_{yyy})w\| , \||\bar{v}|\|_{m, w},  \|\pa_x^{m-1} \bar{v}_{yyy}w\|\Big)\nonumber\\
&\lesssim P\Big(\|\bar{u}_yw\|_{x=0},  [[[\Theta]]]_{m+1,w},  \||\bar{v}|\|_{m, w} ,  \sum_{i=0}^{m-1}\|\pa_x^{i}F w \|\Big),
\end{align}
where we have used \eqref{6.126-0}, \eqref{6.197-0}, \eqref{6.199-0} and  \eqref{6.248}.

Thus, combining \eqref{6.211}-\eqref{6.213},  then using \eqref{6.210},  one gets
\begin{align}\label{6.246-0}
\sum_{i=0}^{m-1}\|\pa_x^{i}\pa_y F w\| &\lesssim   P\big(\|\bar{u}_yw\|_{x=0},  [[[\Theta]]]_{m+1,w},  \||\bar{v}|\|_{m, w}\big).
\end{align}

\subsection{Conclusion} 
\begin{proposition}\label{prop2.16}
Let  the initial data $(\mathscr{U}_0, \mathscr{T}_0)$ satisfy the conditions \eqref{0.1}-\eqref{0.2}, and also generic compatibility conditions at the corner up to order $2\mathfrak{m}_0+6$ so that \eqref{2.12} is well-defined. Let $m\leq \mathfrak{m}_0+1$. There is a small $L>0$ $($independent of $\theta\in(0,1]\, )$ such that there exists a unique solution to \eqref{3.58} in $x\in [0,L]$ with the following estimates
\begin{align}\label{6.195}
&\||\bar{v}|\|_{m, w} + [[[\Theta]]]_{m+1,w} + \|\bar{u}^{\f32} \pa_yS_{m+1}\|_{x=x_0} + \|\bar{u}^2 \pa_x S_{m+1} w\| +  \|\bar{u}\pa_y^2 S_{m+1}w\|\nonumber\\
&\quad + \sum_{k=0}^{m} \Big\{\|\bar{u}^{\f52}\pa_x^k q_{yy}w\|^2_{x=x_0} +  \|\bar{u}^{3}\pa_x^{k+1} q_yw\|^2 + \|\bar{u} \pa_x^k\bar{v}_{yyy}w\|^2 \Big\}
 \lesssim C\big({\bf I}[\mathscr{T}_b, \mathscr{U}_0, \mathscr{T}_0]\big),
\end{align}
\end{proposition}

\noindent{\bf Proof.}  Let $m\leq \mathfrak{m}_0+1$. Recall \eqref{0.2}, it holds that
\begin{equation*}
IN+ \|\bar{u}_yw\|_{x=0} + \sum_{i=0}^{m-1}\|\pa_x^{i}\bar{v}_{yy}w\|_{x=0} + \|\bar{u}_{\theta}^{\f32} \pa_yS_{m+1}\|_{x=0} + \sum_{k=0}^{m}  \|\bar{u}_{\theta}^{\f52}\pa_x^k q_{yy}w\|^2_{x=0}\lesssim {\bf I}[\mathscr{T}_b, \mathscr{U}_0, \mathscr{T}_0],
\end{equation*}
where  ${\bf I}[\mathscr{T}_b, \mathscr{U}_0, \mathscr{T}_0]$ is  independent of $\theta\in (0,1]$.  We divide the proof into two steps. 

\noindent{\it Step 1. Existence and uniform estimates of $\theta$-approximate problem \eqref{5.3}.} Let $(\bar{u}, \bar{v}, T)$ be the local smooth solution of \eqref{5.3}. We point out that the solution must depend on $\theta>0$, but we omit such dependence for simplicity. Then it follows from \eqref{6.139}, \eqref{6.283-00} and \eqref{6.283-0} that 
\begin{align}\label{6.218}
 \||\bar{v}|\|_{m, w}
&\lesssim C\big({\bf I}[\mathscr{T}_b, \mathscr{U}_0, \mathscr{T}_0]\big)   + L^{\f16-}P\Big(\||\bar{v}|\|_{m,w} ,  \sum_{i=0}^{m}\|\pa_x^{i}F w\|, \sum_{i=0}^{m-1} \|\pa_x^{i} F_{y} w\|\Big).
\end{align}

\smallskip

We have from \eqref{6.169-0}, \eqref{6.173-0} and \eqref{6.194} that
\begin{align}\label{6.219-0}
[[[\Theta]]]_{m+1,w}
&\lesssim  C({\bf I}[\mathscr{T}_b, \mathscr{U}_0, \mathscr{T}_0]) + L^{\f16-}P\Big([[[\Theta]]]_{m+1,w},  \||\bar{v}|\|_{m, w}, \|\pa_x^m\bar{v}_{yy}w\|\Big) \nonumber\\
&\quad  + \sqrt{L} \|\Theta_{xyy} w\| + \sqrt{L} \sum_{i=0}^{2+\f{m}{2}} \|\pa_x^i\Theta_{yyy} w\|\nonumber\\
&\lesssim   C({\bf I}[\mathscr{T}_b, \mathscr{U}_0, \mathscr{T}_0]) + L^{\f16-}P\Big([[[\Theta]]]_{m+1,w},  \||\bar{v}|\|_{m, w}, \sum_{i=0}^{m-1}\|\pa_x^{i}F w \|\Big),
\end{align}
where we have used \eqref{6.283-000}, \eqref{6.199-0}, \eqref{6.248} and the following fact
\begin{align}
&\sum_{i=0}^{m} \|\sqrt{\bar{u}_{\theta}}\pa_x^i\Theta w\|^2_{x=0} + \|\sqrt{\bar{u}_{\theta}}S_{m+1} w\|^2_{x=0}  + \sum_{i=0}^{1+\f{m}{2}} \|\pa_x^i\Theta_{yyy}w\|_{x=0} \nonumber\\
&\quad   + \|(\Theta_y, \Theta_{yy}, \Theta_{xy}, \Theta_{xxy}) w\|_{x=0} \lesssim C\big({\bf I}[\mathscr{T}_b, \mathscr{U}_0, \mathscr{T}_0]\big). 
\end{align}
Also we have from \eqref{6.189}  that
\begin{align}\label{6.252}
& \|\bar{u}_{\theta}^{\f32} \pa_yS_{m+1}\|^2_{x=x_0} + \|\bar{u}^2_{\theta} \pa_x S_{m+1} w\|^2 +  \|\bar{u}_{\theta}\pa_y^2 S_{m+1}w\|^2  \nonumber\\
&\lesssim \|\bar{u}_{\theta}^{\f32} \pa_yS_{m+1}\|^2_{x=0} +  P\big({\bf I}[\mathscr{T}_b, \mathscr{U}_0, \mathscr{T}_0], [[[\Theta]]]_{m+1,w},  \||\bar{v}|\|_{m, w}\big),
\end{align}
where we have used \eqref{6.283-000} and \eqref{6.210}.

\smallskip

Combining \eqref{6.218} and \eqref{6.219-0}, then using \eqref{6.209} and \eqref{6.246-0}, one gets that
\begin{align}\label{6.185}
&\||\bar{v}|\|_{m, w} + [[[\Theta]]]_{m+1,w} \nonumber\\
&\lesssim C({\bf I}[\mathscr{T}_b, \mathscr{U}_0, \mathscr{T}_0]) + L^{\f16-}P\Big([[[\Theta]]]_{m+1,w}, \||\bar{v}|\|_{m,w} ,  \sum_{i=0}^{m}\|\pa_x^{i}F w\|, \sum_{i=0}^{m-1} \|\pa_x^{i} F_{y} w\|\Big) \nonumber\\
&\lesssim C({\bf I}[\mathscr{T}_b, \mathscr{U}_0, \mathscr{T}_0]) +   L^{\f16-} P\big(  [[[\Theta]]]_{m+1,w},  \||\bar{v}|\|_{m,w},\,  \|\bar{u}_{\theta}\pa_y^2S_{m+1}w\|\big),
\end{align}

Substituting \eqref{6.252} into \eqref{6.185}, one proves 
\begin{align}\label{6.186}
\||\bar{v}|\|_{m, w} + [[[\Theta]]]_{m+1,w} 
&\lesssim C({\bf I}[\mathscr{T}_b, \mathscr{U}_0, \mathscr{T}_0]) +   L^{\f16-} P\big(  [[[\Theta]]]_{m+1,w},  \||\bar{v}|\|_{m,w}\big).
\end{align}
Based on \eqref{6.186}, the local existence of smooth solution of \eqref{5.3}  established in section \ref{sec2.1} and continuity argument,  there is a small $L>0$ (independent of $\theta>0$) such that the solution of \eqref{5.3} exists in $x\in[0,L]$ with the following uniform estimate
\begin{align}\label{6.186-2}
\||\bar{v}|\|_{m, w} + [[[\Theta]]]_{m+1,w} 
&\lesssim C({\bf I}[\mathscr{T}_b, \mathscr{U}_0, \mathscr{T}_0]).
\end{align}
Then it follows from \eqref{6.252}, \eqref{6.186-2} and  \eqref{6.266} that
\begin{align}\label{6.186-3}
& \|\bar{u}_{\theta}^{\f32} \pa_yS_{m+1}\|^2_{x=x_0} + \|\bar{u}^2_{\theta} \pa_x S_{m+1} w\|^2 +  \|\bar{u}_{\theta}\pa_y^2 S_{m+1}w\|^2 
\lesssim C({\bf I}[\mathscr{T}_b, \mathscr{U}_0, \mathscr{T}_0]),
\end{align}
and
\begin{equation}\label{6.206-1}
\sum_{k=0}^{m} \Big\{\|\bar{u}_{\theta}^{\f52}\pa_x^k q_{yy}w\|^2_{x=x_0} +  \|\bar{u}_{\theta}^{3}\pa_x^{k+1} q_yw\|^2 + \|\bar{u}_{\theta} \pa_x^k\bar{v}_{yyy}w\|^2 \Big\} 
 \lesssim C({\bf I}[\mathscr{T}_b, \mathscr{U}_0, \mathscr{T}_0]).
\end{equation}

\noindent{\it Step 2. Existence and uniform estimate of \eqref{3.58}.} With the $\theta$-approximate solutions established in Step 1, by noting its uniform estimates \eqref{6.186-2}-\eqref{6.186-3}, we can take the limit $\theta\to0+$ to establish the existence of solution of \eqref{3.58} in $x\in [0,L]$. The uniform estimate \eqref{6.195} follows directly from \eqref{6.186-2}-\eqref{6.186-3}. Therefore the proof of Proposition \ref{prop2.16} is completed.

\medskip

\noindent{\bf Proof of Theorem \ref{thm1}:} Noting \eqref{6.195}, one has from \eqref{6.126-0}, \eqref{6.283-000} and  \eqref{6.197-0}-\eqref{6.199-0} that 
\begin{align}\label{6.186-4}
	& \|\bar{u} \pa_x^{m+1}\Theta_{yy}w\| +  \sum_{i=0}^{m+1} \|(\pa_x^i \Theta, \pa_x^i\Theta_y)w\| + \sum_{i=0}^m\|\pa_x^i\Theta_{yy}w\| \nonumber\\
	& + \sum_{i=0}^{m} \|(\pa_x^{i}\bar{v}_y, \pa_x^{i}\bar{v}_{yy})w\| + \sum_{i=0}^{m-1}\|\pa_x^i\bar{v}_{yyy} w\| \lesssim C({\bf I}[\mathscr{T}_b, \mathscr{U}_0, \mathscr{T}_0]),
\end{align}
which yields immediately that 
\begin{align}
\sum_{i=0}^{m} \|(\pa_x^{i+1}\bar{u}, \pa_x^{i+1}\bar{u}_y)w\| +  \sum_{i=1}^{m}\|\pa_x^i\bar{u}_{yy} w\|
 \lesssim  C({\bf I}[\mathscr{T}_b, \mathscr{U}_0, \mathscr{T}_0]).
\end{align}
Also it is clear that 
\begin{align}\label{6.186-8}
\|(\pa_{y}\bar{u}, \pa^2_{y}\bar{u})w \| 
&\lesssim  \sqrt{L} \|(\pa_{y}\bar{u}_0, \pa^2_{y}\bar{u}_0)w \| + L \|(\bar{u}_{xy}, \bar{u}_{xyy})w\| \lesssim \sqrt{L} C({\bf I}[\mathscr{T}_b, \mathscr{U}_0, \mathscr{T}_0]).
\end{align}

Combining the equations \eqref{5.3-0}, \eqref{6.195}, \eqref{6.186-4}-\eqref{6.186-8} and induction argument, one can obtains
\begin{align}\label{6.186-9}
\sum_{0\leq 2i+j\leq 2m+2} \|\pa_x^i \pa_y^{j}\Theta w\| + \sum_{0<2i+j\leq 2m+2} \|\pa_x^i \pa_y^{j}\bar{u} w\| \lesssim C({\bf I}[\mathscr{T}_b, \mathscr{U}_0, \mathscr{T}_0]).
\end{align}
Noting
\begin{align*}
\|\pa_x^i \pa_y^{j}\bar{u} w\la y\ra^{-1}\|_{L^\infty} &\lesssim\|\pa_x^i \pa_y^{j}\bar{u}(0,\cdot) w\la y\ra^{-1}\|_{L^\infty_y} + \sqrt{L} \|\pa_x^{i+1} \pa_y^{j+1}\bar{u} w \|,\\
\|\pa_x^i \pa_y^{j}\Theta w\la y\ra^{-1}\|_{L^\infty} &\lesssim\|\pa_x^i \pa_y^{j}\Theta(0,\cdot) w\la y\ra^{-1}\|_{L^\infty_y} + \sqrt{L} \|\pa_x^{i+1} \pa_y^{j+1}\Theta w \| ,
\end{align*}
which, together with \eqref{6.186-9}, yields that 
\begin{align}\label{211}
\begin{split}
\sum_{0< 2i+j\leq 2m-1}\|\pa_x^i \pa_y^{j}\bar{u}  w\la y\ra^{-1}\|_{L^\infty} 
&\lesssim \sum_{0< 2i+j\leq 2m-1}\|\pa_x^i \pa_y^{j}\bar{u}(0,\cdot) w\la y\ra^{-1}\|_{L^\infty_y} + \sqrt{L} C({\bf I}[\mathscr{T}_b, \mathscr{U}_0, \mathscr{T}_0]),\\
\sum_{0\leq 2i+j\leq 2m-1}\|\pa_x^i \pa_y^{j}\Theta  w\la y\ra^{-1}\|_{L^\infty} 
&\lesssim \sum_{0\leq 2i+j\leq 2m-1}\|\pa_x^i \pa_y^{j}\Theta(0,\cdot) w\la y\ra^{-1}\|_{L^\infty_y} + \sqrt{L}C({\bf I}[\mathscr{T}_b, \mathscr{U}_0, \mathscr{T}_0]).
\end{split}
\end{align}
Also it holds that
\begin{align}\label{212}
\|\pa_x^{m-1}\bar{v}\|_{L^\infty} \lesssim  \|\pa_x^{m-1}\bar{v}(0,\cdot)\|_{L^\infty_y} + \sqrt{L}  \|\pa_x^{m}\bar{v}_y \la y\ra\|\lesssim  \|\pa_x^{m-1}\bar{v}(0,\cdot)\|_{L^\infty_y} + \sqrt{L} C({\bf I}[\mathscr{T}_b, \mathscr{U}_0, \mathscr{T}_0]).
\end{align}
Since $m\leq \mathfrak{m}_0+1$,  we conclude \eqref{0.4} from \eqref{211}-\eqref{212}. Then \eqref{0.3} follows from \eqref{0.4}, \eqref{3.70} and \eqref{6.11-1}. Therefore the proof of Theorem \ref{thm1} is completed.  $\hfill\Box$

 \medskip

\section{Existence for the Steady Linear Compressible Prandtl Layer}\label{sec3}
This section is devoted to proving the existence of  local-in-$x$ solution for the steady linear compressible Prandtl layer equations \eqref{3.6-10}--\eqref{3.6-11}. Recall the smooth solution $(\bar{\rho}^0_p, \bar{u}^0_p, \bar{v}^0_p, \bar{T}^0_p)$ of the nonlinear Prandtl equations \eqref{3.6-00} established in Theorem \ref{thm1}.
\subsection{Reformulation} 
We define 
\begin{align}\label{3.1}
P_{1}:=-\int_y^\infty G_{1}(x,z) dz \quad \& \quad 
P_{2}:= \int_0^y G_{2}(x,z) dz, 
\end{align}
and 
\begin{align}\label{3.2}
	\begin{cases}
		\bar{\bfu}:=\bar{\rho}^0_p u_p + \rho_p \bar{u}^0_p,\\
		\bar{\bfv}:=\bar{\rho}^0_p v_p + \rho_p \bar{v}^0_p  - P_2,
	\end{cases}
 \Longleftrightarrow  \quad 
	\begin{cases}
		u_p =  \bar{T}_p^0 \bar{\bfu} - \rho_p \bar{T}_p^0\bar{u}^0_p,\\
	 v_p   =  \bar{T}_p^0\bar{\bfv} -  \rho_p \bar{T}_p^0\bar{v}^0_p  + \bar{T}_p^0P_2.
	\end{cases}
\end{align}
where we have used the fact $\bar{T}_p^0\bar{\rho}_p^0=1$.

\begin{lemma}\label{lem4.1}
The system \eqref{3.6-10} can be rewritten as the following equivalent system:
\begin{align}\label{3.73}
	\begin{cases}
		\bar{\bfu}_x + \bar{\bfv}_y=0,\\[1mm]
		\bar{u}^0_{p} \pa_x \bar{\bfu}  + [\pa_y \bar{u}^{0}_p - \bar{u}^0_p  \f{\pa_y\bar{T}^0_p}{\bar{T}^0_p}]\, \bar{\bfv}  - \mu \bar{T}_p^0 \pa_{yy}\bar{\bfu} = Q_1(\bar{\bfu}, T_p, G),\\[1mm]
		2\bar{\rho}^0_p \bar{u}^0_p \pa_xT_p + 2\bar{\rho}^0_p  \bar{v}^0_p \pa_yT_p + 2\pa_x\bar{T}^0_p\, \bar{\bfu} + 2\pa_y\bar{T}^0_p\, \bar{\bfv} 
		 = \kappa \pa_{yy}T_p  + H_1(\bar{\bfu}, T_p, G),
	\end{cases}
\end{align}
where 
\begin{align}
H_1(\bar{\bfu}, T_p, G):&= 2 \mu \pa_y \bar{u}^{0}_p \big[\bar{T}_p^0 \pa_y\bar{\bfu} + \pa_y\bar{T}_p^0 \, \bar{\bfu}\big]  + 2 \mu \pa_y \bar{u}^{0}_p \big[ \bar{\rho}^0_p\bar{u}^0_p \, \pa_yT_p + \pa_y(\bar{\rho}^0_p\bar{u}^0_p)\, T_p\big]\nonumber\\
&\quad + G_{4} - 2 \mu \pa_y \bar{u}^{0}_p \pa_y(P_1\bar{u}^0_p)   +  (\bar{U}^0_p\cdot \bar{\nabla})P_1  - \pa_y(P_2 \bar{T}^0_p) -P_2\,\pa_y\bar{T}^0_p ,\\[1mm]
Q_1(\bar{\bfu}, T_p, G):&=\big[2\mu \pa_y\bar{T}_p^0 - \bar{v}^0_p\big] \pa_{y}\bar{\bfu}  + \big[\mu \pa_{yy}\bar{T}_p^0  -  \bar{\rho}^0_p  \bar{v}^0_p  \pa_y\bar{T}_p^0 - \pa_x \bar{u}^{0}_p  \big] \bar{\bfu}  - \f12 \bar{\rho}^0_p  \bar{u}^0_p H_1(\bar{\bfu}, T_p) \nonumber\\
&\quad -  \f12 \kappa \bar{\rho}^0_p  \bar{u}^0_p\pa_{yy}T_p  + \mu \pa_{yy}\big(\bar{\rho}_p^0\bar{u}^0_p  T_p\big)  - \bar{\rho}^0_p  (\bar{U}^0_p\cdot \bar{\nabla}) (\bar{\rho}_p^0\bar{u}^0_p)\, T_p   \nonumber\\
&\quad  + \bar{\rho}^0_p  (\bar{U}^0_p\cdot \bar{\nabla}) (P_1\bar{u}^0_p)     + G_{3} -  \pa_y u^{0}_p  P_2  - P_{1x}  - \mu \pa_{yy}\big(\bar{u}^0_p P_1\big).
\end{align}
The boundary conditions  are 
\begin{align}\label{4.15E}
	\begin{cases}
		\dis (\bar{\bfu}, T_p)|_{x=0}=\big(\bar{\rho}^0_p \tilde{u}_0 + \bar{u}^0_p [\bar{\rho}^0_p P_1 - (\bar{\rho}^0_p)^2 \tilde{T}_0 ],\,  \tilde{T}_0\big)(0,y),\\[1mm]
		\dis (\bar{\bfu}, \bar{\bfv})|_{y=0} =\big(\bar{\rho}^0_p(x,0)\, \tilde{u}_b(x),\, 0\big),\quad \dis \lim_{y\to \infty} (\bar{\bfu}, T_p)=(0,0),\\[1mm]
		\dis \pa_yT_p|_{y=0}= \tilde{T}_b(x)\quad \mbox{for NBC},
		\\[1mm]
		\dis T_p|_{y=0}=\tilde{T}_b(x)\qquad \mbox{for DBC}. 
	\end{cases}
\end{align}
\end{lemma}

\noindent{\bf Proof.} 
1. It follows from $\eqref{3.6-10}_{1,2}$ that
\begin{align*}
	\bar{\rho}^0_p T_p + \bar{T}^0_p \rho_p =  P_1\quad \Longrightarrow\quad 
	\rho_{p}=\f{1}{\bar{T}^0_p} \big[P_1 - \bar{\rho}^0_p T_p\big] 
	=\bar{\rho}^0_p P_1 - (\bar{\rho}^0_p)^2 T_p ,
\end{align*}
and
\begin{align*}
	\bar{\bfu}_x + \bar{\bfv}_y=0.
\end{align*}

2. We  rewrite the heat equation $\eqref{3.6-10}_3$ as
\begin{align*}
	&\bar{\rho}^0_p (\bar{U}^0_p \cdot \bar{\nabla})T_p + \bar{\rho}^0_p\big\{u_p\pa_xT^0_p + v_p \pa_yT^0_p\big\} + (\bar{U}^0_p\cdot \bar{\nabla}) T^{0}_p \, \rho_p   +   [\pa_x u_p + \pa_y v_p]  \nonumber\\
	&= \kappa \pa_{yy}T_p  + 2 \mu \pa_y \bar{u}^{0}_p \cdot \pa_y u_p + G_{4} - P_1 \mbox{div} \bar{U}^0_p .
\end{align*}
A direct calculation shows that 
\begin{align}\label{3.71}
	\begin{split}
		&\bar{\rho}^0_p\big\{u_p\pa_xT^0_p +  v_p \pa_yT^0_p\big\}
		= \pa_x\bar{T}^0_p\, \bar{\bfu} + \pa_y\bar{T}^0_p\, \bar{\bfv} - (\bar{U}^0_p\cdot \bar{\nabla}) \bar{T}^{0}_p \, \rho_p + P_2\,\pa_y\bar{T}^0_p,\\
		&\pa_x u_p + \pa_y v_p
		= \pa_x\bar{T}^0_p\, \bar{\bfu} + \pa_y\bar{T}^0_p\, \bar{\bfv} + \bar{\rho}^0_p (\bar{U}^0_p \cdot \bar{\nabla})T_p - \mbox{div}(P_1\bar{U}^0_p) + \pa_y(P_2 \bar{T}^0_p),\\
		&\pa_y u_p = \bar{T}_p^0 \pa_y\bar{\bfu} + \pa_y\bar{T}_p^0 \, \bar{\bfu} + \bar{\rho}^0_p\bar{u}^0_p \, \pa_yT_p + \pa_y(\bar{\rho}^0_p\bar{u}^0_p)\, T_p- \pa_y(P_1\bar{u}^0_p).
	\end{split}
\end{align}
Substituting \eqref{3.71} into $\eqref{3.6-10}_4$, we have 
\begin{align*}
2\bar{\rho}^0_p \bar{u}^0_p \pa_xT_p + 2\bar{\rho}^0_p  \bar{v}^0_p \pa_yT_p + 2\pa_x\bar{T}^0_p\, \bar{\bfu} + 2\pa_y\bar{T}^0_p\, \bar{\bfv}  
&= \kappa \pa_{yy}T_p  + H_1(\bar{\bfu}, T_p),
\end{align*}
which concludes $\eqref{3.73}_3$.

3. We rewrite the momentum equation $\eqref{3.6-10}_2$ as
\begin{align*}
	\bar{\rho}^0_p  (\bar{U}^0_p\cdot \bar{\nabla}) u_p + \bar{\rho}^0_p \big\{ u_p  \pa_xu^{0}_p + v_p  \pa_yu^0_p\big\}
	=\mu \pa_{yy} u_p - \rho_p \big(\bar{U}^0_p\cdot \bar{\nabla}\big) \bar{u}^{0}_p  + G_{3} - P_{1x}.
\end{align*}
It is noted that 
\begin{align}\label{3.66}
	\bar{\rho}^0_p  (\bar{U}^0_p\cdot \bar{\nabla}) u_p 
	&=(\bar{U}^0_p\cdot \bar{\nabla}) \bar{\bfu} + \bar{\rho}^0_p  (\bar{U}^0_p\cdot \bar{\nabla}) \bar{T}_p^0 \, \bar{\bfu}  + |\bar{\rho}^0_p|^2 \bar{u}^0_p (\bar{U}^0_p\cdot \bar{\nabla}) T_p
	\nonumber\\
	&\quad + \bar{\rho}^0_p  (\bar{U}^0_p\cdot \bar{\nabla}) (\bar{\rho}_p^0\bar{u}^0_p)\, T_p  - \bar{\rho}^0_p  (\bar{U}^0_p\cdot \bar{\nabla}) (P_1\bar{u}^0_p) ,
\end{align}
\begin{align}
\bar{\rho}^0_p \big\{ u^k_p  \pa_xu^{0}_p + v_p \pa_yu^0_p\big\} 
&=\bar{\rho}^0_p \pa_x \bar{u}^{0}_p [\bar{T}_p^0 \bar{\bfu} - \rho_p \bar{T}_p^0\bar{u}^0_p] + \bar{\rho}^0_p \pa_y \bar{u}^{0}_p \,[\bar{T}_p^0\bar{\bfv} -  \rho_p \bar{T}_p^0\bar{v}^0_p  + \bar{T}_p^0P_2] \nonumber\\
	&=\pa_x \bar{u}^{0}_p \, \bar{\bfu} + \pa_y \bar{u}^{0}_p\, \bar{\bfv} - (\bar{U}^0_p \cdot \bar{\nabla}) \bar{u}^{0}_p  \,  \rho_p   +   \pa_y \bar{u}^{0}_p  P_2,
\end{align}
and
\begin{align}\label{3.68-1}
	\pa_{yy} u_p = \bar{T}_p^0 \pa_{yy}\bar{\bfu} + 2\pa_y\bar{T}_p^0 \pa_{y}\bar{\bfu} + \pa_{yy}\bar{T}_p^0 \bar{\bfu} - \pa_{yy}\big(\bar{T}_p^0\bar{u}^0_p\, \rho_p\big).
\end{align}
Substituting \eqref{3.66}-\eqref{3.68-1} into $\eqref{3.6-10}_2$, one gets
\begin{align*}
	&\bar{u}^0_p \pa_x \bar{\bfu}  + \bar{v}^0_p \pa_y \bar{\bfu} + \pa_x \bar{u}^{0}_p \, \bar{\bfu} + \pa_y \bar{u}^{0}_p\, \bar{\bfv} +  |\bar{\rho}^0_p|^2 \bar{u}^0_p (\bar{U}^0_p\cdot \bar{\nabla}) T_p  \nonumber\\
	&=\mu \bar{T}_p^0 \pa_{yy}\bar{\bfu} + 2\mu \pa_y\bar{T}_p^0 \pa_{y}\bar{\bfu} + \big[\mu \pa_{yy}\bar{T}_p^0  -  \bar{\rho}^0_p  (\bar{U}^0_p\cdot \bar{\nabla}) \bar{T}_p^0 \big] \bar{\bfu}   - \mu \pa_{yy}\big(\rho_p \bar{T}_p^0\bar{u}^0_p\big) \nonumber\\
	&\quad - \bar{\rho}^0_p  (\bar{U}^0_p\cdot \bar{\nabla}) (\bar{\rho}_p^0\bar{u}^0_p)\, T_p  + \bar{\rho}^0_p  (\bar{U}^0_p\cdot \bar{\nabla}) (P_1\bar{u}^0_p)     + G_{3} -  \pa_y u^{0}_p  P_2  - P_{1x},
\end{align*}
which, together with $\eqref{3.73}_3$, yields that 
\begin{align*}
&\bar{u}^0_p \pa_x \bar{\bfu}  + \bar{v}^0_p \pa_y \bar{\bfu} + \pa_x \bar{u}^{0}_p \, \bar{\bfu} + \big[\pa_y \bar{u}^{0}_p  - \bar{\rho}^0_p  \bar{u}^0_p \pa_y\bar{T}^0_p\big] \bar{\bfv} \nonumber\\
&=\mu \bar{T}_p^0 \pa_{yy}\bar{\bfu} + 2\mu \pa_y\bar{T}_p^0 \pa_{y}\bar{\bfu} + \big[\mu \pa_{yy}\bar{T}_p^0  -  \bar{\rho}^0_p  (\bar{U}^0_p\cdot \bar{\nabla}) \bar{T}_p^0 \big] \bar{\bfu} + \bar{\rho}^0_p  \bar{u}^0_p \pa_x\bar{T}^0_p\, \bar{\bfu} \nonumber\\
&\quad -  \f12 \kappa \bar{\rho}^0_p  \bar{u}^0_p\pa_{yy}T_p  - \f12 \bar{\rho}^0_p  \bar{u}^0_p H(\bar{\bfu}, T_p)  + \mu \pa_{yy}\big(\bar{\rho}_p^0\bar{u}^0_p  T_p\big)  - \bar{\rho}^0_p  (\bar{U}^0_p\cdot \bar{\nabla}) (\bar{\rho}_p^0\bar{u}^0_p)\, T_p   \nonumber\\
&\quad  + \bar{\rho}^0_p  (\bar{U}^0_p\cdot \bar{\nabla}) (P_1\bar{u}^0_p)     + G_{3} -  \pa_y u^{0}_p  P_2  - P_{1x}  - \mu \pa_{yy}\big(\bar{u}^0_p P_1\big).
\end{align*}
We then immediately deduce the second equality in \eqref{3.73}. 

The boundary conditions \eqref{4.15E} follows directly from \eqref{3.6-11}.
Therefore the proof of Lemma \ref{lem4.1} is completed. $\hfill\Box$

\smallskip

Motivated by \cite{Guo-Iyer-CMP}, we homogenize the boundary conditions in \eqref{4.15E}. We take  $\xi(y) \in C^\infty(\R_+)$ with compact support such that
\begin{align}
	\xi(0)=1,\quad \int_0^\infty \xi(y) dy =0,\quad \xi^{(i)}(0)=0,\,\, \mbox{for}\,\, i\geq 1.
\end{align}
Let $\dis I_{\xi}(y):=-\int_y^{+\infty} \xi(s) ds$. We define 
\begin{align}\label{4.20E}
\begin{cases}
\dis \bfu:=\bar{\bfu} - (\bar{\rho}^0_p \tilde{u}_b)(x,0)\, \xi(y),\\
\dis\bfv:= \bar{\bfv}   - (\bar{\rho}^0_p \tilde{u}_b)_x(x,0)\, I_{\xi}(y),\\
\dis \mathbf{\Theta}:= T_p - \tilde{T}_B,
\end{cases}
\end{align}
where
\begin{align}\label{4.21E}
\tilde{T}_B:=
\begin{cases}
y\chi(y)\tilde{T}_b(x) \quad \mbox{for NBC},\\[1.5mm]
\chi(y) \tilde{T}_b(x)\,\,\, \quad \mbox{for DBC}.
\end{cases}
\end{align}

Then \eqref{3.73} is rewritten as
\begin{align}\label{4.22E} 
\begin{cases}
 \bfu_x + \bfv_y=0,\\[1mm]
\bar{u}^0_{p} \pa_x  \bfu   + [\pa_y \bar{u}^{0}_p -   \bar{u}^0_p \f{\pa_y\bar{T}^0_p}{\bar{T}^0_p}]\,  \bfv   - \mu \bar{T}_p^0 \pa_{yy}\bfu = Q (\bfu, \Theta, G) ,\\[1mm]
2\bar{\rho}^0_p \bar{u}^0_p \pa_x\mathbf{\Theta} + 2\bar{\rho}^0_p  \bar{v}^0_p \pa_y\mathbf{\Theta} + 2\pa_x\bar{T}^0_p\, \bfu + 2\pa_y\bar{T}^0_p\, \bfv - \kappa \pa_{yy}\mathbf{\Theta}
=  H(\bfu, \Theta, G),
	\end{cases}
\end{align}
where 
\begin{align}\label{4.23E}
Q(\bfu, \mathbf{\Theta}, G):&=-\bar{u}^0_{p}  \pa_x(\bar{\rho}^0_p \tilde{u}_b)(x,0)\, \xi(y) - [\pa_y \bar{u}^{0}_p - \bar{\rho}^0_p  \bar{u}^0_p \pa_y\bar{T}^0_p]\, \pa_x(\bar{\rho}^0_p \tilde{u}_b)(x,0) \, I_{\xi}(y) \nonumber\\
&\quad+  \mu \bar{T}_p^0 \pa_x(\bar{\rho}^0_p \tilde{u}_b)(x,0)\, \xi''(y)  + Q_1\big(\bfu + (\bar{\rho}^0_p \tilde{u}_b)(x,0)\, \xi(y), \mathbf{\Theta} + \tilde{T}_B, G\big),\\
H(\bfu, \mathbf{\Theta}, G):&= -2\bar{\rho}^0_p \bar{u}^0_p \pa_x\tilde{T}_B - 2\bar{\rho}^0_p  \bar{v}^0_p \pa_y\tilde{T}_B + \kappa \pa_{yy}\tilde{T}_B - 2\pa_x\bar{T}^0_p\, (\bar{\rho}^0_p \tilde{u}_b)(x,0)\, \xi(y)
 \nonumber\\
&\quad  - 2\pa_y\bar{T}^0_p\, (\bar{\rho}^0_p \tilde{u}_b)_x(x,0)\, I_{\xi}(y) +H_1\big(\bfu +  (\bar{\rho}^0_p \tilde{u}_b)(x,0)\, \xi(y), \mathbf{\Theta} + \tilde{T}_B, G\big).\label{4.24E}
\end{align}
The initial boundary conditions  of \eqref{4.22E} become
\begin{align}\label{4.25E}
\begin{cases}
\bfu |_{x=0}= \big(\bar{\rho}^0_p \tilde{u}_0 + \bar{u}^0_p [\bar{\rho}^0_p P_1 - (\bar{\rho}^0_p)^2 \tilde{T}_0 ]- (\bar{\rho}^0_p \tilde{u}_b)(x,0)\, \xi(y)\big)(0,y),\\[1mm]
\mathbf{\Theta}|_{x=0}=\tilde{T}_0-\tilde{T}_B,\quad 
\dis (\bfu, \bfv)|_{y=0}=(0,0), \quad  \lim_{y\to \infty} (\bfu, \mathbf{\Theta})=(0,0).\\[1mm]
\dis \mathbf{\Theta}_y|_{y=0}= 0\quad \mbox{for NBC}, 
\\[1mm]
\dis \mathbf{\Theta}|_{y=0}=0\quad\,\,  \mbox{for DBC}. 
\\[1mm]
\end{cases}
\end{align}

\subsection{Compatibility conditions}\label{sec3.2}
We shall describe the compatibility conditions of initial boundary data for   $(\bfu, \bfv, \mathbf{\Theta})$. 

We have from $\eqref{4.22E}_2$ that
\begin{align*}
\Big(\f{\bar{T}^0_p\bfv}{\bar{u}^0_p}\Big)_y=-\Big(\f{\bar{T}^0_p}{\bar{u}^0_p}\Big)^2 \Big\{\mu  \pa_{yy}\bfu +  \f{1}{\bar{T}_p^0}Q (\bfu, \mathbf{\Theta}, G)\Big\},
\end{align*}
which implies that 
\begin{align}
\bfv =-\f{\bar{u}^0_p}{\bar{T}^0_p} \int_0^y \Big(\f{\bar{T}^0_p}{\bar{u}^0_p}\Big)^2 \Big\{\mu  \pa_{yy}\bfu +  \f{1}{\bar{T}_p^0}Q (\bfu, \mathbf{\Theta}, G)\Big\} dz.
\end{align}
Then we can define $\bfv|_{x=0}$, and
\begin{align}\label{3.24}
\begin{split}
\pa_x\bfu|_{x=0}&=\bfv_y|_{x=0},\\
\pa_x\mathbf{\Theta} |_{x=0}
&= \f{\bar{T}^0_p}{2\bar{u}^0_p}\big\{\kappa \pa_{yy}\mathbf{\Theta} - 2\bar{\rho}^0_p  \bar{v}^0_p \pa_y\mathbf{\Theta} - 2\pa_x\bar{T}^0_p\, \bfu - 2\pa_y\bar{T}^0_p\, \bfv  +  H(\bfu, \mathbf{\Theta}, G) \big\}.
\end{split}
\end{align}

To make \eqref{3.24}  meaningful, the following compatibility conditions must be satisfied 
\begin{align}\label{3.25}
	\begin{split}
	\Big\{\mu  \pa_{yy}\bfu +  \f{1}{\bar{T}_p^0}Q (\bfu, \mathbf{\Theta}, G)\Big\} \Big|_{x=0} &=O(1) y^2,\quad y\in [0,1], \\
	\big\{\kappa \pa_{yy}\mathbf{\Theta} - 2\bar{\rho}^0_p  \bar{v}^0_p \pa_y\mathbf{\Theta} - 2\pa_x\bar{T}^0_p\, \bfu - 2\pa_y\bar{T}^0_p\, \bfv  +  H(\bfu, \mathbf{\Theta}, G) \big\}\big|_{x=0} &=O(1)y,\quad y\in [0,1].
	\end{split}
\end{align}
Inductively, assuming high order compatibility conditions, we can also define 
\begin{align}
	\pa_x^k\Big(\f{\bar{T}^0_p\bfv}{\bar{u}^0_p}\Big)_y\Big|_{x=0},\quad  
	\pa_x^k\bfu\big|_{x=0}, \quad    \pa_x^k\bfv\big|_{x=0} \quad \mbox{and}\quad \pa_x^k\mathbf{\Theta}\big|_{x=0} \quad \mbox{for}\,\,\, k\geq 0.
\end{align}

\subsection{Existence and uniform estimate for solutions of \eqref{4.22E}}
In this subsection, we aim to establish the existence and uniform estimates for solutions of \eqref{4.22E} and \eqref{4.25E}.

\subsubsection{The $\theta$-approximate problem} To solve $(\bfu, \bfv, T_p)$, we consider the following $\theta$-approximate problem (non-degenerate)
\begin{align}\label{5.114}
	\begin{cases}
		\bfu_x + \bfv_y=0,\\[1mm]
		\bar{u}^0_{p,\theta} \pa_x \bfu  +  [\pa_y \bar{u}^{0}_p - \bar{u}^0_{p,\theta}  \f{\pa_y\bar{T}^0_p}{\bar{T}^0_p}]\,  \bfv  - \mu \bar{T}_p^0 \pa_{yy}\bfu = Q(\bfu, \Theta, G),\\[1mm]
		2\bar{\rho}^0_p \bar{u}^0_{p,\theta} \pa_x\mathbf{\Theta} + 2\bar{\rho}^0_p  \bar{v}^0_p \pa_y\mathbf{\Theta} + 2\pa_x\bar{T}^0_p\, \bfu + 2\pa_y\bar{T}^0_p\, \bfv - \kappa \pa_{yy}\mathbf{\Theta} = H(\bfu, \Theta, G),
	\end{cases}
\end{align}
where $\bar{u}^0_{p,\theta}:=\theta+\bar{u}^0_{p}$ with $\theta\in (0,1]$, the initial  boundary condition is \eqref{4.25E}. It is better to denote the solutions of \eqref{5.114} as $(\bfu^{\t}, \bfv^{\t}, \mathbf{\Theta}^{\t})$, but we shall drop the superscript $\theta$ from $(\bfu^{\t}, \bfv^{\t}, \mathbf{\Theta}^{\t})$ in the remainder of this subsection when no confusion arises.

Noting $\eqref{5.114}_1$, we introduce the stream function 
\begin{align}\label{5.72}
	\psi(x,y):=\int_0^y \bfu(x,z) dz\quad \Longrightarrow\quad \bfu=\psi_y,\quad \bfv=-\psi_x.
\end{align}
We define 
\begin{align}\label{3.29}
\hat{u}^0_{p,\theta}:=\f{\bar{u}^0_{p,\theta}}{\bar{T}^0_p} \cong \bar{u}^0_{p,\theta},
\end{align}
and
\begin{align}\label{5.120}
	\Psi:=\psi_{y} - \f{\pa_y \hat{u}^{0}_{p,\t}}{\hat{u}^{0}_{p,\t}}\, \psi\quad \Longleftrightarrow\quad \psi(x,y)=\hat{u}^0_{p,\t}(x,y) \int_0^y \Big(\f{\Psi}{\hat{u}^0_{p,\t}}\Big)(x,z) dz.
\end{align}
We also assume the high order compatibility conditions on the initial boundary data hold.

\begin{lemma}\label{lem4.2}
The system \eqref{5.114} is equivalent to the following parabolic equations
\begin{align}\label{5.115}
\begin{cases}
\dis \hat{u}^0_{p,\theta} \Psi_x  - \mu \Psi_{yy}
= \f{1}{\bar{T}_p^0 } Q(\psi_y, \mathbf{\Theta}, G) + \mu \Big\{\f{\pa_y \hat{u}^{0}_{p,\t}}{\hat{u}^{0}_{p,\t}}\, \psi\Big\}_{yy} -  \hat{u}^0_{p,\theta}  \Big(\f{\pa_y \hat{u}^{0}_{p,\t}}{\hat{u}^{0}_{p,\t}}\Big)_x \psi ,\\[1mm]
\dis 2\bar{\rho}^0_p \bar{u}^0_{p,\theta} \pa_x \mathbf{\Theta} + 2\bar{\rho}^0_p  \bar{v}^0_p \pa_y\mathbf{\Theta} + 2\pa_x\bar{T}^0_p\, \psi_y - 2\pa_y\bar{T}^0_p\, \psi_x  - \kappa \pa_{yy}\mathbf{\Theta}
=  H(\psi_y, \mathbf{\Theta}, G).
\end{cases}
\end{align}
There exists a small $x_{\t}>0$ such that there exists a smooth solution $(\Psi, \mathbf{\Theta})$ of \eqref{5.115} $($which is equivalent to \eqref{5.114}$)$ in $x\in [0,x_{\t}]$ with sufficient decay in $y\in \mathbb{R}_+$.
\end{lemma}

\noindent{\bf Proof.} 1. 
It is clear that  $\eqref{5.114}$ is equivalent to the following system 
\begin{align}\label{5.116}
\begin{cases}
\dis \bar{u}^0_{p,\theta} \psi_{xy} - [\pa_y \bar{u}^{0}_p -  \bar{u}^0_{p,\t} \f{\pa_y\bar{T}^0_p}{\bar{T}^0_p}]\, \psi_x - \mu \bar{T}_p^0 \psi_{yyy} = Q(\psi_y, \mathbf{\Theta}, G),\\[1mm]
\dis 2\bar{\rho}^0_p \bar{u}^0_{p,\theta} \pa_x\mathbf{\Theta} + 2\bar{\rho}^0_p  \bar{v}^0_p \pa_y\mathbf{\Theta} + 2\pa_x\bar{T}^0_p\, \psi_y - 2\pa_y\bar{T}^0_p\, \psi_x - \kappa \pa_{yy}\mathbf{\Theta}
=  H(\psi_y, \mathbf{\Theta}, G).
\end{cases}
\end{align}

For $\eqref{5.116}_1$, it is not a classical equation. We need to make further formulation. 
It is noted that 
\begin{align*}
\bar{u}^0_{p,\theta} \psi_{xy} - [\pa_y \bar{u}^{0}_p -  \bar{u}^0_{p,\t} \f{\pa_y\bar{T}^0_p}{\bar{T}^0_p}]\, \psi_x
&=\bar{T}^0_p \hat{u}^0_{p,\theta} \psi_{xy} - \bar{T}^0_p \pa_y \hat{u}^0_{p,\theta} \psi_x = \bar{T}^0_p \hat{u}^0_{p,\theta} \Big\{\psi_{xy} - \f{\pa_y \hat{u}^{0}_{p,\t}}{\hat{u}^{0}_{p,\t}}\, \psi_x\Big\}\nonumber\\
&=\bar{T}^0_p \hat{u}^0_{p,\theta} \Psi_x  
+ \bar{T}^0_p \hat{u}^0_{p,\theta}  \Big(\f{\pa_y \hat{u}^{0}_{p,\t}}{\hat{u}^{0}_{p,\t}}\Big)_x \psi, 
\end{align*}
and
\begin{align*}
\begin{split}
\psi_{yyy} &= \Big\{\psi_{y} - \f{\pa_y \hat{u}^{0}_{p,\t}}{\hat{u}^{0}_{p,\t}}\, \psi\Big\}_{yy} + \Big\{\f{\pa_y \hat{u}^{0}_{p,\t}}{\hat{u}^{0}_{p,\t}}\, \psi\Big\}_{yy} = \Psi_{yy} + \Big\{\f{\pa_y \hat{u}^{0}_{p,\t}}{\hat{u}^{0}_{p,\t}}\, \psi\Big\}_{yy} ,
\end{split}
\end{align*}
then  $\eqref{5.116}_1$ becomes
\begin{align}\label{5.119}
\bar{T}^0_p \hat{u}^0_{p,\theta} \Psi_x  - \mu \bar{T}_p^0 \Psi_{yy}
 = Q(\psi_y, \mathbf{\Theta}, G) + \mu \bar{T}_p^0 \Big\{\f{\pa_y \hat{u}^{0}_{p,\t}}{\hat{u}^{0}_{p,\t}}\, \psi\Big\}_{yy} -   \bar{T}^0_p \hat{u}^0_{p,\theta}  \Big(\f{\pa_y \hat{u}^{0}_{p,\t}}{\hat{u}^{0}_{p,\t}}\Big)_x \psi.
\end{align}
Thus we have proved \eqref{5.115}.

2. For the linear uniform parabolic equations  \eqref{5.115}, by applying the standard parabolic theory,  we can obtain the local  solution  $(\Psi, \mathbf{\Theta})$  in $x\in [0,x_{\theta}]$ with enough decay rate and  the lifespan $x_{\t}>0$ depending on $\theta\in (0,1]$.  Therefore the proof of Lemma \ref{lem4.2} is completed. $\hfill\Box$

\smallskip

\begin{remark}
Thus, to extend the local solution in Lemma \ref{4.20E} to a time interval independent of $\theta$, it is essential to establish some estimates that are uniform in $\theta$. We also note that the solution $(\Psi, \mathbf{\Theta})$ of \eqref{5.15} (or equivalently, the solution $(\bfu, \bfv)$ of \eqref{5.114}) necessarily depends on $\theta$, which we will omit in the sequel when no confusion arises.
\end{remark}

\subsubsection{Estimates on quotient}
Recall \eqref{3.29}, we define the corresponding  quotient
\begin{align} 
\bfq=\f{\psi}{\hat{u}^0_{p,\theta}}.
\end{align}
Then the momentum equation $\eqref{5.114}_2$ becomes
\begin{align}\label{5.75}
(\hat{u}^0_{p,\theta})^2 \bfq_{xy} - \mu \psi_{yyy} \equiv \bar{\rho}_p^0 (\bar{u}^0_{p,\theta})^2 \mathbf{q}_{xy} - \mu \pa_{yyy}\big[\bar{u}^0_{p,\theta} \bfq\big] = R,
\end{align}
where
\begin{align}\label{5.76}
R:&= - \hat{u}^0_{p,\theta} \pa_x \hat{u}^0_{p,\t} \bfq_y - \big[\hat{u}^0_{p,\theta} \pa_{xy}\hat{u}^0_{p,\t}  -  \pa_y\hat{u}^0_{p,\t}\,  \pa_x\hat{u}^0_{p,\t} \big] \bfq + \f{1}{\bar{T}_p^0 } Q(\psi_y, \mathbf{\Theta}, G).
\end{align}
Noting $(\bfu, \bfv)|_{y=0}$, we have 
\begin{align}\label{4.38E}
\bfq|_{y=0}=\bfq_y|_{y=0}=0,\quad \psi|_{y=0}=\psi_{y}|_{y=0}=0.
\end{align}

For quotient, we define the norms
\begin{align}
	\begin{split}
		\|\bfq_y\|^2_{\MBX_{\t, w}}:&=  \|\sqrt{\bar{u}^0_{p,\t}} \bfq_{yy}w\|^2 + \sup_{x_0\in[0,L]}\|\bar{u}^0_{p,\theta} \bfq_yw\|^2_{x=x_0},\\
		\|\bfq_y\|^2_{\MBX^k_{\t, w}}:&= \sum_{i=0}^k \|\pa_x^iq_y\|^2_{\MBX_{\t, w}}
		\equiv \sum_{i=0}^k \Big\{\|\sqrt{\bar{u}^0_{p,\theta}} \pa_{x}^i \bfq_{yy}w\|^2 + \sup_{x_0\in[0,L]}\|\bar{u}^0_{p,\theta} \pa_{x}^i \bfq_yw\|^2_{x=x_0}\Big\}.
	\end{split}
\end{align}
For temperature, we define norms
\begin{align}
	\begin{split}
		\|h\|^2_{\MBY_{\t, w}} &:= \kappa \|h_{y}w\|^2 + \sup_{x_0\in[0,L]}\|\sqrt{\bar{u}^0_{p,\t}}h w\|^2_{x=x_0},\\
		\|h\|^2_{\MBY^k_{\t, w}}&:=\sum_{i=0}^k \|\pa_x^ih\|^2_{\MBY_{\t, w}}
		=\sum_{i=0}^k \Big\{\kappa \|\pa_x^ih_{y}w\|^2 + \sup_{x_0\in[0,L]}\|\sqrt{\bar{u}^0_{p,\t}}\pa_x^ih w\|^2_{x=x_0}\Big\}.
	\end{split}
\end{align}
For $\theta=0$, we denote $\mathbb{X}_{w}=\mathbb{X}_{0,w}$,  $\mathbb{X}^k_w=\mathbb{X}^k_{0, w}$ and  $\mathbb{Y}_w=\mathbb{Y}_{0, w}$,  $\mathbb{Y}^k_w=\mathbb{Y}^k_{0, w}$.  We also note that $\MBX_{\t, w}\cong \MX_{\t, w}$ and $\MBY_{\t, w}\cong \MY_{\theta, w}$.

\begin{lemma}\label{lem4.3}
For $0\leq k\leq n$, it holds that 
\begin{align}\label{5.78}
\|\bfq_y\|^2_{\MBX^n_{\t, w}}
&\lesssim \sum_{k=0}^{n} \|\bar{u}^0_{p,\t}\pa_x^k \bfq_y w\|^2_{x=0} +  L^{\f12}\sum_{i=0}^{n-1}\|\pa_x^{i} \bfq_{yy}w\bar{\chi}\|^2_{x=0}   + L^{\f16-} \big\{\|\bfq_{y}\|^2_{\MBX^n_{\t, w}}  + \sum_{i=0}^{n} \|\pa_x^i R w\|^2\big\}.
\end{align}
\end{lemma}

\noindent{\bf Proof.} The proof is similar to Lemma \ref{lem6.8}, but it is much easier now since we are dealing with a linear problem. 

1. Applying $\pa^k_x$ to \eqref{5.75}, one gets
\begin{align}\label{5.79}
	\pa^k_x\big((\hat{u}^0_{p,\theta})^2 \bfq_{xy}\big)  - \mu \,\pa_x^k\pa_{yyy} [\hat{u}^0_{p,\t}\bfq] = \pa_x^kR,\,\,\,\mbox{for}\,\, 0\leq k\leq n.
\end{align}
Multiplying \eqref{5.79} by $\pa_x^k \bfq_{y}w^2$, one has that 
\begin{align}\label{5.80}
&\big\la \pa^k_x\big( (\hat{u}^0_{p,\t})^2  \bfq_{xy}\big),\, \pa_x^k\bfq_{y}w^2\big\ra -  \big\la \pa_x^k [\hat{u}^0_{p,\t} \bfq]_{yyy},\,  \pa^k_x \bfq_y w^2\big\ra \nonumber\\
&\leq L^{\f16-} \|\pa_x^k R w\|^2 + L^{\f16-}\|\pa_x^k \bfq_{y}\|^2_{\MBX_w},
\end{align}
where we have used 
\begin{align*}
	\big\la \pa_x^k R,\, \pa_x^k\bfq_{y}w^2\big\ra \lesssim \|\pa_x^k R w\|\cdot \|\pa_x^k\bfq_{y}w\| \lesssim L^{\f16-} \|\pa_x^k R w\|^2 + L^{\f16-}\|\pa_x^k \bfq_{y}\|^2_{\MBX_w}.
\end{align*}

2. Integrating by parts in $x$  and using \eqref{6.48}, one obtains
\begin{align}\label{5.81}
	&\big\la \pa^k_x\big( (\hat{u}^0_{p,\t})^2  \bfq_{xy}\big),\, \pa_x^k \bfq_{y}w^2\big\ra \nonumber\\
	&=\big\la (\hat{u}^0_{p,\t})^2  \pa^k_x\bfq_{xy},\, \pa_x^k\bfq_{y}w^2\big\ra + \sum_{i=0}^{k-1} C_{k}^i \big\la  \pa_x^{k-i}( (\hat{u}^0_{p,\t})^2) \, \pa_x^i \bfq_{xy} ,\, \pa_x^k\bfq_{y}w^2\big\ra\nonumber\\
	&\geq \f12 \|\hat{u}^0_{p,\t}\pa_x^k \bfq_y w\|^2_{x=x_0} - \f12 \|\hat{u}^0_{p,\t} \pa_x^k \bfq_y w\|^2_{x=0}  - L^{\f13-} \|\bfq_{y}\|^2_{\MBX^k_{\t, w}}.
\end{align}

3.   It is clear that 
\begin{align}\label{5.82}
-  \big\la \pa^k_x [\hat{u}^0_{p,\t} \bfq]_{yyy},\, \pa_x^k\bfq_{y} w^2\big\ra  
&= -  \big\la [\hat{u}^0_{p,\t} \pa^k_x\bfq]_{yyy},\, \pa_x^k \bfq_{y}w^2\big\ra -  \sum_{i=1}^{k} C_k^i  \big\la [\pa^i_x\hat{u}^0_{p,\t} \cdot \pa^{k-i}_x\bfq]_{yyy},\,  \pa_x^k \bfq_{y}w^2\big\ra .
\end{align}
We point out that the second term vanishes  for $k=0$.

\smallskip

3.1. For the first term on RHS of \eqref{5.82}, noting $\pa^k_x\bfq_{y}|_{y=0}=0$, we have 
\begin{align}\label{5.83}
	-\big\la  [\hat{u}^0_{p,\t}  \pa^k_x\bfq]_{yyy},\,   \pa^k_x\bfq_{y}w^2\big\ra 
	&=\big\la [\hat{u}^0_{p,\t} \pa^k_x\bfq]_{yy},\,   \pa^k_x\bfq_{yy}w^2\big\ra + \big\la   [\hat{u}^0_{p,\t}  \pa^k_x\bfq]_{yy},\,   2\pa^k_x\bfq_{y}ww_y\big\ra \nonumber\\
	&=\|\sqrt{\hat{u}^0_{p,\t}}\pa^k_x \bfq_{yy}w\|^2 + \big\la 2 \pa_y\hat{u}^0_{p,\t} \, \pa^k_x\bfq_{y} + \pa_{yy}\hat{u}^0_{p,\t}\, \pa^k_x\bfq,\,   \pa^k_x\bfq_{yy} w^2\big\ra \nonumber\\
	&\quad  + \big\la \hat{u}^0_{p,\t}\pa^k_x\bfq_{yy} + 2 \pa_y\hat{u}^0_{p,\t} \, \pa^k_x\bfq_{y} + \pa_{yy}\hat{u}^0_{p,\t}\, \pa^k_x\bfq,\, 2\pa^k_x\bfq_{y} ww_y\big\ra.
\end{align}

Integrating by parts in $y$, one obtains
\begin{align}\label{5.84}
\big\la 2 \pa_y\hat{u}^0_{p,\t} \pa^k_x\bfq_{y} ,\,   \pa^k_x\bfq_{yy} w^2\big\ra
= - \big\la  (\pa_y\hat{u}^0_{p,\t}w^2)_y,\,   |\pa^k_x\bfq_{y}|^2 \big\ra 
&\lesssim L^{\f13-}  \|\bfq_{y}\|^2_{\MBX^k_w},\\
\big\la \hat{u}^0_{p,\t}\pa^k_x\bfq_{yy} ,\, 2\pa^k_x\bfq_{y} ww_y\big\ra + \big\la 2 \pa_y\hat{u}^0_{p,\t} \pa^k_x\bfq_{y},\, 2\pa^k_x\bfq_{y} ww_y \big\ra 
&\lesssim  L^{\f13-} \|\bfq_{y}\|^2_{\MBX^k_{\t, w}}.\label{5.85}
\end{align}

Noting the decay property of $\pa_y^l\bar{u}^0_{p}$ ($l\geq 1$), we have 
\begin{align}\label{5.86}
\big\la \pa_{yy}\hat{u}^0_{p,\t}\pa^k_x\bfq,\,   \pa^k_x\bfq_{yy} w^2\big\ra 
&= - \big\la \pa_{yy}\hat{u}^0_{p,\t}\pa^k_x\bfq_y,\,   \pa^k_x\bfq_{y} w^2\big\ra - \big\la (\pa_{yy}\hat{u}^0_{p,\t}w^2)_y \pa^k_x\bfq,\,   \pa^k_x\bfq_{y} \big\ra \nonumber\\
& \lesssim L^{\f13-} \|\bfq_{y}\|^2_{\MBX^k_{\t, w}},\\
\big|\big\la \pa_{yy}\hat{u}^0_{p,\t}\pa^k_x\bfq,\, 2\pa^k_x\bfq_{y} ww_y \big\ra\big| 
&\lesssim L^{\f13-} \|\bfq_{y}\|^2_{\MBX^k_{\t, w}}.\label{5.87}
\end{align}
Thus, substituting \eqref{5.84}-\eqref{5.87} into \eqref{5.83}, one obtains that
\begin{align}\label{5.88}
	-\big\la [\hat{u}^0_{p,\t} \pa^k_x\bfq]_{yyy} ,\,   \pa^k_x \bfq_{y}w^2\big\ra   
	&\geq \f78 \|\sqrt{\hat{u}^0_{p,\t}}\pa^k_x \bfq_{yy}w\|^2 - L^{\f13-}\|\bfq_y \|^2_{\MBX^k_w}.
\end{align}

3.2. For the second term (vanishing when $k=0$) on RHS of \eqref{5.82}, we note that ($1\leq i\leq k$)
\begin{align}\label{5.89}
	\big\la [\pa^i_x\hat{u}^0_{p,\t} \pa^{k-i}_x\bfq]_{yyy},\,  \pa_x^k\bfq_{y}w^2\big\ra 
	&=\big\la \pa^i_x \pa^3_{y}\hat{u}^0_{p,\t} \pa^{k-i}_x\bfq,\,  \pa_x^k\bfq_{y}w^2\big\ra  + \big\la 3\pa^i_x \pa^2_{y}\hat{u}^0_{p,\t}\pa^{k-i}_x\bfq_y,\,  \pa_x^k\bfq_{y}w^2\big\ra
	\nonumber\\
	&\,\,\, + \big\la 3\pa^i_x \pa_y\hat{u}^0_{p,\t}  \pa^{k-i}_x\bfq_{yy},\,  \pa_x^k\bfq_{y}w^2\big\ra  + \big\la \pa^i_x\hat{u}^0_{p,\t} \pa^{k-i}_x\bfq_{yyy},\,  \pa_x^k\bfq_{y}w^2\big\ra.
\end{align}
It is clear that 
\begin{align}
\sum_{i=1}^{k}\big\la \pa^i_x\pa^3_{y}\hat{u}^0_{p,\t} \pa^{k-i}_x\bfq,\,  \pa_x^k\bfq_{y}w^2\big\ra
&\lesssim \sum_{i=1}^{k}\|\pa_x^{k-i}\bfq\|_{L^2_xL^\infty_y} \|\pa^i_x\pa^3_{y}\hat{u}^0_{p,\t} w\|_{L^\infty_xL^2_y}\cdot \|\pa_x^k\bfq_{y}w\| \nonumber\\
&\lesssim L^{\f13-}  \|\pa_x^k\bfq_{y}\|^2_{\MBX_w} , \\
\sum_{i=1}^{k}\big\la 3\pa^i_x\pa^2_{y}\hat{u}^0_{p,\t} \pa^{k-i}_x\bfq_y,\,  \pa_x^k\bfq_{y}w^2\big\ra
&\lesssim \sum_{i=1}^{k}\|\pa^i_x\pa_y^2\hat{u}^0_{p,\t}\|_{L^\infty} \|\pa^{k-i}_x\bfq_yw\|\cdot  \|\pa_x^k\bfq_{y}w\| \nonumber\\ 
&\lesssim L^{\f13-}  \|\pa_x^k\bfq_{y}\|^2_{\MBX_w},
\end{align}
and
\begin{align} \label{5.92}
&\sum_{i=1}^{k}\big\la 3\pa^i_x\pa_y\hat{u}^0_{p,\t} \pa^{k-i}_x\bfq_{yy} + \pa^i_x\hat{u}^0_{p,\t} \pa^{k-i}_x\bfq_{yyy},\,  \pa_x^k\bfq_{y}w^2\big\ra  \nonumber\\
&\lesssim L^{\f16-}  \|\bfq_y\|^2_{\MBX^k_w} + L^{\f16-}  \sum_{i=0}^{k-1}  \|(\hat{u}^0_{p,\t}\pa_x^{i}\bfq_{yyy}, \pa_x^{i}\bfq_{yy})w\|^2.
\end{align}
Thus we have from \eqref{5.89}-\eqref{5.92} that
\begin{align}\label{5.93}
-  \sum_{i=1}^{k} C_k^i  	\big\la [\pa^i_x\hat{u}^0_{p,\t} \pa^{k-i}_x\bfq]_{yyy},\,  \pa_x^k\bfq_{y}w^2\big\ra 
	& \lesssim  L^{\f16-}  \|\bfq_y\|^2_{\MBX^k_w} + L^{\f16-}  \sum_{i=0}^{k-1}  \|(\hat{u}^0_{p,\t} \pa_x^{i}\bfq_{yyy}, \pa_x^{i}\bfq_{yy})w\|^2.
\end{align}

3.3. Substituting \eqref{5.93} and \eqref{5.88} into \eqref{5.82}, one obtains that 
\begin{align}\label{5.94}
  -  \big\la \pa^k_x [\hat{u}^0_{p,\t} \bfq]_{yyy},\, \pa_x^k \bfq_{y} w^2\big\ra 
&\geq \f78 \|\sqrt{\hat{u}^0_{p,\t}}\pa^k_x \bfq_{yy}w\|^2 - L^{\f16-}  \|\bfq_y\|^2_{\MBX^k_w} \nonumber\\
&\quad  - L^{\f16-}  \sum_{i=0}^{k-1}  \|(\hat{u}^0_{p,\t} \pa_x^{i}\bfq_{yyy}, \pa_x^{i}\bfq_{yy})w\|^2.
\end{align}

4. By similar arguments as the proof of \eqref{6.38-0}, \eqref{6.130-0} and \eqref{6.126-0}, one gets that 
\begin{align}
\sum_{i=0}^{n-1}	\|\pa_x^{i} \bfq_{yy}w\|
&\lesssim L^{\f12}\sum_{i=0}^{n-1}\|\pa_x^{i} \bfq_{yy}w\bar{\chi}\|_{x=0}   + L^{\f12-} \|\bfq_y\|_{\MBX^{n}_w} +  \sum_{i=1}^{n-1}  \|\pa_x^i\psi_{yyy} \la y\ra^2\| ,\label{5.95-0} \\
\sum_{i=0}^{n-1}	\|\hat{u}^0_{p,\t} \pa_x^{i} \bfq_{yyy}w\|
&\lesssim  \|\bfq_y\|_{\MBX^n_{\t, w}} + \sum_{i=0}^{n-1}\|\pa_x^{i}\psi_{yyy}w\| + \sum_{j=0}^{n-1} \|\pa_x^{j} \bfq_{yy}w\|,\label{5.95-01}
\end{align}
and 
\begin{align}\label{5.95-00}
\sum_{i=0}^{n-1} \|\pa_x^i\psi_{yyy}w\| & \lesssim L^{\f16-}  \|\bfq_y\|_{\MBX^n_{\t, w}} + \sum_{i=0}^{n-1} \|\pa_x^i R w\|.
\end{align}

5. Substituting \eqref{5.94}  and \eqref{5.81} into \eqref{5.80}, then using \eqref{5.95-0}-\eqref{5.95-00}, one 
concludes \eqref{5.78}.
Therefore the proof of Lemma \ref{lem4.3} is completed. $\hfill\Box$

\medskip

\begin{lemma}\label{lem5.9}
It holds that 
\begin{align}
\sum_{k=0}^n\|\pa_x^k\psi_yw\| &\lesssim L^{\f16-} \|\bfq_y\|_{\MBX^n_{\t, w}},\label{5.95}\\
\sum_{k=0}^n\|\pa_x^k\psi_{yy}w\|&\lesssim  \|\bfq_y\|_{\MBX_{\theta, w}^{n}},\label{5.95-1}\\
\sum_{i=0}^{n-1}\|\pa_x^{i}\psi_{yy}w\|^2_{x=x_0} 
&\lesssim  \sum_{i=0}^{n-1}\|\pa_x^{i}\psi_{yy}w\|^2_{x=0} + L^{\f16-} \sum_{i=0}^{n-1} \|\pa_x^{i}R w\|^2 + L^{\f16-}  \|\bfq_y\|^2_{\MBX_{\theta, w}^{n}},\label{5.95-2}
\end{align}
and
\begin{align}\label{5.96}
	\begin{split}
	\sum_{i=0}^{n-1}\|\pa_x^{i}\psi_{yyyy}w\| &\lesssim  \|\bfq_{y}\|_{\MBX^n_{\theta, w}}  + \sum_{i=0}^{n-1}\|\pa_x^{i}R_y w\|.
\end{split}
\end{align}
\end{lemma}

\noindent{\bf Proof.}   For $0\leq k\leq n$, we have 
\begin{align}
\begin{split}
\|\pa_x^k\psi_yw\| 
&\lesssim \sum_{i=0}^k \big\{\|\pa_x^i\hat{u}^0_{p,\t}\pa_x^{k-i}\bfq_yw\|+  \|\pa_x^i\pa_y\hat{u}^0_{p,\t}\pa_x^{k-i}\bfq w\|\big\}  
\lesssim L^{\f16-}  \|\bfq_y\|_{\MBX^n_{\theta, w}},\\
\|\pa_x^k\psi_{yy}w\|
&\lesssim   \sum_{i=0}^k \|(\pa_x^i \hat{u}^0_{p,\t}\pa_x^{k-i}\bfq_{yy}, \pa_x^i\pa_y\hat{u}^0_{p,\t}\pa_x^{k-i} \bfq_{y}, \pa_x^i\pa^2_y\hat{u}^0_{p,\t}\pa_x^{k-i} \bfq) w\|
\lesssim  \|\bfq_y\|_{\MBX_{\theta, w}^{k}}.
\end{split}
\end{align}
 
\smallskip

 For $0\leq k\leq n$, we have from \eqref{6.17} that 
\begin{align}\label{5.99} 
|\hat{u}^0_{p,\t}|^2 \pa_x^{k}\bfq_y - \pa_x^{k-1}\psi_{yyy}  + \sum_{i=1}^{k-1} C_{k-1}^i \pa^i_x\big(|\hat{u}^0_{p,\t}|^2\big)\pa_x^{k-i} \bfq_{y} = \pa_x^{k-1}R.
\end{align}
Multiplying \eqref{5.99} by $\pa_x^{k}\psi_yw^2$, we obtain
\begin{align}\label{5.101}
\f12 \|\pa_x^{k-1}\psi_{yy}w\|^2_{x=x_0} 
&\leq \f12 \|\pa_x^{k-1}\psi_{yy}w\|^2_{x=0} + \big\{\|\pa_x^{k-1}R w\| + \|\bfq_y\|_{\MBX_{\theta, w}^{k}}\big\}\|\pa_x^k\psi_yw\| \nonumber\\
&\leq\f12 \|\pa_x^{k-1}\psi_{yy}w\|^2_{x=0} + L^{\f16-}  \big\{\|\pa_x^{k-1}R w\|^2 + \|\bfq_y\|^2_{\MBX_{\theta, w}^{k}}\big\}.
\end{align}

\smallskip
 
Also, it follows from \eqref{5.99} that 
\begin{align}\label{5.100}
	\begin{split}
		\|\pa_x^{k-1}\psi_{yyyy}w\| &\lesssim  \|\bfq_{y}\|_{\MBX^k_{\theta, w}}  + \|\pa_x^{k-1}R_y w\|.
	\end{split}
\end{align}
Therefore the proof Lemma \ref{lem5.9}  is completed. $\hfill\Box$

\medskip
\begin{lemma}[Elliptic estimate]\label{lem4.5}
It holds that 
\begin{align}\label{5.53-0E}
& \sum_{k=0}^{n} \Big\{\|(\bar{u}^0_{p,\theta})^{\f52} \pa_x^k \bfq_{yy}w\|^2_{x=x_0} +  \|(\bar{u}^0_{p,\theta})^{3}\pa_x^{k+1} \bfq_yw\|^2 + \|\bar{u}^0_{p,\theta} \pa_x^k\psi_{yyy}w\|^2 \Big\}  \nonumber\\
&\lesssim \sum_{k=0}^{n}  \|(\bar{u}^0_{p,\theta})^{\f52}\pa_x^k \bfq_{yy}w\|^2_{x=0} +  \|\bfq_{y}\|^2_{\MBX^n_{\theta, w}} + \sum_{i=0}^{n} \|\pa^i_x R w\|^2.
\end{align}
The elliptic estimation \eqref{5.53E} play no role for us to close the uniform-in-$\theta$ estimate, it only help us to know that $\pa_x^{n+1}\bfq_y$ is well-defined in some sense.
\end{lemma}

\noindent{\bf Proof.} Let $0\leq k\leq n$. We have from \eqref{5.75} that 
\begin{align}\label{5.53E}
&\|(\hat{u}^0_{p,\t})^{3}\pa_x^{k+1} \bfq_yw\|^2 + \f34 \mu^2\|\hat{u}^0_{p,\t}\pa_x^k\psi_{yyy}w\|^2 - 2\mu \big\la (\hat{u}^0_{p,\t})^{4} \pa_x^{k+1} \bfq_y,\, \pa_x^k\psi_{yyy}w^2\big\ra \nonumber\\
	&\lesssim \|\pa_x R w\|^2 + L^{\f13-}\|\bfq_{y}\|^2_{\MBX^k_{\theta, w}}.
\end{align}

Integrating by parts, one obtains 
\begin{align}\label{6.269E}
	& - \big\la (\hat{u}^0_{p,\t})^{4} \pa_x^{k+1} \bfq_y,\, \pa_x^k\psi_{yyy}w^2\big\ra \nonumber\\
	&=\big\la (\hat{u}^0_{p,\t})^{4} \pa_x^{k+1} \bfq_{yy},\, \pa_x^k(\hat{u}^0_{p,\t} \bfq_{yy}) w^2\big\ra + \big\la  (\hat{u}^0_{p,\t})^{4} \pa_x^{k+1} \bfq_{yy},\, \pa_x^k(2\bar{u}^0_{py} \bfq_{y}) w^2\big\ra \nonumber\\
	&\quad  + \big\la (\hat{u}^0_{p,\t})^{4} \pa_x^{k+1} \bfq_{yy},\, \pa_x^k(\bar{u}^0_{pyy} \bfq)w^2\big\ra + \big\la \big\{ (\hat{u}^0_{p,\t})^{4}w^2\big\}_y \pa_x^{k+1} \bfq_y,\, \pa_x^k\psi_{yy} \big\ra.
\end{align}

For the first term on RHS of \eqref{6.269E}, 
a direct calculation shows that 
\begin{align}\label{6.277E}
&\big\la (\hat{u}^0_{p,\t})^{4} \pa_x^{k+1} \bfq_{yy},\, \pa_x^k(\hat{u}^0_{p,\t} \bfq_{yy}) w^2\big\ra \nonumber\\
&=\big\la (\hat{u}^0_{p,\t})^{4} \pa_x^{k+1} \bfq_{yy},\, \hat{u}^0_{p,\t} \pa_x^{k}\bfq_{yy} w^2\big\ra 
+ \sum_{i=1}^k C_k^i \big\la (\hat{u}^0_{p,\t})^{4} \pa_x^{k+1} \bfq_{yy} ,\, \pa_x^i\bar{u}^0_{p} \pa_x^{k-i} \bfq_{yy}w^2\big\ra\nonumber\\
&\geq \f12 \|(\hat{u}^0_{p,\t})^{\f52}\pa_x^k \bfq_{yy}w\|^2_{x=x_0} - \f12\|(\hat{u}^0_{p,\t})^{\f52}\pa_x^k \bfq_{yy}w\|^2_{x=0} - \f18 \|(\hat{u}^0_{p,\t})^{3}\pa_x^{k+1} \bfq_yw\|^2 \nonumber\\
&\quad - C \|\bfq_y\|^2_{\MBX^k_{\theta, w}}  - C \sum_{i=0}^{k-1} \|\pa_x^i R w\|^2,
\end{align}
where we have used 
\begin{align*}
& \sum_{i=1}^k C_k^i |\big\la (\hat{u}^0_{p,\t})^{4} \pa_x^{k+1} \bfq_{yy} ,\, \pa_x^i\bar{u}^0_{p} \pa_x^{k-i} \bfq_{yy}w^2\big\ra| \nonumber\\
&=\sum_{i=1}^k C_k^i  |\big\la (\hat{u}^0_{p,\t})^{4}  \pa_x^i\bar{u}^0_{p}  \pa_x^{k+1} \bfq_{y} ,\, \pa_x^{k-i} \bfq_{yyy}w^2\big\ra| + \sum_{i=1}^k C_k^i |\big\la \big\{(\hat{u}^0_{p,\t})^{4}\pa_x^i\bar{u}^0_{p} w^2\big\}_y \pa_x^{k+1} \bfq_{y} ,\,  \pa_x^{k-i} \bfq_{yy} \big\ra| \nonumber\\
&\geq \f18 \|(\hat{u}^0_{p,\t})^{3}\pa_x^{k+1} \bfq_yw\|^2 + C \|\bfq_y\|^2_{\MBX^k_{\theta, w}} + C  \sum_{i=0}^{k-1}\|(\hat{u}^0_{p,\t})^{2} \pa_x^i \bfq_{yyy} w\|^2,
\end{align*}
and
\begin{align*}
\sum_{i=0}^{k-1}\|(\hat{u}^0_{p,\t})^2 \pa_x^{i} \bfq_{yyy}w\|
&\lesssim  \|\bfq_y\|_{\MBX^k_{\theta, w}} + \sum_{i=0}^{k-1}\|\pa_x^{i}\psi_{yyy}w\| \lesssim  \|\bfq_y\|_{\MBX^k_{\theta, w}} + \sum_{i=0}^{k-1} \|\pa_x^i R w\|.
\end{align*}
 
\smallskip

For the second term on RHS of \eqref{6.269E}, it is clear that
\begin{align}\label{6.274-0E}
&|\big\la (\hat{u}^0_{p,\t})^{4} \pa_x^{k+1} \bfq_{yy},\, \pa_x^k(2\bar{u}^0_{py} \bfq_{y}) w^2\big\ra| \nonumber\\
&\leq  |\big\la  (\hat{u}^0_{p,\t})^{4} \pa_x^{k+1} \bfq_{y},\, \pa_x^k(2\bar{u}^0_{py} \bfq_{y})_y w^2\big\ra| + |\big\la \big\{ (\hat{u}^0_{p,\t})^{4} w^2\big\}_y \pa_x^{k+1} \bfq_{y},\, \pa_x^k(2\bar{u}^0_{py} \bfq_{y}) \big\ra| \nonumber\\
&\leq   \f18 \| (\hat{u}^0_{p,\t})^{3}\pa_x^{k+1} \bfq_yw\|^2 + C \|\bfq_y \|^2_{\MBX^k_{\theta, w}}.
\end{align}
For 3th term on RHS of \eqref{6.269E}, one gets
\begin{align}\label{6.275-0E}
	&|\big\la  (\hat{u}^0_{p,\t})^{4} \pa_x^{k+1} \bfq_{yy},\, \pa_x^k(\bar{u}^0_{pyy} \bfq)w^2\big\ra|\nonumber\\
	&\leq |\big\la (\hat{u}^0_{p,\t})^{4} \pa_x^{k+1} \bfq_{y},\, \pa_x^k(\bar{u}^0_{pyy} \bfq)_y w^2\big\ra| + |\big\la \big\{ (\hat{u}^0_{p,\t})^{4}w^2\big\}_y \pa_x^{k+1} \bfq_{y},\, \pa_x^k(\bar{u}^0_{pyy} \bfq) \big\ra| \nonumber\\
	&\leq  \f18 \| (\hat{u}^0_{p,\t})^{3}\pa_x^{k+1} \bfq_yw\|^2 + C \|\bfq_y \|^2_{\MBX^k_{\theta, w}}.
\end{align}
 
For the last two terms on RHS of \eqref{6.269E}, we have 
\begin{align}\label{6.282E}
 \big\la \big\{ (\hat{u}^0_{p,\t})^{4}w^2\big\}_y \pa_x^{k+1} \bfq_y,\, \pa_x^k\psi_{yy} \big\ra
&\leq  \f18 \| (\hat{u}^0_{p,\t})^{3}\pa_x^{k+1} \bfq_yw\|^2 + C \|\bfq_y \|^2_{\MBX^k_{\theta, w}},
\end{align}
where we  have used \eqref{5.95-1}.

Substituting \eqref{6.269E}-\eqref{6.282E} into \eqref{5.53E}, one 
concludes \eqref{5.53-0E}. Therefore the proof of Lemma \ref{lem4.5} is completed. $\hfill\Box$

\medskip

\subsubsection{Estimates on temperature $\mathbf{\Theta}$} Applying $\pa_x^k$ to  $\eqref{5.115}_2$, one has 
\begin{align}\label{5.103}
	&2\bar{\rho}^0_p \bar{u}^0_{p,\theta} \pa_x^{k+1} \mathbf{\Theta} + 2\bar{\rho}^0_p  \bar{v}^0_p \pa_y\pa_x^k \mathbf{\Theta} - \kappa \pa^2_{y}\pa_x^k \mathbf{\Theta} + 2\pa_x^k\big(\pa_x\bar{T}^0_p\, \psi_y\big) - 2\pa_x^k\big(\pa_y\bar{T}^0_p\,\psi_x\big)\nonumber\\
	&=\pa_x^kH - 2\sum_{i=1}^{k}  \big\{\pa_x^i\big(\bar{\rho}^0_p \bar{u}^0_{p,\theta}\big)\, \pa_x^{k+1-i} \mathbf{\Theta} +   \pa_x^i\big(\bar{\rho}^0_p  \bar{v}^0_p\big) \pa_y\pa_x^{k-i} \mathbf{\Theta}\big\} .
\end{align}

\begin{lemma}\label{lem4.6}
It holds that 
\begin{align}\label{6.65E}
\|\mathbf{\Theta}\|^2_{\MBY^{n-1}_{\t, w}} &\lesssim \sum_{i=0}^{n-1}\|\sqrt{\bar{u}^0_{p,\t}} \pa_x^k\mathbf{\Theta} w\|^2_{x=0} + L^{\f16-} \Big\{ \|\mathbf{\Theta}\|^2_{\MBY^{n-1}_{\t, w}} +\|\bfq_yw\|^2_{\MX^n_{\theta}} + \sum_{k=0}^{n-1} \|\pa_x^k H w\|^2\Big\}.
\end{align}
\end{lemma}

\noindent{\bf Proof.}  Let $0\leq k\leq n-1$. 
Multiplying $\eqref{5.103}$ by $\pa_x^k\mathbf{\Theta} w^2$, we obtain 
\begin{align}\label{6.67E}
&\|\sqrt{\bar{\rho}^0_p \bar{u}^0_{p,\t}} \pa_x^k\mathbf{\Theta} w\|^2_{x=x_0} + \kappa \|\pa_{y}\pa_x^k\mathbf{\Theta} w\|^2 \nonumber\\
&\lesssim \|\sqrt{\bar{\rho}^0_p \bar{u}^0_p} \pa_x^k\mathbf{\Theta} w\|^2_{x=0} + L^{\f16-} \Big\{ \|\mathbf{\Theta}\|^2_{\MBY^{k}_{\t, w}} + \sum_{i=0}^{k+1}\|\pa_x^{i}\psi_y w\|^2 + \|\pa_x^k H w\|^2\Big\},
\end{align}
where we have used 
\begin{align*}
|\big\la 2\pa_x^k\big(\pa_x\bar{T}^0_p\, \psi_y\big),\, \pa_x^k\mathbf{\Theta} w^2\big\ra|
	&\lesssim  L^{\f16-} \|T_p\|_{\MBY^k_{\t, w}} \sum_{i=0}^{k}\|\pa_x^{i}\psi_y w\|,
\end{align*}
and
\begin{align*}
\sum_{i=0}^{k} |\big\la 2\pa_x^k\big(\pa_y\bar{T}^0_p\,\psi_x\big),\, \pa_x^k\mathbf{\Theta} w^2\big\ra|
&\lesssim L^{\f16-} \|\mathbf{\Theta}\|_{\MBY^k_{\t, w}}  \sum_{i=0}^{k}\|\pa_x^i\pa_y\bar{T}^0_p w\|_{L^\infty_xL^2_y} \|\pa_x^{k+1-i}\psi\|_{L^2_xL^\infty} \nonumber\\
&\lesssim L^{\f16-} \|T_p\|_{\MBY^k_{\t, w}} \sum_{i=1}^{k+1}\|\pa_x^{i}\psi_y w\|.
\end{align*}
Thus it is clear to conclude \eqref{6.65E} from \eqref{6.67E} and \eqref{5.95}. Therefore the proof of Lemma \ref{lem4.6} is completed. $\hfill\Box$

\medskip

Noting the definition of $R$ in \eqref{5.76}, we need to control $\|\bar{u}^0_p \pa_x^n \pa_y^2 \mathbf{\Theta} w\|$. To control the highest order derivative estimate on  $\pa_x^n\mathbf{\Theta}$, we have to deal with the term $\pa_x^{n+1}\psi$ which requires high space weight. 
Similarly as in section \ref{sec2.4.2},  we   define a pseudo entropy
\begin{align}\label{5.109}
	\mathbf{S}_k:=\pa_x^{k}\mathbf{\Theta}  -  \pa_y\bar{T}^0_p \pa_x^{k} \bfq\quad \Longleftrightarrow \quad \pa_x^{k}\mathbf{\Theta} = \mathbf{S}_k +  \pa_y\bar{T}^0_p \pa_x^{k} \bfq.
\end{align}
Noting
\begin{align*}
\bar{\rho}^0_p \bar{u}^0_{p,\theta} \pa_x^{k+1}\mathbf{\Theta} - \pa_y\bar{T}^0_p\, \pa_x^{k+1}\psi 
&= \hat{u}^0_{p,\theta} \big\{\pa_x^{k}\mathbf{\Theta}  -  \pa_y\bar{T}^0_p \pa_x^{k} \bfq\big\}_x +  \hat{u}^0_{p,\theta}  \pa_{xy}\bar{T}^0_p \,  \pa_x^{k}\bfq  \nonumber\\
&\quad  - \pa_y\bar{T}^0_p \sum_{i=1}^{k+1} C_{k+1}^i \pa_x^i \hat{u}^0_{p,\t}\, \pa_x^{k+1-i}\bfq,
\end{align*}
which, together with \eqref{5.103}, yields that
\begin{align}\label{5.110}
	&2\hat{u}^0_{p,\theta} \pa_x\mathbf{S}_k + 2\bar{\rho}^0_p  \bar{v}^0_p \pa_y\mathbf{S}_k - \kappa \pa^2_{y} \mathbf{S}_k \nonumber\\
	&= \kappa \pa^2_{y}\big( \pa_y\bar{T}^0_p \pa_x^{k} \bfq\big) + 2\bar{\rho}^0_p  \bar{v}^0_p \pa_y\big(\pa_y\bar{T}^0_p \pa_x^{k} \bfq\big) - 2\pa_y\bar{T}^0_p \sum_{i=1}^{k+1} C_{k+1}^i \pa_x^i \hat{u}^0_{p,\t}\, \pa_x^{k+1-i}\bfq\nonumber\\
	&\quad  + 2\hat{u}^0_{p,\theta}  \pa_{xy}\bar{T}^0_p \,  \pa_x^{k}\bfq   - 2\pa_x^k\big(\pa_x\bar{T}^0_p\, \psi_y\big)  + \sum_{i=1}^{k} 2C_k^i \pa_x^i\pa_y\bar{T}^0_p\, \pa_x^{k+1-i}\psi  \nonumber\\
	& \quad + \pa_x^kH - 2\sum_{i=1}^{k}  \big\{\pa_x^i \hat{u}^0_{p,\theta} \, \pa_x^{k+1-i} \mathbf{\Theta} +   \pa_x^i\big(\bar{\rho}^0_p  \bar{v}^0_p\big) \pa_y\pa_x^{k-i}\mathbf{\Theta} \big\}.
\end{align}
It follows from \eqref{4.25E}  and \eqref{4.38E} that 
\begin{align}\label{4.84E}
\begin{split}
\pa_y \mathbf{S}_k \big|_{y=0}=0, \quad \mbox{for NBC},\\
\mathbf{S}_k\big|_{y=0}=0.\quad \mbox{for DBC}.
\end{split}
\end{align}
\begin{lemma}\label{lem4.7}
It holds that 
\begin{align}\label{5.71E}
\|\mathbf{S}_{n}\|^2_{\MBY_{\theta, w}} \lesssim \|\sqrt{\bar{u}^0_{p,\theta}} \mathbf{S}_n w\|^2_{x=0} +  L^{\f16-} \Big\{\|\mathbf{S}_n\|^2_{\MBY_{\t, w}} + \|\bfq_y\|^2_{\MBX^n_{\t, w}} + \|\pa_x^n H w\|^2\Big\}.
\end{align}
\end{lemma}

\noindent{\bf Proof.} 
Multiplying \eqref{5.110} by $\mathbf{S}_k w^2$, we obtain
\begin{align}\label{5.72E}
&\|\sqrt{\hat{u}^0_{p,\theta}} \mathbf{S}_kw\|^2_{x=x_0} + \f34\kappa\|\pa_y \mathbf{S}_k w\|^2 \nonumber\\
&\lesssim \|\sqrt{\hat{u}^0_{p,\theta}} \mathbf{S}_kw\|^2_{x=0}  + L^{\f16-} \Big\{\|\mathbf{S}_k\|^2_{\MBY_{\t, w}} + \|\bfq_y\|^2_{\MBX^k_{\t, w}} + \|\pa_x^k H w\|^2\Big\},
\end{align}
where we have used 
\begin{align*}
\big\la \pa^2_{y}\big(\pa_y\bar{T}^0_p \pa_x^{k} \bfq\big) ,\, \mathbf{S}_k w^2\big\ra 
&= \big\la \pa_{y}\big(\pa_y\bar{T}^0_p \pa_x^{k} \bfq\big) ,\, \pa_y\mathbf{S}_k w^2\big\ra + \big\la \pa_{y}\big( \pa_y\bar{T}^0_p \pa_x^{k} \bfq\big) ,\, 2\mathbf{S}_k ww_y\big\ra \nonumber\\
&\leq \f18 \|\pa_y \mathbf{S}_k w\|^2 + L^{\f16-} \|\mathbf{S}_k\|^2_{\MBY_{\t, w}} + L^{\f16-} \|\bfq_y\|^2_{\MBX^k_{\t, w}},
\end{align*}
and
\begin{align*}
\sum_{i=1}^{k} 2C_k^i |\big\la \pa_x^i\pa_y\bar{T}^0_p\, \pa_x^{k+1-i}\psi ,\, \mathbf{S}_k w^2\big\ra|
 &\lesssim L^{\f16-}\|\mathbf{S}_k\|_{\MBY_{\t, w}} \sum_{i=1}^{k} \|\pa_x^{i} \psi \|_{L^2_xL^\infty_y}
\lesssim  L^{\f16-}\|\mathbf{S}_k\|_{\MBY_{\t, w}} \sum_{i=1}^{k} \|\pa_x^{i} \psi_y \la y\ra \|\nonumber\\
&\lesssim L^{\f16-}\|\mathbf{S}_k\|^2_{\MBY_{\t, w}}  + L^{\f16-} \|\bfq_y\|^2_{\MBX^k_{\t, w}}.
\end{align*}
Taking $k=n$ in \eqref{5.72E}, we conclude \eqref{5.71E}.
Therefore the proof of Lemma \ref{lem4.7} is completed. $\hfill\Box$

\smallskip

\begin{lemma}[Elliptic estimate]\label{lem4.8}
It holds that 
\begin{align}\label{5.75E}
&\|\bar{u}^0_{p,\theta} \pa^2_{y} \mathbf{S}_nw\|^2 +  \|(\bar{u}^0_{p,\theta})^2 \pa_x\mathbf{S}_nw\|^2 + \|  (\bar{u}^0_{p,\theta})^{\f32}\pa_{y} \mathbf{S}_n\|^2_{x=x_0} \nonumber\\
&\lesssim \|(\bar{u}^0_{p,\theta})^{\f32}\pa_{y} \mathbf{S}_n\|^2_{x=0} +  \|\mathbf{S}_{n}\|^2_{\MBY_{\theta, w}} + \|\bfq_y\|^2_{\MBX^n_{\t, w}} + \|\pa_x^n H w\|^2.
\end{align}
\end{lemma}

\noindent{\bf Proof.} We square on both side of \eqref{5.110} with $k=n$ to obtain
\begin{align}\label{5.76E}
&\kappa ^2 \|\hat{u}^0_{p,\theta} \pa^2_{y} \mathbf{S}_nw\|^2 + 4\|(\hat{u}^0_{p,\theta})^2 \pa_x\mathbf{S}_nw\|^2 - 4\kappa\big\la \pa^2_{y} \mathbf{S}_n,\,   (\hat{u}^0_{p,\theta})^3 \pa_x\mathbf{S}_n w^2 \big\ra \nonumber\\
&\lesssim \|\mathbf{S}_{n}\|^2_{\MBY_{\theta, w}} + \|\bfq_y\|^2_{\MBX^n_{\t, w}} + \|\pa_x^n H w\|^2.
\end{align}
Integrating by parts, one has
\begin{align*}
&- 4\kappa\big\la \pa^2_{y} \mathbf{S}_n,\,  (\hat{u}^0_{p,\theta})^3 \pa_x\mathbf{S}_n\big\ra  = 4\kappa\big\la \pa_{y} \mathbf{S}_n,\,  (\hat{u}^0_{p,\theta})^3  \pa_{xy}\mathbf{S}_n w^2\big\ra
+  4\kappa\big\la \pa_{y} \mathbf{S}_n,\, \big\{(\hat{u}^0_{p,\theta})^3 w^2\big\}_y \pa_x\mathbf{S}_n\big\ra\nonumber\\
&\geq 2\kappa \|\sqrt{(\hat{u}^0_{p,\theta})^3}\pa_{y} \mathbf{S}_n\|^2_{x=x_0} -  2\kappa \|\sqrt{ (\hat{u}^0_{p,\theta})^3}\pa_{y} \mathbf{S}_n\|^2_{x=0} - \f18 \|(\bar{u}^0_{p,\theta})^2 \pa_x\mathbf{S}_nw\|^2  - C \|\mathbf{S}_{n}\|^2_{\MBY_{\theta, w}},
\end{align*}
which, together with \eqref{5.76}, concludes \eqref{5.75E}. Therefore the proof of Lemma \ref{lem4.8} is completed. $\hfill\Box$

\subsubsection{Conclusion} 
\begin{lemma}\label{lem4.10}
Assume \eqref{0.17}-\eqref{0.18} and high order compatibility conditions as in Section \ref{sec3.2}. Let $n\leq \mathfrak{n}_0+2$. There is a suitably small constant $L>0$ such that the solution of \eqref{4.22E} and \eqref{4.25E} exists in $x\in[0,L]$ satisfying 
\begin{align}\label{4.88E-2}
\sum_{2i+j\leq 2n} \|(\pa^i_x\pa_y^{j}\bfu, \pa^i_x\pa_y^{j}\mathbf{\Theta}) w\|^2& \lesssim C(\tilde{\bf I}[\tilde{u}_b, \tilde{T}_b, \tilde{u}_0, \tilde{T}_0], \hat{\bf I}[G]).
\end{align}
\end{lemma}

\noindent{\bf Proof.} {\it Step 1.} Recall \eqref{3.1}-\eqref{3.2}, \eqref{4.20E} and $n\leq \mathfrak{n}_0+2$, then we have from  \eqref{0.17}-\eqref{0.18} that
\begin{align}\label{4.88E-1}
&\sum_{k=0}^{n} \|\bar{u}^0_{p,\t}\pa_x^k \bfq_y w\|^2_{x=0} + \sum_{i=0}^{n-1}\|\pa_x^{i}\psi_{yy}w\|^2_{x=0}  +  \sum_{i=0}^{n-1}\|\sqrt{\bar{u}^0_{p,\t}} \pa_x^k\mathbf{\Theta} w\|^2_{x=0}  + \sum_{i=0}^{n-1}\|\pa_x^{i} \bfq_{yy}w\bar{\chi}\|^2_{x=0}\nonumber\\
	&\quad + \sum_{k=0}^{n}  \|(\bar{u}^0_{p,\theta})^{\f52}\pa_x^k \bfq_{yy}w\|^2_{x=0}  +  \|\sqrt{\bar{u}^0_{p,\theta}} \mathbf{S}_n w\|^2_{x=0}  + \|(\bar{u}^0_{p,\theta})^{\f32}\pa_{y} \mathbf{S}_n\|^2_{x=0}  \lesssim C(\tilde{\bf I}, \hat{\bf I}),
\end{align}

{\it Step 2.} It follows from \eqref{5.109} that 
\begin{align}\label{4.89E}
\begin{split}
\|\pa_x^{n}\mathbf{\Theta} w\| &\lesssim  L^{\f16-}\|\mathbf{S}_n\|_{\MBY_{\t, w}} + L^{\f16-}\|\bfq_y \|_{\MBX^{n}_{\t, w}},\\
\|\pa_x^{n}\mathbf{\Theta}_y w\| &\lesssim  \|\mathbf{S}_n\|_{\MBY_{\t, w}} +  L^{\f16-}\|\bfq_y \|_{\MBX^{n}_{\t, w}},\\
\|\bar{u}^0_p \pa_x^{n}\mathbf{\Theta}_{yy} w\| & \lesssim  \|\bar{u}^0_p \pa^2_y\mathbf{S}_n\| + \|\bfq_y\|_{\MBX^{n}_{\t, w}}.
\end{split}
\end{align}
Then it is clear that
\begin{align}\label{4.90E}
\sum_{i=0}^{n} \|(\pa_x^i\mathbf{\Theta}, \pa_x^i\mathbf{\Theta}_y) w\| 
&\lesssim \|\mathbf{\Theta}\|_{\MBY^{n-1}_{\t, w}} + \|\mathbf{S}_n\|_{\MBY_{\t, w}} + L^{\f16-}\|\bfq_y \|_{\MBX^{n}_{\t, w}}.
\end{align}

We have from  \eqref{4.24E}  that 
\begin{align}\label{4.91E}
\sum_{i=0}^{n} \|\pa_x^i H w\|&\lesssim \|(\pa_x^{n+1}G, \pa_x^{n}\pa_yG) w \la y\ra \|+ \sum_{i=0}^{n}\big\{\|(\pa_x^i\psi_{y}, \pa_x^i\psi_{yy}) w\| + \|(\pa_x^i\mathbf{\Theta}, \pa_x^i\mathbf{\Theta}_y) w\|\big\}\nonumber\\
&\lesssim \hat{\bf I} + \|\mathbf{\Theta}\|_{\MBY^{n-1}_{\t, w}} + \|\mathbf{S}_n \|_{\MBY_{\t, w}} + \|\bfq_y \|_{\MBX^{n}_{\t, w}},
\end{align}
where we have used \eqref{4.90E}, and  \eqref{5.95}-\eqref{5.95-1}.

\smallskip

For later use, we note from \eqref{5.103} that 
\begin{align}\label{4.92E}
\sum_{i=0}^{n-1} \|\pa_y^2 \pa^i_x \mathbf{\Theta} w\| 
&\lesssim \sum_{i=0}^{n} \|(\pa_x^i\mathbf{\Theta}, \pa_x^i\mathbf{\Theta}_y) w\| + \sum_{i=0}^{n-1} \|(\pa_x^i\psi_{y}, \pa_x^i H) w\|\nonumber\\
&\lesssim \|\mathbf{\Theta}\|_{\MBY^{n-1}_{\t, w}} + \|\mathbf{S}_n\|_{\MBY_{\t, w}} + \|\bfq_y \|_{\MBX^{n}_{\t, w}}.
\end{align}
Then it follows from \eqref{5.76} that 
\begin{align}\label{4.93E}
\sum_{i=0}^{n} \|\pa_x^i R w\|
&\lesssim \sum_{i=0}^{n}\big\{\|(\pa_x^i\psi_{y} ,\pa_x^i\psi_{yy} )w\| + \|(\pa_x^i\mathbf{\Theta}, \pa_x^i\mathbf{\Theta}_y)  w\| + \|\bar{u}^0_p \pa_x^i \mathbf{\Theta}_{yy} w\|\big\}\nonumber\\
&\quad + \|(\pa_x^{n+1}G, \pa_x^{n}\pa_yG) w \la y\ra \| \nonumber\\
&\lesssim  \hat{\bf I} + \|\bfq_y\|_{\MBX^n_{\t, w}} + \|\mathbf{\Theta}\|_{\MBY^{n-1}_{\t, w}}+ \|\mathbf{S}_{n}\|_{\MBY_{\theta, w}}+ \|\bar{u}^0_p \pa^2_y\mathbf{S}_n\|,
\end{align}
where we have used \eqref{5.95}-\eqref{5.95-1}, \eqref{4.90E}, \eqref{4.92E} and $\eqref{4.89E}_3$.

\smallskip

 It follows from \eqref{5.78} and \eqref{4.93E} that
\begin{align}\label{4.95E}
	\|\bfq_y\|^2_{\MBX^n_{\t, w}}
&\lesssim C(\tilde{\bf I}, \hat{\bf I})  +  L^{\f16-}\Big\{\|\bfq_{y}\|^2_{\MBX^n_{\theta, w}}  + \|\mathbf{\Theta}\|^2_{\MBY^{n-1}_{\t, w}}+ \|\mathbf{S}_{n}\|^2_{\MBY_{\theta, w}} \Big\}
\end{align}
where we have also used \eqref{5.75E}.

We have from \eqref{6.65E} and  \eqref{5.71E} that 
\begin{align}\label{4.96E}
\|\mathbf{\Theta}\|^2_{\MBY^{n-1}_{\t, w}} + \|\mathbf{S}_{n}\|^2_{\MBY_{\theta, w}} 
& \lesssim C(\tilde{\bf I}, \hat{\bf I})  + L^{\f16-} \Big\{\|\mathbf{\Theta} w\|^2_{\MBY^{n-1}_{\t}} + \|\mathbf{S}_n \|^2_{\MBY_{\t, w}} +\|\bfq_y\|^2_{\MBX^n_{\t, w}} \Big\}
\end{align}

Combining \eqref{4.95E}-\eqref{4.96E}, and taking $L>0$ suitably small (independent of $\theta\in(0,1]$), then we get 
\begin{align}\label{4.98E}
\|\mathbf{\Theta}\|^2_{\MBY^{n-1}_{\t, w}} + \|\mathbf{S}_n\|^2_{\MBY_{\t, w}} +\|\bfq_y\|^2_{\MBX^n_{\t, w}} &\lesssim C(\tilde{\bf I}, \hat{\bf I}). 
\end{align}
With above uniform estimate and Lemma \ref{lem4.2}, we can extend the local solution $(\bfu, \bfv, \mathbf{\Theta})$ of \eqref{5.114} to $x\in [0,L]$.

\medskip

Moreover, it follows from \eqref{4.98E}, \eqref{5.95-00} and Lemmas \ref{lem5.9} \& \ref{lem4.5} that 
\begin{align}
&\sum_{i=0}^n \Big\{\|(\bar{u}^0_{p,\theta})^{\f52} \pa_x^n \bfq_{yy}w\|^2_{x=x_0} +  \|(\bar{u}^0_{p,\theta})^{3}\pa_x^{i+1} \bfq_yw\|^2 + \|\bar{u}^0_{p,\theta} \pa_x^i\psi_{yyy}w\|^2  \Big\} \nonumber\\
&\quad + \|\bar{u}^0_{p,\theta} \pa^2_{y} \mathbf{S}_nw\|^2 +  \|(\bar{u}^0_{p,\theta})^2 \pa_x\mathbf{S}_nw\|^2 + \|  (\bar{u}^0_{p,\theta})^{\f32}\pa_{y} \mathbf{S}_n\|^2_{x=x_0}  \lesssim  C(\tilde{\bf I}, \hat{\bf I}),
\end{align}
and
\begin{align}\label{4.100E}
\begin{split}
L^{-\f18}	\sum_{k=0}^n\|\pa_x^k\psi_yw\|^2 + \sum_{k=0}^n\|\pa_x^k\psi_{yy}w\|^2 + \sum_{i=0}^{n-1}\|\pa_x^{i}\psi_{yyy}w\|^2 \lesssim C(\tilde{\bf I}, \hat{\bf I}),\\
\sum_{i=0}^{n-1}\|\pa_x^{i}\psi_{yy}w\|^2_{x=x_0}  \lesssim C(\tilde{\bf I}, \hat{\bf I}) +  \sum_{i=0}^{n-1}\|\pa_x^{i}\psi_{yy}w\|^2_{x=0},
\end{split}
\end{align}

Hence, noting $\psi_y=\bfu$, we have from \eqref{4.90E}, \eqref{4.92E}, \eqref{4.100E}  that 
\begin{align}\label{4.100E-1}
	&\sum_{i=0}^n\big\{\|(\pa_x^i \bfu, \pa_x^i \bfu_y,\, \pa_x^i\mathbf{\Theta}, \pa_x^i\mathbf{\Theta}_y) w\|^2\big\} + \sum_{i=0}^{n-1} \| (\pa_x^i\pa_y^2 \bfu, \pa^i_x\pa_y^2 \mathbf{\Theta})  w\|^2 
	\lesssim C(\tilde{\bf I}, \hat{\bf I}).
\end{align}

\medskip

We have from  $\eqref{5.114}_2$ that 
\begin{align}\label{4.101E}
\sum_{i=0}^{n-1} \|\pa_x^i \pa_y^3\mathbf{\Theta} w\|^2 &\lesssim \sum_{i=0}^{n} \|(\pa_x^i\mathbf{\Theta}, \pa_x^i\mathbf{\Theta}_y) w\|^2 + \sum_{i=0}^{n-1} \|\pa_y^2 \pa^i_x \mathbf{\Theta} w\| + \sum_{i=0}^{n}  \|\pa_x^i\psi_{y}w\|^2 \nonumber\\
&\quad  + \sum_{i=0}^{n-1}  \|(\pa_x^i\psi_{yy}, \pa_x^{i}\psi_{yyy})w\|^2 + \sum_{i=0}^{n-1} \|\pa_x^i H_y w\|^2  
\lesssim C(\tilde{\bf I}, \hat{\bf I}). 
\end{align}
Then it follows from  \eqref{5.76}, \eqref{4.90E}, \eqref{4.92E}, \eqref{4.98E} and \eqref{4.101E} that 
\begin{align}\label{4.102E}
\sum_{i=0}^{n-1}\|\pa_x^{i}R_y w\|^2
&\lesssim C(\tilde{\bf I}, \hat{\bf I}) + \sum_{i=0}^{n-1} \big\{ \|(\pa_x^i\psi_{y}, \pa_x^i\psi_{yy}, \pa_x^{i}\psi_{yyy})w\|^2   +  \|\bar{u}^0_p\pa_y^3\mathbf{\Theta} w \|^2\big\}\nonumber\\
&\quad +  \|\mathbf{\Theta}\|^2_{\MBY^{n-1}_{\t, w}} + \|\mathbf{S}_n\|^2_{\MBY_{\t, w}} + \|\bfq_y\|^2_{\MBX^{n}_{\t, w}} \lesssim C(\tilde{\bf I}, \hat{\bf I}),
\end{align}
which, together with \eqref{5.96}, yields that 
\begin{align}\label{4.103E}
\sum_{i=0}^{n-1}\|\pa_x^{i}\psi_{yyyy}w\|^2 &\lesssim C(\tilde{\bf I}, \hat{\bf I}).
\end{align}

Then we have from  \eqref{4.101E} and \eqref{4.103E} that 
\begin{align*}
 \sum_{i=0}^{n-1} \|(\pa_x^i \pa^3_y\bfu, \pa^i_x \pa_y^3\mathbf{\Theta} ) w\|^2   \lesssim C(\tilde{\bf I}, \hat{\bf I}).
\end{align*}
Similarly, by induction arguments, we can obtain
\begin{align}\label{4.105E}
\sum_{2i+j\leq 2n} \|(\pa^i_x\pa_y^{j}\bfu, \pa^i_x\pa_y^{j}\mathbf{\Theta}) w\|^2& \lesssim C(\tilde{\bf I}, \hat{\bf I}).
\end{align}

{\it Step 3.} With the solution of \eqref{5.114}  established in Step 2, we take the limit $\theta\to 0+$  to derive the existence of solution for \eqref{4.22E} and \eqref{4.25E} in $x\in [0,L]$. The uniform estimate \eqref{4.88E-2} follows directly  from \eqref{4.105E}. Therefore the proof of Lemma \ref{lem4.10} is completed. $\hfill\Box$

\medskip

 
\subsection{\bf Proof of Theorem \ref{thm2}.} 
Noting Lemma \ref{lem4.10}, and the equivalence of \eqref{4.22E} and \eqref{3.6-10}, we have solved the original linear Prandtl layer problem \eqref{3.6-10}-\eqref{3.6-11}.

Recall $n\leq \mathfrak{n}_0 +2$, and $w_1:=(1+y)^{l_1}$. For $0\leq 2i+j\leq 2n-3$, it follows from  \eqref{4.88E-2} that 
\begin{align*}
\begin{split}
\|\pa_x^{n-2}\bfv \|_{L^\infty}  &\lesssim \|\pa_x^{n-2}\bfv(0,\cdot)\|_{L^\infty}  + \sqrt{L} C(\tilde{\bf I}, \hat{\bf I}), \\
\|(\pa_x^i \pa_y^{j}\bfu, \, \pa_x^i \pa_y^{j}\mathbf{\Theta})  w_1\|_{L^\infty} 
&\lesssim \|(\pa_x^i \pa_y^{j}\bfu, \, \pa_x^i \pa_y^{j}\Theta)(0,\cdot)  w_1\|_{L^\infty}  + \sqrt{L} C(\tilde{\bf I}, \hat{\bf I}),
\end{split}
\end{align*}
which, together with \eqref{3.2} and  \eqref{4.20E}, yields that
\begin{align*}
	&\|\pa_x^{n-2}v_p\|_{L^\infty} + \sum_{0\leq 2i+j\leq 2n-4} \|(\pa_x^i \pa_y^{j}u_p, \, \pa_x^i \pa_y^{j}T_p)  w_1\|_{L^\infty} \nonumber\\
	&\lesssim C({\bf I})\,  \Big\{ \sqrt{L} C(\tilde{\bf I}, \hat{\bf I})  +  \sum_{j=0}^{n-2}\|\pa_x^{j}v_p(0,\cdot)\|_{L^\infty_y} + \sum_{0\leq 2i+j\leq 2n-4}  \|(\pa_x^i \pa_y^{j} \tilde{u}_0, \, \pa_x^i \pa_y^{j} \tilde{T}_0)  w_1\|_{L^\infty}   \nonumber\\
	&\quad  + \sum_{i=0}^{n-2}\|(\pa_x^i\tilde{T}_b, \pa_x^i\tilde{u}_b)\|_{L^\infty} + \sum_{i=0}^{n-2}\|\pa_x^i G  w_1 \la y\ra \|_{L^\infty}\Big\}.
\end{align*}
Hence we conclude \eqref{1.10}. Therefore the proof of Theorem \ref{thm3} is completed. $\hfill\Box$

\medskip

\section{Construction of the Compressible Euler Layer} \label{sec4}
In this section,  we aim to establish for the existence of solution for the linear compressible Euler layer problem \eqref{2.23}-\eqref{2.4-1}. Noting the uniform subsonic condition \eqref{1.3-6} and Lemma \ref{lem2.2}, we first consider  the existence of an elliptic boundary value problem in Section \ref{sec4.1}.

\subsection{Existence of  an elliptic problem}\label{sec4.1}
We study the following elliptic boundary value problem:
\begin{align}\label{5.1}
	\begin{cases}
		\dis -\pa_Y (a_{11} \pa_{Y}V) -\pa_{x} (a_{22} \pa_x V)   + \pa_Y(b_1 V)=\pa_Y F,\\[2mm]
		\dis V|_{Y=0}=-v_p(x,0),\quad V |_{x=0}=V_{0}(Y),\quad V |_{x=L}=V_{L}(Y),
	\end{cases}
\end{align}
where $v_p(x,0), V_{0}(Y), V_{L}(Y)$ are given smooth data, and 
\begin{align}
	\begin{split}
		\dis a_{11}(Y)&:=\f{\rho_e^0 T_e^0 u^0_e }{2T^0_e - |u_e^0|^2} >0,  \qquad a_{22}(Y):=\f12 \rho_e^0 u_e^0>0,\\
		\dis b_1(Y)&:=\f{T_e^0[\rho_e^0  \pa_Yu^0_e  - u^0_e \pa_{Y}\rho_e^0  -  2\f{\rho_e^0}{u_e^0}\pa_Y T^0_e]}{2T^0_e - |u_e^0|^2} + \f{\rho_e^0}{u_e^0}\pa_Y T^0_e.
	\end{split}
\end{align}
For the boundary data in \eqref{5.1}, we need the following compatibility conditions
\begin{align}
	V_0(0)=-v_p(0,0),\quad V_L(0)=-v_p(L,0).
\end{align}

For the case of $v_p(0,0)\neq 0$ and $v_p(L,0)\neq 0$, we define an auxiliary function
\begin{align}
	S_1(x,Y):&= (1-\f{x}{L}) \f{V_0(Y)}{v_p(0,0)} v_p(x,0) + \f{x}{L} \f{V_L(Y)}{v_p(L,0)} v_p(x,0)\nonumber\\
	&\equiv v_p(x,0) \f{V_0(Y)}{v_p(0,0)}  + \f{x}{L} \Big[\f{V_L(Y)}{v_p(L,0)} - \f{V_0(Y)}{v_p(0,0)}\Big]v_p(x,0).
\end{align}
We assume the following uniform-in-$L$ bound on $V_0(Y), V_L(Y)$:
\begin{align}\label{5.6}
	\begin{split}
	\Big|\pa_Y^k\Big(\f{V_0(Y)}{V_0(0)}, \f{V_L(Y)}{V_L(0)}\Big) \la Y\ra^{l}\Big| + \f{1}{L}\Big|\pa_Y^k\Big(\f{V_0(Y)}{V_0(0)}- \f{V_L(Y)}{V_L(0)}\Big) \la Y\ra^{l}\Big|&\lesssim 1,\,\, k=0,1,\cdots,
	\end{split}
\end{align}
where $l\gg1$.  Then it is follows from \eqref{5.6} that 
\begin{align}\label{4.6}
\begin{split}
\big|\pa_Y^j S_1\, \la Y\ra^{l}\big|&\lesssim |v_p(x,0)|,\quad \mbox{for}\,\, j\geq 0,\\
\big|\pa_x^i \pa_Y^j S_1\, \la Y\ra^{l}\big|&\lesssim |(\pa_x^i v_p, \pa_x^{i-1} v_p)(x,0)|,\quad \mbox{for}\,\, i\geq1,\,\, j\geq0.
\end{split}
\end{align}

\smallskip

For the case of $v_p(0,0)\equiv V_0(0)=0$ and $v_p(L,0)\equiv-V_L(0)\neq 0$, we need to reset the auxiliary function 
as
\begin{align}
S_1(x,Y):&= (1-\f{x}{L}) \big[V_0(Y) - \chi(Y)\,v_p(x,0)\big] + \f{x}{L} \f{V_L(Y)}{v_p(L,0)} v_p(x,0)\nonumber\\
&=(1-\f{x}{L}) V_0(Y) - \chi(Y)\,v_p(x,0)  - \f{x}{L} \Big\{\f{V_L(Y) }{V_L(0)}- \chi(Y)\Big\}v_p(x,0).
\end{align}
We assume 
\begin{align}\label{5.9-1}
\f{1}{L}\Big|\pa_Y^k\Big(\f{V_L(Y) }{V_L(0)} - \chi(Y)\Big)\, \la Y\ra^{l}\Big| + \f{1}{L} |\pa_Y^kV_0(Y) \,\la Y\ra^{l}|\lesssim 1.
\end{align}
which yields the estimates as in \eqref{4.6}.
Similarly, one can also modify the function $S_1$ if $v_p(L,0)=0$, the detail is omitted for simplicity of presentation. Therefore, under our assumption \eqref{5.6} and \eqref{5.9-1},  we always have \eqref{4.6}.

\medskip

We denote
\begin{align}
\begin{split}
\bar{V}:&=V-S_1,\\
F_b:&=\pa_Y (a_{11} \pa_{Y}S_1) + \pa_{x} (a_{22} \pa_x S_1)   - \pa_Y(b_1 S_1),
\end{split}
\end{align}
then it follows from \eqref{5.1} that
\begin{align}\label{5.7}
	\begin{cases}
		-\pa_Y (a_{11} \pa_{Y}\bar{V}) -\pa_{x} (a_{22} \pa_x \bar{V})   + \pa_Y(b_1 \bar{V}) 
		 =\pa_Y F + F_b,\\[2mm]
		\bar{V}|_{x=0,L}=0,\quad \bar{V}|_{Y=0}=0.
	\end{cases}
\end{align}


\begin{lemma}\label{lem5.1}
	There exists a unique weak solution $\bar{V}\in H^1_0$ of \eqref{5.7} satisfying
	\begin{align}\label{5.13-1}
		\|(\bar{V},\bar{V}_Y,\bar{V}_x)w\| \lesssim \|Fw\| +  \|v_p(\cdot,0)\|_{W^{1,\infty}_x}.
	\end{align}
\end{lemma}

\noindent{\bf Proof.} Multiplying \eqref{5.7} by $\varphi \in H^{1}_0$, one obtains the following bilinear form
\begin{align}\label{5.9}
	\mathfrak{B}(\bar{V},\varphi)
	&=\int_0^L\int_0^\infty \big\{\pa_YF + F_b\big\} \varphi dYdx,
\end{align}
where 
\begin{align*}
	\mathfrak{B}(\bar{V},\varphi):&=\int_0^L\int_0^\infty \big\{a_{11} \pa_{Y}\bar{V} \, \pa_Y \varphi  + a_{22} \pa_x \bar{V} \, \pa_x\varphi - b_1\bar{V}\, \pa_Y \varphi\big\} dydx.
\end{align*}

It is clear
\begin{align}\label{5.12}
\Big|\int_0^L\int_0^\infty \big\{\pa_YF + F_b\big\} \varphi dYdx\Big| \lesssim \|(F, S_1,  S_{1x}, S_{1Y})\|\cdot \|\varphi\|_{H^1_0}
\end{align}
Taking $\varphi=\bar{V}$ in \eqref{5.9}, we have 
\begin{align}\label{5.13}
	\mathfrak{B}(\bar{V},\bar{V})
	&=\int_0^L\int_0^\infty \Big\{a_{11} |\pa_{Y}\bar{V}|^2 + a_{22} |\pa_x \bar{V}|^2 - \f12 b_1 \pa_Y(|\bar{V}|^2)\Big\} dydx \nonumber\\
	&\geq c_0\|(\pa_Y\bar{V},\pa_x\bar{V})\|^2-C\|\bar{V}\|^2 \geq (c_0-CL^2) \|\bar{V}\|^2_{H^1_0}.
\end{align}
With \eqref{5.12}-\eqref{5.13}, applying the standard elliptic theory, we can obtain a unique weak solution $\bar{V}\in H^1_0$ of \eqref{5.9}. 

Taking $\varphi=\bar{V}w^2$ in \eqref{5.9}. It is clear to see
\begin{align}\label{5.15}
	&\Big|\int_0^L\int_0^\infty \big\{\pa_YF + F_b\big\} \bar{V} w^2 dYdx\Big|\nonumber\\
	&\leq \Big|\int_0^L\int_0^\infty \big[F + a_{11} \pa_{Y}S_1   - b_1 S_1\big] \big[\pa_{Y}\bar{V}  w^2 + 2 \bar{V}   w \pa_{Y}w \big] dYdx\Big| \nonumber\\
	&\quad + \Big|\int_0^L\int_0^\infty  a_{22} \pa_x S \pa_x\bar{V} w^2  dYdx\Big|\nonumber\\
	&\lesssim \|(\pa_x\bar{V},\pa_Y\bar{V})w\|\cdot \|(F, S_1,  S_{1x}, S_{1Y})w\|,
\end{align}
and
\begin{align}\label{5.16}
	\mathfrak{B}(\bar{V},\varphi):&=\int_0^L\int_0^\infty \big\{[a_{11} \pa_{Y}\bar{V} - b_1 \bar{V}] [\pa_Y \bar{V} w^2 + 2 \bar{V}  w\pa_Y w]  + a_{22} |\pa_x \bar{V}|^2 \, w^2 \big\} dydx\nonumber\\
	&\geq c_0 \|(\pa_x\bar{V}, \pa_Y\bar{V})w\|^2.
\end{align}
Then we have from \eqref{5.15}-\eqref{5.16} that
\begin{align*}
\|(\bar{V},\bar{V}_Y,\bar{V}_x)w\| \lesssim \|(F, S_1,  S_{1x}, S_{1Y})w\| \lesssim \|Fw\| +  \|v_p(\cdot,0)\|_{W^{1,\infty}_x},
\end{align*}
where we have used \eqref{4.6}. Thus \eqref{5.13-1} is proved. 
Therefore the proof of Lemma \ref{lem5.1} is completed. $\hfill\Box$

\medskip

%


\medskip

We rewrite $\eqref{5.7}_1$ as
\begin{align}\label{5.18}
	\begin{cases}
		- \big(a_{11} \bar{V}_{YY} -  [b_1 - \pa_Y a_{11}] \bar{V}_{Y}\big)_{Y}  - a_{22} \bar{V}_{Yxx}  - \pa_Y a_{22} \bar{V}_{xx}  + \big(\pa_Yb_1 \bar{V}\big)_{Y}
		= \big(\pa_{Y} F + F_b\big)_{Y},\\
		\big[a_{11} \bar{V}_{YY} + (b_1 - \pa_Y a_{11}) \bar{V}_{Y}\big]  \big|_{Y=0} 
		=-\pa_Y F - F_b.
	\end{cases}
\end{align}

\begin{lemma}\label{lem5.3}
It holds that 
\begin{align}\label{5.19}
\|(V,V_Y,V_x, V_{YY},V_{xY}, V_{xx})w\| \lesssim  \|(F, \pa_YF)w\| + \|v_p(\cdot,0)\|_{W^{2,\infty}_x}.
\end{align}
\end{lemma}

\noindent{\bf Proof.} Multiplying \eqref{5.18} by $\bar{V}_Y$, one has
\begin{align}\label{5.20}
	&-\big\la \big\{a_{11} \bar{V}_{YY} -  [b_1 - \pa_Y a_{11}] \bar{V}_{Y}+\pa_{Y} F + F_b\big\}_{Y},\, \bar{V}_Y \big\ra - \big\la a_{22} \bar{V}_{Yxx},\, \bar{V}_Y\big\ra\nonumber\\
	&= \big\la \pa_Y a_{22} \bar{V}_{xx},\, \bar{V}_Y\big\ra - \big\la \big(\pa_Yb_1 \bar{V}\big)_{Y},\, \bar{V}_Y\big\ra 
\end{align}
Integrating by parts in $Y$ and using $\eqref{5.18}_2$, one gets
\begin{align}\label{5.21}
	&-\big\la \big\{a_{11} \bar{V}_{YY} - [b_1 - \pa_Y a_{11}] \bar{V}_{Y}+\pa_{Y} F + F_b\big\}_{Y},\, \bar{V}_Y \big\ra\nonumber\\
	&= \big\la a_{11} \bar{V}_{YY} - [b_1 - \pa_Y a_{11}] \bar{V}_{Y}+\pa_{Y} F + F_b ,\, \bar{V}_{YY} \big\ra \nonumber\\
	&\geq c_0\|\bar{V}_{YY}\|^2 -C \|(\bar{V}_{Y}, \pa_YF, F_b)\|^2.
\end{align}
Integrating by parts in $x$ and using $\bar{V}\big|_{x=0,L}=0$, one obtains
\begin{align*}
	- \big\la a_{22} \bar{V}_{Yxx},\, \bar{V}_Y\big\ra= \big\la a_{22} \bar{V}_{Yx},\, \bar{V}_{Yx}\big\ra\geq c_0\|\bar{V}_{Yx}\|^2.
\end{align*}
For the rest terms in \eqref{5.20}, it is clear 
\begin{align}\label{5.23}
	& \big\la \pa_Y a_{22} \bar{V}_{xx},\, \bar{V}_Y\big\ra - \big\la \big(\pa_Yb_1 \bar{V}\big)_{Y},\, \bar{V}_Y\big\ra \nonumber\\
	&\lesssim \|\bar{V}_Y\| \cdot \|(\bar{V}, \bar{V}_Y, \bar{V}_{xx})\| \lesssim \|\bar{V}_Y\|\cdot \|(\bar{V}, \bar{V}_Y, \bar{V}_{YY}, \pa_YF,F_b)\|
\end{align}
where we have used $\eqref{5.7}_1$ to deal with $\bar{V}_{xx}$.

Substituting \eqref{5.21}-\eqref{5.23} into \eqref{5.20}, one has
\begin{align}\label{5.24}
	\|(\bar{V}_{YY},\bar{V}_{Yx})\|^2 \lesssim \|(\bar{V}, \bar{V}_{Y}, \pa_YF, F_b)\|^2 \lesssim  \|(\pa_YF, F_b, S_1,  S_{1x}, S_{1Y})\|^2.
\end{align}
It follows from $\eqref{5.7}_1$ and \eqref{5.24} that 
\begin{align*}
	\|\bar{V}_{xx}\|^2\lesssim \|(\pa_YF, F_b, S_1,  S_{1x}, S_{1Y})\|^2.
\end{align*}
Similarly, we can also establish the weighted estimate
\begin{align*}
\|(\bar{V}_{YY},\bar{V}_{Yx}, \bar{V}_{xx})w\| \lesssim  \|(\pa_YF, F_b, S_1,  S_{1x}, S_{1Y})w\|,
\end{align*}
which, together with \eqref{5.13-1}, yields that 
\begin{align*}
\|(V,V_Y,V_x, V_{YY},V_{xY}, V_{xx})w\| &\lesssim  \|(F, \pa_YF)w\| + \sum_{j=0}^2 \|\nabla^j S_1 w\|\nonumber\\
&\lesssim \|(F, \pa_YF)w\| + \|v_p(\cdot,0)\|_{W^{2,\infty}_x},
\end{align*}
where we have used \eqref{4.6} in the last inequality. Therefore the proof of Lemma \ref{lem5.3} is completed. $\hfill\Box$

\medskip

It is noted that 
\begin{align*}
	\bar{V}|_{x=0,L}=\bar{V}|_{Y=0}=0\quad \mbox{and}\quad 
	\big[a_{11} \bar{V}_{YY} - (b_1 - \pa_Y a_{11}) \bar{V}_{Y}\big]  \big|_{Y=0} 
	=-\pa_Y F - F_b.
\end{align*}
or equivalently
\begin{align}\label{5.26}
	\begin{split}
		&\pa^2_YV(0,Y)\equiv V_0''(Y), \qquad\qquad \pa^2_YV(L,Y)\equiv V_L''(Y),\\
		&a_{11}(0)\pa^2_Y V(x,0) - [b_1-\pa_Y a_{11}](0)  V_Y(x,0)  = \big\{ -a_{22} v_{pxx} + \pa_Yb_1 v_p - \pa_YF \big\}(x,0).
	\end{split}
\end{align}
Then, to get higher-order derivatives estimates,  we need to assume the following compatibility conditions on $V_0 ~ \&~ V_L$ so that 
\begin{align}\label{5.28}
	\big( \pa_Y F +  F_b\big)(0,0)=\big( \pa_Y F +  F_b\big)(L,0)=0.
\end{align}

Denote
\begin{align}\label{5.28-0}
	\begin{cases}
		W_0(Y):=a_{11}(Y) V''_{0}(Y) - [b_1 - \pa_Y a_{11}](Y) V'_{0}(Y),\\
		W_L(Y):=a_{11}(Y) V''_{L}(Y) - [b_1 - \pa_Y a_{11}](Y) V'_{L}(Y),
	\end{cases}
\end{align}
then the compatibility condition \eqref{5.28} is equivalent to the following
\begin{align}\label{5.28-1}
	\begin{cases}
		W_0(0)
		= - \big\{a_{22} v_{pxx} -\pa_Yb_1 v_p + \pa_YF \big\}(0,0),\\
		W_L(0)
		=- \big\{a_{22} v_{pxx} -\pa_Yb_1 v_p + \pa_YF \big\}(L,0),
	\end{cases}
\end{align}

Now we define
\begin{align}\label{5.35-0}
	W(x,Y)&:=a_{11}(Y) V_{YY}(x,Y) - [b_1 - \pa_Y a_{11}](Y) V_Y(x,Y).
\end{align}
If $W_0(0)\neq 0$ and $W_L(0)\neq 0$, we introduce a  new auxiliary function
\begin{align}\label{5.36-0}
	S_2(x,Y):&=\big(1-\f{x}{L}\big) \f{W_0(Y)}{W_0(0)} \big\{ -a_{22} v_{pxx} + \pa_Yb_1 v_p - \pa_YF \big\}(x,0) \nonumber\\
	&\quad + \f{x}{L} \f{W_L(Y)}{W_L(0)} \big\{ -a_{22} v_{pxx} + \pa_Yb_1 v_p - \pa_YF \big\}(x,0).
\end{align}
If $W_0(0)=0$, we modify \eqref{5.36-0} as
\begin{align}
	S_2(x,Y):&=\big(1-\f{x}{L}\big) \Big[W_0(Y) + \chi(Y)\big\{ -a_{22} v_{pxx} + \pa_Yb_1 v_p - \pa_YF \big\}(x,0)\Big] \nonumber\\
	&\quad + \f{x}{L} \f{W_L(Y)}{W_L(0)} \big\{ -a_{22} v_{pxx} + \pa_Yb_1 v_p - \pa_YF\big\}(x,0).
\end{align}
One can also define $S_2$ for other cases. Similarly as \eqref{5.6} and \eqref{5.9-1},  we need to assume some compatibility conditions on $W_0, W_L$ so that 
\begin{align}\label{5.38-0}
	\begin{split}
		\big|\pa_Y^j S_2\, \la Y\ra^{l}\big|&\lesssim |(\pa_x^2v_p, v_p)(x,0)|,\quad \mbox{for}\,\, j\geq 0,\\
		\big|\pa_x^i \pa_Y^j S_2\, \la Y\ra^{l}\big|&\lesssim |(\pa_x^{i+2} v_p, \pa_x^{i+1} v_p, \pa_x^iv_p,  \pa_x^{i-1} v_p)(x,0)|,\quad \mbox{for}\,\, i\geq1,\,\, j\geq0.
	\end{split}
\end{align}

Define
\begin{align*}
	\bar{W}:=W - S_2,
\end{align*}
then we have from  $\eqref{5.1}_1$ that
\begin{align}\label{5.33}
	\begin{cases}
		\dis -\bar{W}_{YY} -  \f{a_{22}}{a_{11}}\bar{W}_{xx}
		=\pa_Y^3 F  +2\pa_Ya_{22} V_{xxY} + \f{a_{22}}{a_{11}}[b_1-\pa_Ya_{11}]V_{xxY} + \pa_Y^2a_{22} V_{xx} \\
		\dis \hspace{3.5cm}   - \big(\pa_Yb_1 V\big)_{YY} + S_{2YY} + \f{a_{22}}{a_{11}}S_{2xx},\\
		\bar{W}|_{x=0,L}=\bar{W}|_{Y=0}=0.
	\end{cases}
\end{align}

\begin{lemma}\label{lem5.4}
	Let the compatibility conditions described above hold, then we have
	\begin{align}\label{5.38}
		\|(V_{YYY},V_{xYY},V_{xxY},V_{xxx})w\| 
		&\lesssim \|v_p(\cdot,0)\|_{W^{4,\infty}_x} + \sum_{j=0}^2 \|\nabla^j F w\|.
	\end{align}
\end{lemma}

\noindent{\bf Proof.} Multiplying \eqref{5.33} by $\bar{W}$, one obtains 
\begin{align*}
	\|(\bar{W}_Y,\bar{W}_x)\|^2 
	&\lesssim \|(S_{2YY}, S_{2xx}, V, V_Y V_{YY}, V_{xx}, \pa_Y^2F)\|^2,
\end{align*}
which yields immediately that 
\begin{align}\label{5.35}
	\|(V_{YYY}, V_{xYY})\| \lesssim \|(S_{2YY}, S_{2xx}, V, V_Y, V_{YY}, V_{xY},  V_{xx}, \pa_Y^2F)\|.
\end{align}

Applying $\pa_Y$ and $\pa_x$ to $\eqref{5.1}_1$ respectively , we can derive
\begin{align}\label{5.36}
	\begin{split}
		\|V_{xxY}\|&\lesssim \|V_{YYY}\| + \|(V_{YY}, V_{xx}, V, \pa^2_YF)\|,\\
		\|V_{xxx}\|&\lesssim \|V_{xYY}\| + \|(V_{YY}, V_{xY}, V_Y,V_x, \pa_{xY}F)\|.
	\end{split}
\end{align}
Thus it follows from \eqref{5.35}-\eqref{5.36} that
\begin{align}\label{5.34}
	\|(V_{YYY},V_{xYY},V_{xxY},V_{xxx})\|&\lesssim \|(S_{2YY}, S_{2xx}, V, V_Y, V_{YY}, V_{xY},  V_{xx}, \pa_Y^2F,\pa_{xY}F)\|\nonumber\\
	&\lesssim \|(S_1, S_{1x},  S_{1Y}, S_{1YY}, S_{1xY}, S_{1xx}, S_{2YY}, S_{2xx})\|\nonumber\\
	&\quad + \|(F, \pa_YF, \pa_Y^2F,\pa_{xY}F)\|\nonumber\\
	&\lesssim \sum_{j=0}^2 \|\nabla^j F\| +  \|v_p(\cdot,0)\|_{W^{4,\infty}_x}.
\end{align}
Similarly, we can also  derive the corresponding weighted estimate. Therefore the proof of Lemma \ref{lem5.4} is completed. $\hfill\Box$

\medskip

\begin{lemma}\label{lem5.5}
	It holds that
	\begin{align}\label{5.46-0}
		\|(V,V_Y, V_x) \la Y\ra^{l}\|^2_{L^\infty} &\lesssim \|v_p(\cdot,0)\|_{W^{1,\infty}_x} + L \|v_p(\cdot,0)\|_{W^{4,\infty}_x} + L \sum_{j=0}^2 \|\nabla^j F \la Y\ra^{l+1}\|_{L^\infty}.
	\end{align}
\end{lemma}

\noindent{\bf Proof.} Noting \eqref{6.162}, one has
\begin{align}
	|\la Y\ra^{l}\bar{V}(x,Y)|^2&\lesssim \int_0^\infty | \bar{V}_Y(x,Y)|^2\la Y\ra^{2l+1} dY \lesssim L\int_0^\infty \int_0^L \la Y\ra^{2l+1} |\bar{V}_{xY}(x,Y)|^2 dx dY\nonumber\\
	&\lesssim L \|(V_{xY}, S_{1xY}) \la Y\ra^{l+\f12}\|^2,
\end{align}
and
\begin{align}
	|\la Y\ra^{l}\bar{V}_Y(x,Y)|^2&\lesssim \int_0^\infty | \bar{V}_{YY}(x,Y)|^2\la Y\ra^{2l+1} dY \lesssim L\int_0^\infty \int_0^L \la Y\ra^{2l+1} |\bar{V}_{xYY}(x,Y)|^2 dx dY\nonumber\\
	&\lesssim  L \|(V_{xYY}, S_{1xYY}) \la Y\ra^{l+\f12}\|^2.
\end{align}
Due to $\bar{V}_Y(0,Y)=\bar{V}_Y(L,Y)=0$, there exists a point $x_*(Y)\in [0,L]$ so that $\bar{V}_{xY}(x_{*}(Y),Y)=0$. Then it holds that 
\begin{align}
	|\la Y\ra^{l}\bar{V}_{x}(x,Y)|^2 &\lesssim \int_0^\infty | \bar{V}_{xY}(x,Y)|^2\la Y\ra^{2l+1} dY \lesssim L\int_0^\infty \int_0^L \la Y\ra^{2l+1} |\bar{V}_{xxY}(x,Y)|^2 dx dY\nonumber\\
	&\lesssim L \|(V_{xxY}, S_{1xxY})\la Y\ra^{l+\f12}\|^2.
\end{align}
Then it is direct to derive 
\begin{align*}
	|(V,V_Y, V_x)(x,Y)\, \la Y\ra^{l}|^2 &\lesssim |(S_1, S_{1Y}, S_{1x})(x,Y)\, \la Y\ra^{l}|^2 +  L \|(S_{1xY}, S_{1xYY},  S_{1xxY}) \la Y\ra^{l+\f12}\|^2 \nonumber\\
	&\quad  + L \|(V_{xY}, V_{xYY}, V_{xxY}) \la Y\ra^{l+\f12}\|^2\nonumber\\
	&\lesssim \|v_p(\cdot,0)\|_{W^{1,\infty}_x} + L \|v_p(\cdot,0)\|_{W^{4,\infty}_x} + L \sum_{j=0}^2 \|\nabla^j F \la Y\ra^{l+2}\|_{L^\infty},
\end{align*}
which concludes \eqref{5.46-0}. Therefore the proof of Lemma \ref{lem5.5} is completed. $\hfill\Box$

\medskip

The compatibility conditions can be assumed to arbitrary order by iterating above
process, and thus we can obtain
\begin{proposition}\label{prop5.6}
	Assume the higher order compatibility conditions are satisfied at corners $(0,0)$ and $(L,0)$. Then the solution of \eqref{5.1} satisfies
	\begin{align}\label{5.45}
		\sum_{i=0}^{k} \|\nabla^i V\, \la Y\ra^{l}\|_{L^\infty} \lesssim \|v_p(\cdot,0)\|_{W^{k,\infty}_x} + L \|v_p(\cdot,0)\|_{W^{k+3,\infty}_x} + L \sum_{j=0}^{k+1} \|\nabla^j F \la Y\ra^{l+2}\|_{L^\infty}.
	\end{align}
\end{proposition}


\subsection{Existence of the compressible Euler layers} 

\begin{proposition}\label{prop5.7}
For the Euler layer problem \eqref{2.23}-\eqref{2.4-1}, assume $v^{k-1}_p$,  $\tilde{v}_e^k|_{x=0,L}$  and the forcing terms $F^k_{i}, i=1,2,3,4$ satisfying
compatibility conditions as described in Section \ref{sec4.1}. There exist solutions $(\rho_e^k, u_e^k, v_e^k, T_e^k), \, k=1,2,\cdots, N$ to  \eqref{2.23}-\eqref{2.4-1} with the following  estimations
\begin{align}\label{5.46}
	\sum_{j=0}^{m_k} \|\nabla^j v_e^k\, \la Y\ra^{l_k}\|_{L^\infty}  &\lesssim L \sum_{j=0}^{1+ m_k} \|\nabla^j F^k \la Y\ra^{2+l_k}\|  +  L \|v^{k-1}_p(\cdot,0)\|_{W^{3+m_k,\infty}_x}\nonumber\\
	&\quad  +  \|v^{k-1}_p(\cdot, 0)\|_{W_x^{m_k,\infty}}.
\end{align}

\begin{align}\label{5.46-1}
	\sum_{j=0}^{m_k-1} \|\pa_Y^j (\rho^k_e, T^k_e, u^k_e)\, \la Y\ra^{l_k}\|_{L^\infty} 
	&\lesssim  \sum_{j=0}^{m_k-1} \|\pa_Y^j (\tilde{\rho}^k_e, \tilde{T}^k_e, \tilde{u}^k_e)\, \la Y\ra^{l_k}\|_{L^\infty} + L \sum_{j=0}^{1+ m_k} \|\nabla^j F^k \la Y\ra^{2+l_k}\|\nonumber\\
	&\quad   +  L \|v^{k-1}_p(\cdot,0)\|_{W^{3+m_k,\infty}_x},
\end{align}
and
\begin{align}\label{5.46-2}
\sum_{j=0}^{m_k-2} \|\nabla^j \pa_x(\rho^k_e, T^k_e, u^k_e)\, \la Y\ra^{l}\|_{L^\infty} 
&\lesssim  \sum_{i=0}^{m_k-2}\|\nabla^j F^{k}\, \la Y\ra^{l}\|_{L^\infty}   +  L \sum_{j=0}^{1+ m_k} \|\nabla^j F^k \la Y\ra^{2+l_k}\|\nonumber\\
&\quad   +  L \|v^{k-1}_p(\cdot,0)\|_{W^{3+m_k,\infty}_x}  +  \|v^{k-1}_p(\cdot, 0)\|_{W_x^{m_k,\infty}}.
\end{align}
\end{proposition}

\noindent{\bf Proof.} Noting Lemma \ref{lem2.2}, then applying Proposition \ref{prop5.6}, we can  obtain the existence of $v_e^i$ with 
\begin{align}\label{5.53}
	\sum_{j=0}^{m_k} \|\nabla^j v_e^k\, \la Y\ra^{l_k}\|_{L^\infty}  &\lesssim L \sum_{j=0}^{1+ m_k} \|\nabla^j F^k \la Y\ra^{2+l_k}\|  +  L \|v^{k-1}_p(\cdot,0)\|_{W^{3+m_k,\infty}_x}\nonumber\\
	&\quad  +  \|v^{k-1}_p(\cdot, 0)\|_{W_x^{m_k,\infty}}.
\end{align}
which yields \eqref{5.46}.

It follows from $\eqref{2.35}_1$ and \eqref{2.4-1} that
\begin{align}
	\rho^k_e(x,Y)
	& = \tilde{\rho}^k_e(Y) + \int_0^x \f{1}{2T^0_e- |u^0_e|^2}\Big[\rho_e^0 u_e^{0}\pa_Y v^k_e -  \big[\rho_e^0 \pa_Yu_e^0  - u_e^{0} \pa_{Y}\rho^0_e  - 2\f{\rho_e^0}{u_e^0}\pa_Y T^0_e\big] v_e^{k} \nonumber\\
	& \quad  + F_{2}^k -   u_e^{0} F_{1}^k -  \f{1}{u^0_e}\big[ F_{4}^k - T_e^0F_{1}^k \big] \Big](s,Y) ds,\label{5.54}\\
	T_e^k(x,Y) &= \tilde{T}_e^k(Y) + \int_0^x \Big(\f{T_e^0}{\rho_e^0} \pa_x \rho_e^k  -  \f{2}{u_e^0}\pa_Y T^0_e\, v_e^k + \f{1}{\rho^0_e u^0_e}\big[ F_{4}^k - T_e^0F_{1}^k \big]\Big)(s,Y) ds,\label{5.55}
\end{align}
and 
\begin{align}\label{5.56}
	u^k_e(x,Y)&= \tilde{u}^k_e(Y) + \int_0^x \f{1}{\rho_e^0(Y)}\Big(- \f{u_e^0}{\rho_e^0}\, \pa_x \rho_e^k - \pa_{Y}v_{e}^k - \f{1}{\rho_e^0}\pa_{Y}\rho_e^0 \, v_{e}^k +  \f{1}{\rho_e^0} F_{1}^k\Big)(s,Y) ds.
\end{align}
Using \eqref{5.54}-\eqref{5.56}, we can establish the existence  of $(\rho^i_k, T^k_e, u^k_e)$ satisfying the following estimates
\begin{align}\label{5.57}
\sum_{j=0}^{m_k-1} \|\pa_Y^j (\rho^k_e, T^k_e, u^k_e)\, \la Y\ra^{l}\|_{L^\infty} 
&\lesssim  \sum_{j=0}^{m_k-1} \|\pa_Y^j (\tilde{\rho}^k_e, \tilde{T}^k_e, \tilde{u}^k_e)\, \la Y\ra^{l}\|_{L^\infty} +  L \sum_{i=0}^{m_k-1}\|\pa_Y^j F^{k}\, \la Y\ra^{l}\|_{L^\infty} \nonumber\\
&\quad + L \sum_{i=0}^{m_k}  \|\pa_Y^jv_e^k  \la Y\ra^{l}\|_{L^\infty},
\end{align}
and
\begin{align}\label{5.58}
\sum_{j=0}^{m_k-2} \|\nabla^j \pa_x(\rho^k_e, T^k_e, u^k_e)\, \la Y\ra^{l}\|_{L^\infty} 
&\lesssim  \sum_{i=0}^{m_k-2}\|\nabla^j F^{k}\, \la Y\ra^{l}\|_{L^\infty}   + \sum_{i=0}^{m_k-1}  \|\pa_Y^jv_e^k  \la Y\ra^{l}\|_{L^\infty}.
\end{align}
Applying \eqref{5.53} to  \eqref{5.57}-\eqref{5.58}, one concludes \eqref{5.46-1}-\eqref{5.46-2}.
Therefore the proof of Proposition \ref{prop5.7} is completed. $\hfill\Box$

\section{Proof of Theorem \ref{thm3}} \label{sec5}
Applying Proposition \ref{prop5.7} to  \eqref{2.23}-\eqref{2.4-1} and Theorem \ref{thm2} to \eqref{3.6-1A}-\eqref{3.6-2A}, step by step, we can establish the existence and uniform estimates of solutions $(\rho_e^k, u_e^k, v_e^k, T_e^k)$ ($1\leq k\leq N$) of \eqref{2.23}-\eqref{2.4-1},  and the solutions  $(\rho_p^k, u_p^k, v_p^k, T_p^k)$ ($1\leq k\leq N-1$) of \eqref{3.6-1A}-\eqref{3.6-2A} (Here $v^k_p=\mathfrak{v}^k_p-\mathfrak{v}^k_p(x,\infty)$, then $v^k_p\to 0$ with decay rate  as $y\to\infty$).

Next, we apply Theorem \ref{thm2} to obtain the solution $(\rho_p^N,\, u^{N}_p,\, \mathfrak{v}^{N}_p,\, T^{N}_p)$  of \eqref{3.6-1A}--\eqref{3.6-2A} for $k=N$ with uniform estimates as in \eqref{1.10}.  To ensure that the vertical velocity $v_s$ satisfies $v_s|{y=0}=0$ and $\displaystyle \lim_{y\to\infty}v_s=0$, we define
\begin{align}\label{05.1}
v^N_p:&= \mathfrak{v}^{N}_p  \chi(\sqrt{\v}y).
\end{align}
Then it is direct to check that $(\rho_p^N,\, u^{N}_p,\, v^{N}_p,\, T^{N}_p)$ satisfy 
\begin{align}
	\begin{cases}
		\big(\bar{\rho}^0_p T^N_p + \bar{T}^0_p \rho^N_p\big)_y = G_{1}^N,	\\[1.5mm]
		\big(\bar{\rho}^0_p u^N_p + \rho^N_p \bar{u}^0_p \big)_x + \big(\bar{\rho}^0_p v^N_p + \rho^N_p \bar{v}^0_p \big)_y = G_{2}^N - \big(\bar{\rho}^0_p \mathfrak{v}^N_p \bar{\chi}(\sqrt{\v}y)\big)_y,\\[1.5mm]
		\bar{\rho}^0_p  (\bar{U}^0_p\cdot \bar{\nabla}) u^N_p + \bar{\rho}^0_p \big\{ u^N_p  \pa_x\bar{u}^{0}_p + \rho_p^N \big(\bar{U}^0_p\cdot \bar{\nabla}\big) \bar{u}^{0}_p  +  \pa_y\bar{u}^0_p\, v^N_p \big\}  + \big(\bar{\rho}^0_p T^N_p + \bar{T}^0_p \rho^N_p \big)_x \\
		\qquad\qquad\qquad\qquad\qquad\qquad\qquad\qquad\qquad 
		=\mu \pa_{yy} u^N_p + G_{3}^N - \bar{\chi}(\sqrt{\v}y) \,\bar{\rho}^0_p \, \pa_y\bar{u}^0_p\, \mathfrak{v}^N_p ,\\[1.5mm]
		\bar{\rho}^0_p(\bar{U}^0_p \cdot \bar{\nabla})T^N_p + \bar{\rho}^0_p\big\{u^N_p\pa_x \bar{T}^0_p + v^N_p \pa_y \bar{T}^0_p\big\}   + \rho_p^N \, (\bar{U}^0_p\cdot \bar{\nabla}) \bar{T}^{0}_p +  \big\{\bar{\rho}_p^0 T^N_p + \bar{T}^0_p \rho_p^N\big\} \mbox{\rm div} \bar{U}^0_p \\[1mm]
		\qquad\quad +  \bar{\rho}^0_p \bar{T}^0_p  \big\{\pa_x u^N_p + \pa_y v^N_p \big\} \\
		 \quad = \kappa \pa_{yy}T^N_p  + 2 \mu \pa_y u^N_p \cdot \pa_y \bar{u}^{0}_p  + G_{4}^N - \bar{\chi}(\sqrt{\v}y)\, \bar{\rho}^0_p \pa_y \bar{T}^0_p  \mathfrak{v}^N_p - \bar{\rho}^0_p \bar{T}^0_p   \pa_y \big(\bar{\chi}(\sqrt{\v}y) \mathfrak{v}^N_p\big).
	\end{cases}
\end{align}
Therefore the proof of Theorem \ref{thm3} is completed.  $\hfill\Box$

\appendix
\section{ Formal expansions}\label{App-A}
We hope to construct approximate solution of the  compressible Navier-Stokes equations \eqref{1.1} in the following expansion:
\begin{align}\label{3.4}
	\begin{split}
		\rho^{\v}&
		=\rho_e^0(Y)  + \rho_p^0(x,y) + \sum_{i=1}^{\infty} \sqrt{\v}^{i} [\rho_{e}^{i}(x,Y) + \rho_{p}^{i}(x,y)],\\
		u^{\v}&
		= u_e^0(Y) + u_p^0(x,y) + \sum_{i=1}^{\infty} \sqrt{\v}^{i} [u_{e}^{i}(x,Y) + u_{p}^{i}(x,y)],\\
		v^{\v}&
		=\sqrt{\v} [v_e^1(x,Y) + v_p^0(x,y)] + \sum_{i=1}^{\infty} \sqrt{\v}^{i+1} [v_{e}^{i+1}(x,Y) + v_{p}^{i}(x,y)],\\
		T^{\v}&
		= T_e^0(Y) + T_p^0(x,y) + \sum_{i=1}^{\infty} \sqrt{\v}^{i} \, [T_{e}^{i}(x,Y) + T_{p}^{i}(x,y)].
	\end{split}
\end{align}


\smallskip

\begin{lemma}\label{lem2.1}
\noindent{\it Part 1}: For $k\geq 1$, the linear compressible Euler layers $(\rho_e^k, u_e^k, v_e^k, T_e^k)$  satisfy
\begin{align} \label{2.4}
		\begin{cases}
			\dis \mbox{\bf div} \big(\rho_e^0 U_e^{k} + \rho_e^{k} U_e^0\big)=F_{1}^k,\\
			\dis \rho_e^0 (U_e^{0}\cdot\nabla)u_e^k + \rho_e^0 (U_e^{k}\cdot\nabla)u_e^0 +   \big(\rho_e^{0} T_e^k + \rho_e^k T_e^0\big)_x = F_{2}^k,\\
			\dis \rho_e^0 (U_e^0\cdot\nabla) v_e^{k}  + \big(\rho_e^{0} T_e^k + \rho_e^k T_e^0\big)_Y = F_{3}^k,\\
			\dis \rho_e^0 (U_e^{0}\cdot\nabla)T_e^k + \rho_e^0 U_e^k\cdot \nabla T^0_e  + \rho_e^{0} T_e^{0} \mbox{\bf div}U_e^{k} = F_{4}^k,
		\end{cases}
	\end{align}
	where $U_e^k:=(u_e^k, v_e^k)^{t}$ and
	\begin{align}
		F_{1}^k:&=- \sum_{\substack{i+j=k\\ i, j\geq1}} \mbox{\bf div}\big(\rho_e^i U_e^j\big),\\
		F_{2}^k:&= \mu \, \Delta u_e^{k-2}  + \lambda\, \pa_x \mbox{\bf div}U_e^{k-2} - \sum_{\substack{i+j+l=k\\ l, i, j\leq k-1}} \rho_e^i (U_e^j\cdot \nabla) u_e^l
		- \sum_{\substack{i+j=k\\ i, j\leq k-1}} \big(\rho_e^i T_e^j\big)_x,\\
		F_{3}^k:&= \mu\, \Delta v_e^{k-2} + \lambda\, \pa_Y \mbox{\bf div} U_e^{k-2} - \sum_{\substack{i+j+l=k,\\ i,j,l\leq k-1}} \rho_e^i U_e^j\cdot \nabla v_e^l - \sum_{\substack{i+j=k\\ i, j\leq k-1}} \big(\rho_e^i T_e^j\big)_Y,\\
		F_{4}^k:&=\kappa\, \Delta T_e^{k-2} + 2\mu   \sum_{\substack{i+j=k-2,\\i,j\geq0}}  S(U_e^i):S(U_e^j)  +  (\lambda-\mu)   \sum_{\substack{i+j=k-2,\\i,j\geq1}} \mbox{\bf div}U^i_e\cdot \mbox{\bf div}U^j_e\nonumber\\
		&\qquad -  \sum_{\substack{i+j+l=k,\\ 0\leq i,j,l\leq k-1}} \rho_e^i U_e^j\cdot \nabla T_e^l - \sum_{\substack{i+j+l=k,\\ i,j\geq0, k-1\geq l}} \rho_e^{i} T_e^{j}  \mbox{\bf div}U_e^{l}.
	\end{align}
It is direct to check that $F_i^1=0,\, i=1,2,3,4$.

\medskip

\noindent{Part 2}:  We denote
\begin{align}\label{3.5-0}
	\begin{cases}
		\bar{\rho}^0_p:= \rho_e^0(0) + \rho^0_p(x,y),\\
		\bar{u}^0_p:= u_e^0(0) + u^0_p(x,y),\\
		\bar{v}^0_p:= v_e^1(x,0) + v^0_p(x,y),\\
		\bar{T}^0_p:= T_e^0(0) + T^0_p(x,y),
	\end{cases}
\end{align}
then the leading order nonlinear compressible Prandtl layer satisfies
\begin{align}\label{3.6-0}
	\begin{cases}
		\pa_y\big(\bar{\rho}^0_p \bar{T}^0_p\big)=0,\\
		(\bar{\rho}^0_p \bar{u}^0_p)_x + (\bar{\rho}^0_p \bar{v}^0_p)_y =0,\\
		\bar{\rho}^0_p[\bar{u}^0_p \, \pa_x + \bar{v}^0_p \, \pa_y] \bar{u}^0_p + (\bar{\rho}^0_p \bar{T}^0_p)_x = \mu \pa_{yy} \bar{u}^0_p, \\
		\bar{\rho}^0_p[\bar{u}^0_p \, \pa_x + \bar{v}^0_p \, \pa_y] \bar{T}^0_p + \bar{\rho}^0_p \bar{T}^0_p [\pa_x \bar{u}^0_p + \pa_y \bar{v}^0_p] = \kappa \pa_{yy}\bar{T}^0_{p} + \mu |\pa_y \bar{u}^0_p|^2.
	\end{cases}
\end{align}

\smallskip

\noindent{Part 3}: For $k\geq 1$, the linear compressible Prandtl layers $(\rho_p^k, u_p^k, v_p^k, T_p^k)$ satisfy
\begin{align}\label{3.6-1}
	\begin{cases}
		\big(\bar{\rho}^0_p T^k_p + \bar{T}^0_p \rho^k_p\big)_y = G_{1}^k,	\\[1.5mm]
		\big(\bar{\rho}^0_p u^k_p + \rho^k_p \bar{u}^0_p \big)_x + \big(\bar{\rho}^0_p[ v^k_p + v^{k+1}_e(x,0)] + \rho^k_p \bar{v}^0_p \big)_y = G_{2}^k,\\[1.5mm]
		\rho_p^k \big(\bar{U}^0_p\cdot \bar{\nabla}\big) u^{0}_p  +   
		\bar{\rho}^0_p  (\bar{U}^0_p\cdot \bar{\nabla}) u^{k}_p + \bar{\rho}^0_p \big\{ u^k_p  \pa_xu^{0}_p + [v^k_p +  v^{k+1}_e(x,0)] \pa_yu^0_p\big\} \\
		\qquad\qquad\qquad\qquad\qquad\qquad\qquad\qquad\,\,  + \big(\bar{\rho}^0_p T^k_p + \bar{T}^0_p \rho^k_p \big)_x 
		=\mu \pa_{yy} u^{k}_p + G_{3}^k,\\[1.5mm]
		\bar{\rho}^0_p(\bar{U}^0_p \cdot \bar{\nabla})T^k_p + \bar{\rho}^0_p\big\{u^k_p\pa_xT^0_p + [v^k_p + v^{k+1}_e(x,0)]\pa_yT^0_p\big\}   + \rho_p^k \, (\bar{U}^0_p\cdot \bar{\nabla}) T^{0}_p \\[1mm]
		\qquad +  \big\{\bar{\rho}_p^0 T^k_p + \bar{T}^0_p \rho_p^k\big\} \mbox{\rm div} U^0_p  +  \bar{\rho}^0_p \bar{T}^0_p  \big\{\pa_x u^k_p + \pa_y  v^k_p \big\}
		= \kappa \pa_{yy}T^k_p  + 2 \mu \pa_y u^{k}_p \cdot \pa_y u^{0}_p  + G_{4}^k,
	\end{cases}
\end{align}
where $\mbox{\rm div} = \pa_x + \pa_y$, $\bar{\nabla}=(\pa_x,\pa_y)^t$,  $G_{4}^k=G_{41}^k + G_{42}^k$ and
{\small
	\begin{align}
		G_{1}^k:&=  (\mu +\lambda) \pa_{yy}v^{k-2}_p + \lambda \pa_{xy} u_p^{k-2} + \mu   \pa_{xx}v^{k-4}_p  -  \sum_{\substack{i+j=k\\ 0\leq i,j\leq k-1}} \big(\rho_p^i T^j_p\big)_y \nonumber\\
		&\quad- \sum_{\substack{i+j+l=k\\ i,j,l\geq0, j\leq k-1}} \Big(\f{y^{l}}{l!} \big[\pa_Y^{l}\rho^i_e(x,0)\cdot T^j_p +  \pa_Y^{l}T^i_e(x,0)\cdot \rho^j_p\big] \Big)_y \nonumber\\
		&\quad -  \sum_{n+m=k-2} \Big(\rho_p^n + \sum_{\substack{i+l=n\\ i,l\geq0}} \f{y^{l}}{l!} \, \pa_Y^{l}\rho^i_e(x,0)\Big) \bigg\{ \sum_{\substack{i+j+l=m\\ i,j,l\geq0}} (U^i_p\cdot \bar{\nabla}) \Big(\f{y^l}{l!} \, \pa^l_{Y}v_e^{j+1}(x,0)\Big) \nonumber \\
		&\qquad + \sum_{\substack{i+j=m\\ i,j\geq0}} (U^i_p\cdot \bar{\nabla}) v^{j}_p + \sum_{\substack{i+j+l=m\\i,j,l\geq0}} \f{y^{l}}{l!} \Big( \pa_Y^{l}u^i_e(x,0)\cdot \pa_x v_p^j + \pa^l_{Y}v_e^{i+1}(x,0) \cdot \pa_y v_p^j \Big)\bigg\}\nonumber \\
		&\quad - \sum_{\substack{i+j+l+n+h=k-2\\ i,j,l,n,h\geq0}} \rho_{p}^{h} \f{y^l}{l!}\big[\pa_Y^{l}u^i_e(x,0) \pa_x +  \pa^l_{Y}v_e^{i+1}(x,0) \, \pa_y\big]\Big(\f{y^n}{n!} \, \pa^n_{Y}v_e^{j+1}(x,0)\Big),
\end{align}}
{\small
	\begin{align}
		G_{2}^k:&= -\sum_{\substack{i+j=k\\ 0\leq i,j\leq k-1}} \Big\{\big(\rho_p^j u^i_p\big)_x + \big(\rho_p^j v^i_p\big)_y\Big\} - \sum_{\substack{i+j+l=k\\ i,j,l\geq0,\, j\leq k-1}} \f{y^{l}}{l!} \big[  \pa_Y^{l}\rho^i_e(x,0) \cdot  u^j_p + \pa_Y^{l} u^i_e(x,0) \cdot  \rho^j_p\big]_x  \nonumber\\
		&\quad - \sum_{\substack{i+j+l=k\\ i,j,l\geq0, \, j\leq k-1}} \f{y^{l}}{l!} \big[  \pa_Y^{l}\rho^i_e(x,0) \cdot v^j_p + \mathbf{1}_{\{i\leq k-1\}} \pa_Y^{l} v^{i+1}_e(x,0) \cdot  \rho^j_p\big]_y ,
\end{align}}
{\small
	\begin{align}
		G_{3}^k:&= (\mu + \lambda) \pa_{xx} u^{k-2}_p +  \lambda \pa_{xy} v_p^{k-2} - \big(\bar{U}^0_p\cdot \bar{\nabla}\big) u^{0}_p \cdot  \sum_{\substack{i+l=k\\ i,l\geq0}} \f{y^{l}}{l!} \, \pa_Y^{l}\rho^i_e(x,0) \nonumber\\
		&\quad - \sum_{\substack{i+j+l=k\\ i,j,l\geq0, i\leq k-1}}  \bar{\rho}^0_p  
		\cdot  (U^i_p\cdot \bar{\nabla}) \Big(\f{y^{l}}{l!} \, \pa_Y^{l}u^j_e(x,0)\Big)  - \sum_{\substack{i+j=k\\ 0\leq i,j\leq k-1}}  \bar{\rho}^0_p  (U^i_p\cdot \bar{\nabla}) u^{j}_p \\
		&\quad- \sum_{\substack{i+j+l=k\\i,j,l\geq0, j\leq k-1}}  \bar{\rho}^0_p  \f{y^{l}}{l!} \Big( \pa_Y^{l}u^i_e(x,0)\cdot \pa_x u_p^j + \mathbf{1}_{\{i\leq k-1\}}\pa^l_{Y}v_e^{i+1}(x,0) \cdot \pa_y u_p^j \Big) \nonumber\\
		&\quad - \sum_{\substack{n+m=k\\ 0\leq m,n\leq k-1}}  \bigg\{ \Big(\rho_p^n + \sum_{\substack{i+l=n\\ i,l\geq0}} \f{y^{l}}{l!} \, \pa_Y^{l}\rho^i_e(x,0)\Big) 
		\cdot \bigg[\sum_{\substack{i+j+l=m\\ i,j,l\geq0}} (U^i_p\cdot \bar{\nabla}) \Big(\f{y^{l}}{l!} \, \pa_Y^{l}u^j_e(x,0)\Big)\nonumber \\
		&\qquad +  \sum_{\substack{i+j=m\\ i,j\geq0}} (U^i_p\cdot \bar{\nabla}) u^{j}_p + \sum_{\substack{i+j+l=m\\i,j,l\geq0}} \f{y^{l}}{l!} \Big( \pa_Y^{l}u^i_e(x,0)\cdot \pa_x u_p^j + \pa^l_{Y}v_e^{i+1}(x,0) \cdot \pa_y u_p^j \Big)\bigg] \bigg\} \nonumber\\
		&\quad -  \sum_{\substack{m+i+j+l+n=k\\ h,i,j,l,n\geq0, \, m\leq k-1}}  \rho_{p}^{m} \f{y^{l}}{l!} \, \Big[ \pa_Y^{l}u^i_e(x,0)\cdot \pa_x  + \pa^l_{Y}v_e^{i+1}(x,0)\cdot \pa_y \Big] \Big(\f{y^{n}}{n!} \, \pa_Y^{n} u^j_e(x,0)\Big) \nonumber\\
		&\quad - \sum_{\substack{i+j=k\\ 0\leq i,j\leq k-1}} \big(\rho_p^i T^j_p\big)_x  - \sum_{\substack{i+j+l=k\\ i,j,l\geq0, j\leq k-1}} \f{y^{l}}{l!} \Big[\pa_Y^{l}\rho^i_e(x,0)\cdot T^j_p +  \pa_Y^{l}T^i_e(x,0)\cdot \rho^j_p\Big]_x,
\end{align}}
{\small
	\begin{align}
		G_{41}^k:&= -\sum_{\substack{m+n=k\\ 0\leq n,m\leq k-1}} \Big(\rho_p^m + \sum_{\substack{i+l=m\\ i,l\geq0}} \f{y^{l}}{l!} \, \pa_Y^{l}\rho^i_e(x,0)\Big) \cdot \bigg\{ \sum_{\substack{i+j+l=n\\ i,j,l\geq0}} (U^i_p\cdot \bar{\nabla}) \Big(\f{y^l}{l!} \, \pa^l_{Y}T_e^j(x,0)\Big) \nonumber\\
		&\qquad + \sum_{\substack{i+j=n\\ i,j\geq0}} (U^i_p\cdot \bar{\nabla}) T^{j}_p
		+  \sum_{\substack{i+j+l=n\\i,j,l\geq0}} \f{y^{l}}{l!} \Big( \pa_Y^{l}u^i_e(x,0)\cdot \pa_x T_p^j + \pa^l_{Y}v_e^{i+1}(x,0) \cdot \pa_y T_p^j \Big) \bigg\} \nonumber\\
		&\quad -  \sum_{\substack{i+j+l+n+m=k\\ i,j,l,n,m\geq0, m\leq k-1}} \rho_{p}^{m}  \f{y^{l}}{l!} \Big( \pa_Y^{l}u^i_e(x,0)\cdot \pa_x + \mathbf{1}_{\{i\leq k-1\}} \pa^l_{Y}v_e^{i+1}(x,0) \cdot \pa_y \Big) \Big(\f{y^n}{n!} \, \pa^n_{Y}T_e^j(x,0)\Big)\nonumber\\
		&\quad - \sum_{\substack{n+m=k\\ 0\leq n,m\leq k-1}} \Big\{ \sum_{\substack{i+j=n\\ i,j\geq0}} \rho_p^i T^j_p + \sum_{\substack{i+j+l=n\\ i,j,l\geq0}} \f{y^{l}}{l!} \Big[\pa_Y^{l}\rho^i_e(x,0)\cdot T^j_p +  \pa_Y^{l}T^i_e(x,0)\cdot \rho^j_p\Big]\Big\}\nonumber\\
		&\qquad\quad  \times\bigg\{\sum_{\substack{i+l=m\\ i,l\geq0}} \Big[\Big(\f{y^{l}}{l!} \, \pa_Y^{l}u^i_e(x,0) \Big)_x + \Big(\f{y^l}{l!} \, \pa^l_{Y}v_e^{i+1}(x,0)\Big)_y\Big] + \mbox{\rm div} U^m_p \bigg\} \nonumber\\
		&\quad - \Big\{\sum_{\substack{i+j=k\\ 0\leq i,j\leq k-1}} \rho_p^i T^j_p  + \sum_{\substack{i+j+l=k\\ i,j,l\geq0, j\leq k-1}} \f{y^{l}}{l!} \Big[\pa_Y^{l}\rho^i_e(x,0)\cdot T^j_p +  \pa_Y^{l}T^i_e(x,0)\cdot \rho^j_p\Big]\Big\}  \mbox{\rm div} U^0_p\nonumber\\
		&\quad - \big(\bar{\rho}^0_p\bar{T}^0_p - \rho^0_e(0) T^0_e(0)  \big) \sum_{\substack{i+l=k\\ i,l\geq0}} \Big[\Big(\f{y^{l}}{l!} \, \pa_Y^{l}u^i_e(x,0) \Big)_x + \mathbf{1}_{\{i\leq k-1\}}\Big(\f{y^l}{l!} \, \pa^l_{Y}v_e^{i+1}(x,0)\Big)_y\Big] \nonumber\\
		&\quad  -  \sum_{\substack{i+j+l+n+m=k\\ i,j,l,n\geq0, m\leq k-1}} \f{y^{l}}{l!}\f{y^n}{n!} \pa_Y^{l}\rho^i_e(x,0) \cdot \pa^n_{Y}T_e^i(x,0) \cdot\mbox{\rm div} U^m_p \nonumber\\
		&\quad - \sum_{\substack{i+l=k\\ i,l\geq0}} \f{y^{l}}{l!} \, \pa_Y^{l}\rho^i_e(x,0)\cdot 
		(\bar{U}^0_p\cdot \bar{\nabla}) T^{0}_p ,
\end{align}}
and
{\small
	\begin{align}
		G_{42}^k:&=  \kappa \pa_{xx} T_{p}^{k-2}   +  (\lambda-\mu) \sum_{\substack{i+j=k-2\\ i,j\geq0}} [ \pa_{x} u_p^i + \pa_y v_p^i] \cdot [ \pa_{x} u_p^j + \pa_y v_p^j]  + 2 \mu  \pa_x v_p^{k-2} \cdot \pa_y u^{0}_p\nonumber\\
		&\quad + 2  (\lambda-\mu)  \sum_{\substack{i+l+m=k-2\\ i,l,m\geq0}} \Big[\Big(\f{y^{l}}{l!} \, \pa_Y^{l}u^i_e(x,0) \Big)_x + \Big(\f{y^l}{l!} \, \pa^l_{Y}v_e^{i+1}(x,0)\Big)_y\Big] \cdot \mbox{div} U^m_p \nonumber\\
		&\quad + 4\mu \sum_{\substack{i+j+l=k-2\\ i,j,l\geq0}} \bigg\{ \Big(\f{y^{l}}{l!} \, \pa_Y^{l}u^i_e(x,0)\Big)_x \cdot \pa_x u_p^j + \Big(\f{y^l}{l!} \, \pa^l_{Y}v_e^{i+1}(x,0)\Big)_y \, \pa_y v_p^j \bigg\} \nonumber\\
		&\quad + 2 \mu \Big\{\sum_{\substack{i+l+m=k\\ i,l,m\geq0,m\leq k-1}} \Big(\f{y^{l}}{l!} \, \pa_Y^{l}u^i_e(x,0)\Big)_y  + \sum_{\substack{i+l+m=k-2\\ i,l\geq0}} \Big(\f{y^l}{l!} \, \pa^l_{Y}v_e^{i+1}(x,0)\Big)_x\Big\} \nonumber\\
		&\qquad \times  \big\{\pa_y u^{m}_p +  \pa_x v_p^{m-2}\big\}   + 2\mu  \sum_{\substack{i+j=k-2\\ i,j\geq0}} \big[\pa_x u_p^i\cdot  \pa_x u_p^j + \pa_y v_p^i\cdot \pa_y v_p^j\big] \nonumber\\
		&\quad   + \mu \sum_{\substack{i+j=k\\ 0\leq i,j\leq k-1}}  \big\{\pa_y u^{i}_p +  \pa_x v_p^{i-2}\big\} \cdot  \big\{\pa_y u^{j}_p +  \pa_x v_p^{j-2}\big\}.
\end{align}}

\noindent\underline{\it Part 4: Boundary conditions.} The following conditions hold:
\begin{align}\label{3.57-0}
\bar{u}^0_p\equiv u_{e}^{0}(x,0) + u_{p}^{0}(x,0)=0\quad \& \quad \bar{v}^0_p\equiv v_{e}^{1}(x,0) + v_{p}^{0}(x,0)=0,
\end{align}
\begin{align}
u_{e}^{i}(x,0) + u_{p}^{i}(x,0)&=0, \quad i=1,2,\cdots,\label{3.57}\\
v_{e}^{i+1}(x,0) + v_{p}^{i}(x,0)&=0,\quad i=1,2,\cdots,\label{3.57-1}
\end{align}  
and
\begin{align}\label{3.59-1}
		\begin{cases}
			T_{e}^{0}(x,0) + T_{p}^{0}(x,0)=\mathscr{T}_b(x),\\[1.5mm]
			T_{e}^{i}(x,0) + T_{p}^{i}(x,0)=0,\quad i=1,2,\cdots,
		\end{cases}
		\,\, \mbox{for DBC},
\end{align}
or
	\begin{align}\label{3.59-2}
		\begin{cases}
			\pa_yT^0_p(x,0)=0, \\[1.5mm]
			\pa_yT^i_p(x,0) + \pa_Y T_e^{i-1}(x,0)=0,\quad i=1,2,\cdots,
		\end{cases}
		\,\, \mbox{for NDC}.
	\end{align}
We remark that the boundary conditions \eqref{3.57-0}, $\eqref{3.59-1}_1$, and $\eqref{3.59-2}_1$ are used to solve the nonlinear Prandtl equations \eqref{3.6-0}; while \eqref{3.57}, \eqref{3.57-1}, $\eqref{3.59-1}_2$, and $\eqref{3.59-2}_2$ serve as the boundary conditions for the linear compressible Prandtl layer equations \eqref{3.6-1}. In addition, \eqref{3.57-1} is used to determine the next-order linear Euler layers \eqref{2.4}.
\end{lemma}

\noindent{\bf Proof.}  Substituting \eqref{3.4} into the  compressible Navier-Stokes equations \eqref{1.1},  we first compare the order of $\sqrt{\v}$ for interior Euler layers, then we compare the order of $\sqrt{\v}$ for the Prandtl layers (For the interaction of Euler and Prandtl layers, we need to use Taylor expansion for Euler layers so that they can be uniformly represented by $y$ variable), then the linear Euler layer equations \eqref{2.4}, the nonlinear Prandtl equations \eqref{3.6-0}, and the linear Prandtl layer equations \eqref{3.6-1} follow directly. We omit the details for simplicity. 

Noting the formal expansion \eqref{3.4}, to satisfy  \eqref{1.2}-\eqref{1.3-1}, it is clear to obtain \eqref{3.57}-\eqref{3.59-2}. Therefore the proof of Lemma \ref{lem2.1} is completed. $\hfill\Box$

\bigskip

To solve the linear compressible Euler layer equations \eqref{2.4}, we assume  the  uniform subsonic condition
\begin{align}\label{1.3-6}
	\sup_{Y\in \mathbb{R}_+}\big[2T_e^0(Y) - |u^0_e(Y)|^2\big]\geq \tilde{c}_0>0.
\end{align}
then we reformulate it in the following lemma.
\begin{lemma}[Reformulation of linear compressible Euler layers]\label{lem2.2}
For $k\geq 1$,  the system \eqref{2.4} is equivalent to the following equations
\begin{align}\label{2.23}
&-\pa_Y\Big(\f{\rho_e^0 T_e^0 u^0_e }{2T_0 - |u_e^0|^2}\, \pa_{Y}v_{e}^k\Big) -\pa_{x}\big(\f12 \rho_e^0 u_e^0 \pa_x v_e^k \big) \nonumber\\
&\quad   + \pa_Y\Big(\Big[\f{T_e^0[\rho_e^0  \pa_Yu^0_e  - u^0_e \pa_{Y}\rho_e^0  -  2\f{\rho_e^0}{u_e^0}\pa_Y T^0_e]}{2T_0 - |u_e^0|^2} + \f{\rho_e^0}{u_e^0}\pa_Y T^0_e\Big] v_e^k\Big)\nonumber\\
&=\pa_Y\Big(T^0_eF_{2}^k -  T^0_e u_e^{0} F_{1}^k -  \f{T^0_e}{u^0_e}\big[ F_{4}^k - T_e^0F_{1}^k \big]\Big) - \pa_Y\Big(\f{1}{2u^0_e}\big[ F_{4}^k - T_e^0F_{1}^k \big]\Big),
\end{align}
and
\begin{align}\label{2.35}
	\begin{cases}
	\dis	\pa_x \rho^k_e 
		 = \f{1}{2T^0_e- |u^0_e|^2}\Big\{\rho_e^0 u_e^{0}\pa_Y v^k_e -  \big[\rho_e^0 \pa_Yu_e^0  - u_e^{0} \pa_{Y}\rho^0_e  - 2\f{\rho_e^0}{u_e^0}\pa_Y T^0_e\big] v_e^{k} \\
		\dis \qquad\qquad\qquad\qquad\qquad  + F_{2}^k -   u_e^{0} F_{1}^k -  \f{1}{u^0_e}\big[ F_{4}^k - T_e^0F_{1}^k \big] \Big\},\\[1.5mm]
		\dis \pa_xT_e^k  = \f{T_e^0}{\rho_e^0} \pa_x \rho_e^k  -  \f{2}{u_e^0}\pa_Y T^0_e\, v_e^k + \f{1}{\rho^0_e u^0_e}\big[ F_{4}^k - T_e^0F_{1}^k \big],\\[1.5mm]
		\dis \rho_e^0 \pa_x u^k_e   =  - \f{u_e^0}{\rho_e^0}\, \pa_x \rho_e^k - \pa_{Y}v_{e}^k - \f{1}{\rho_e^0}\pa_{Y}\rho_e^0 \, v_{e}^k +  \f{1}{\rho_e^0} F_{1}^k.
	\end{cases}
\end{align}
Under the condition \eqref{1.3-6}, it is clear that \eqref{2.23} is uniform elliptic.
Then, for \eqref{2.23}-\eqref{2.35}, we impose the following initial and boundary conditions
\begin{align}\label{2.4-1}
	\begin{cases}
		\dis v_e^k|_{Y=0}=-v^{k-1}_p(x,0),\quad   v_e^{k}|_{x=0}=\tilde{v}_{e,0}^k(Y),\quad v_e^{k}|_{x=L}=\tilde{v}_{e,L}^k(Y),\\[1.5mm]
		\dis (\rho^k_e, u^k_e, T^k_e)|_{x=0}=(\tilde{\rho}^k_e, \tilde{u}^k_e, \tilde{T}^k_e)(Y), \quad \lim_{Y\to\infty} (\rho^k_e, u^k_e, v^k_e, T^k_e) =(0, 0, 0, 0).
	\end{cases}
\end{align}
We remark that the forcing terms $F^k_i$ {\rm(}$i=1,2,3,4${\rm)} are specific functions (see Lemma \ref{lem2.1} for details) which are determined by the $j$-th {\rm(}$j\leq k-1${\rm)} Euler layers, and should therefore be regarded as known. 
\end{lemma}

\noindent{\bf Proof.}    It follows from \eqref{2.4} that 
\begin{align}
	&\rho_e^0 \pa_x u^k_e = - \rho_e^0 \pa_Y v^k_e - \pa_{Y}\rho^0_e\, v^{k}_e - u^0_e \pa_x \rho^k_e  + F_{1}^k,\label{2.31-0}\\
	&\rho_e^0 \pa_xT_e^k = T_e^0 \pa_x \rho^k_e - 2\f{\rho_e^0}{u_e^0}\pa_Y T^0_e\, v_e^k + \f{1}{u^0_e}\big[ F_{4}^k - T_e^0F_{1}^k \big].\label{2.31-1}
\end{align}
Substituting \eqref{2.31-0}-\eqref{2.31-1} into $\eqref{2.4}_2$,   one has
\begin{align}\label{2.31-6}
	&- \rho_e^0 u_e^{0}\pa_Y v^k_e  + \big[\rho_e^0 \pa_Yu_e^0  - u_e^{0} \pa_{Y}\rho^0_e  - 2\f{\rho_e^0}{u_e^0}\pa_Y T^0_e\big] v_e^{k} + \big[2T^0_e- |u^0_e|^2\big] \pa_x \rho^k_e  \nonumber\\
	&\quad  = F_{2}^k -   u_e^{0} F_{1}^k -  \f{1}{u^0_e}\big[ F_{4}^k - T_e^0F_{1}^k \big],
\end{align}
which implies that 
\begin{align}\label{2.31-2}
	&-\pa_Y\Big(\f{\rho_e^0 T_e^0 u^0_e }{2T_0 - |u_e^0|^2}\, \pa_{Y}v_{e}^k\Big)  + \pa_Y\Big(\f{T_e^0[\rho_e^0  \pa_Yu^0_e  - u^0_e \pa_{Y}\rho_e^0  -  2\f{\rho_e^0}{u_e^0}\pa_Y T^0_e]}{2T_0 - |u_e^0|^2} v_{e}^k\Big)  + \big(T_e^0\pa_{x} \rho_e^k\big)_Y \nonumber\\
	&=\pa_Y\Big(T^0_eF_{2}^k -  T^0_e u_e^{0} F_{1}^k -  \f{T^0_e}{u^0_e}\big[ F_{4}^k - T_e^0F_{1}^k \big]\Big).
\end{align}

Applying $\pa_x$ to $\eqref{2.4}_3$, one obtains that 
\begin{align*}
	\pa_{x}\big(\rho_e^0 u_e^0 \pa_x v_e^k\big) + \big(\rho_e^{0} \pa_x T_e^k\big)_{Y} + \big(T_e^0 \pa_x \rho_e^k \big)_{Y} = \pa_{Y} F_{3}^k,
\end{align*}
which, together with \eqref{2.31-1}, yields that 
\begin{align}\label{2.31-3}
	\f12 \pa_{x}\big(\rho_e^0 u_e^0 \pa_x v_e^k\big) + \big(T_e^0 \pa_x \rho_e^k \big)_{Y} -  \big(\f{\rho_e^0}{u_e^0}\pa_Y T^0_e\, v_e^k\big)_{Y} = \pa_Y\Big(\f{1}{2u^0_e}\big[ F_{4}^k - T_e^0F_{1}^k \big]\Big).
\end{align}
Combining \eqref{2.31-2} and \eqref{2.31-3}, one obtains that
\begin{align}\label{2.31-4}
	&-\pa_Y\Big(\f{\rho_e^0 T_e^0 u^0_e }{2T_0 - |u_e^0|^2}\, \pa_{Y}v_{e}^k\Big)   - \f12 \pa_{x}\big(\rho_e^0 u_e^0 \pa_x v_e^k\big) \nonumber\\
	&\quad  + \pa_Y\Big(\Big[\f{T_e^0[\rho_e^0  \pa_Yu^0_e  - u^0_e \pa_{Y}\rho_e^0  -  2\f{\rho_e^0}{u_e^0}\pa_Y T^0_e]}{2T_0 - |u_e^0|^2} + \f{\rho_e^0}{u_e^0}\pa_Y T^0_e\Big]\, v_e^k\Big) \nonumber\\
	&=\pa_Y\Big(T^0_eF_{2}^k -  T^0_e u_e^{0} F_{1}^k -  \f{T^0_e}{u^0_e}\big[ F_{4}^k - T_e^0F_{1}^k \big]\Big) - \pa_Y\Big(\f{1}{2u^0_e}\big[ F_{4}^k - T_e^0F_{1}^k \big]\Big).
\end{align}
which concludes \eqref{2.23}.  It is easy to derive \eqref{2.35} from \eqref{2.31-0}-\eqref{2.31-6}.


The boundary condition  $v_e^k|_{Y=0}$  follows directly from \eqref{3.57-1}.
Under the uniform subsonic condition \eqref{1.3-6}, it is clear that \eqref{2.23} is uniform elliptic, so we need to provide boundary data on both sides $v^k_e|_{x=0,L}$, while we need only to give $(\rho^k_e, u^k_e, T^k_e)|_{x=0}$. So the given data in \eqref{2.4-1} are the five functions $\tilde{v}_{e,0}^k, \tilde{v}_{e,L}^k, \tilde{\rho}^k_e, \tilde{u}^k_e, \tilde{T}^k_e$.
Therefore the proof of Lemma  \ref{lem2.2} is completed. $\hfill\Box$

\medskip

\begin{lemma}[Reformulation of linear compressible Prandtl layers]
For  $k\geq 1$, we denote 
\begin{align}\label{0.23}
	\mathfrak{v}^k_p:=v^k_p + v^{k+1}(x,0),
\end{align}
then the linear compressible Prandtl layer equations \eqref{3.6-1} become equivalent to
\begin{align}\label{3.6-1A}
	\begin{cases}
		\big(\bar{\rho}^0_p T^k_p + \bar{T}^0_p \rho^k_p\big)_y = G_{1}^k,	\\[1.5mm]
		\big(\bar{\rho}^0_p u^k_p + \rho^k_p \bar{u}^0_p \big)_x + \big(\bar{\rho}^0_p \mathfrak{v}^k_p + \rho^k_p \bar{v}^0_p \big)_y = G_{2}^k,\\[1.5mm]
		\bar{\rho}^0_p  (\bar{U}^0_p\cdot \bar{\nabla}) u^{k}_p + \bar{\rho}^0_p \big\{ u^k_p  \pa_x\bar{u}^{0}_p + \rho_p^k \big(\bar{U}^0_p\cdot \bar{\nabla}\big) \bar{u}^{0}_p  +  \pa_y\bar{u}^0_p \mathfrak{v}^k_p \big\}  + \big(\bar{\rho}^0_p T^k_p + \bar{T}^0_p \rho^k_p \big)_x \\
		\qquad\qquad\qquad\qquad\qquad\qquad\qquad\qquad\qquad\qquad\qquad\qquad\,\, 
		=\mu \pa_{yy} u^{k}_p + G_{3}^k,\\[1.5mm]
		\bar{\rho}^0_p(\bar{U}^0_p \cdot \bar{\nabla})T^k_p + \bar{\rho}^0_p\big\{u^k_p\pa_x\bar{T}^0_p + \mathfrak{v}^k_p \pa_y\bar{T}^0_p\big\}   + \rho_p^k \, (\bar{U}^0_p\cdot \bar{\nabla}) \bar{T}^{0}_p +  \big\{\bar{\rho}_p^0 T^k_p + \bar{T}^0_p \rho_p^k\big\} \mbox{\rm div} \bar{U}^0_p \\[1mm]
		\qquad\qquad\qquad\qquad\qquad\qquad  +  \bar{\rho}^0_p \bar{T}^0_p  \big\{\pa_x u^k_p + \pa_y \mathfrak{v}^k_p \big\}  = \kappa \pa_{yy}T^k_p  + 2 \mu \pa_y u^{k}_p \cdot \pa_y \bar{u}^{0}_p  + G_{4}^k.
	\end{cases}
\end{align}
For \eqref{3.6-1A}, the initial and boundary conditions are 
\begin{align}\label{3.6-2A}
	\begin{cases}
		\dis (u^k_{p}, T^k_{p})|_{x=0}=\big(\tilde{u}^k_p, \tilde{T}^k_{p}\big)(y), \quad \lim_{y\to \infty} (u^k_p,  T^k_p)=(0, 0), \\[1mm]
		\dis (u^k_{p}, \mathfrak{v}^k_{p})|_{y=0} = -\big(u^k_e, 0\big)(x,0),\\[1mm]
		\dis \pa_yT^k_p |_{y=0} =  - \pa_YT^{k-1}_e(x,0) \quad \mbox{for NBC},\\[1mm]
		T^k_{p} |_{y=0} = - T^k_e(x,0)\qquad \mbox{for DBC}.
	\end{cases}
\end{align}
\end{lemma}

\noindent{\bf Proof.} In view of \eqref{0.23}, \eqref{3.6-1}, \eqref{3.57}, \eqref{3.57-1}, $\eqref{3.59-1}_2$, and $\eqref{3.59-2}_2$, the lemma follows directly. $\hfill\Box$

\begin{remark}
It is noted that \eqref{3.6-1A} is a parabolic system for $(u^k_p, T^k_p)$ with $x$ as the evolution variable.  The given data in \eqref{3.6-2A} are the two functions $\tilde{u}^k_p, \tilde{T}^k_{p}$.  The forcing terms $G^k_{i}, i=1,2,3,4$ are specifically functions (see Lemma \ref{lem2.1} for details) which are determined by $i$-th ($i\leq k$) Euler layers and $i$-th ($i\leq k-1$) linear compressible Prandtl layers, and should be regarded as known for the purposes of stating the result. 
Then it is clear that $(\rho^k_p, u^k_p, \mathfrak{v}^k_p, T^k_p)$ satisfy the linear Prandtl layer equations \eqref{3.6-10}. Since we require $\dis\lim_{y\to \infty}v^k_p=0$, we deduce $\mathfrak{v}^k_p(x,\infty)=v^{k+1}_e(x,0)$, that means once \eqref{3.6-1A} is solved, we can then determine $v^k_p=\mathfrak{v}^k_p - v^{k+1}_e(x,0)$. Hence we shall apply Theorem \ref{thm2} to construct  $(\rho^k_p, u^k_p, v^k_p, T^k_p)$ for $1\leq k\leq N-1$ in \eqref{3.4-0}.
 
For $k = N$, we first apply Theorem \ref{thm2} to obtain $(\rho_p^N, u_p^N, \mathfrak{v}_p^N, T_p^N)$. However, unlike the cases with $1 \leq k \leq N-1$, we cannot define $v^N_p$ in the same manner. To ensure the boundary conditions $v_s|_{y=0} = 0$ and $\displaystyle \lim_{y \to \infty} v_s = 0$ are satisfied, we define
\begin{align}
	v^N_p = \chi(\sqrt{\v} y), \mathfrak{v}_p^N.
\end{align}
Then $(\rho_p^N, u_p^N, v_p^N, T_p^N)$ will satisfy a system similar to \eqref{3.6-1A}, but with additional higher-order corrections of order $\sqrt{\v}$.
\end{remark}

\smallskip


\section{The Approximate Solution: $(\rho_s, u_s, v_s, T_s)$}\label{App-B}
\noindent{\it 1. The equations for $(\rho_s, u_s, v_s, T_s)$ with high-order $\v$ remainders:}
Noting Theorem \ref{thm3},  it is direct to obtain the equations of $(\rho_s, u_s, v_s, T_s)$:
{\small 
\begin{align}\label{B.1}
	\begin{cases}
	\dis	\pa_x(\rho_s u_s) + \pa_Y(\rho_s v_s) =  R_{s0},\\
	\dis	\rho_s  \big[u_s \pa_x u_s + v_s \pa_Y u_s \big]  + \pa_x p_s  = \mu \, \v {\bf\Delta} u_s  + \lambda\, \v\, \pa_x\mbox{\bf div} U_s + R_{s1},\\
	\dis	\rho_s  \big[u_s \pa_x v_s + v_s \pa_Y v_s \big]  + \pa_Y p_s  = \mu \, \v {\bf\Delta} v_s  + \lambda\, \v\, \pa_Y\mbox{\bf div} U_s + \sqrt{\v} R_{s2},\\
	\dis	\rho_s  \big[u_s \pa_x T_s + v_s \pa_Y T_s \big] + p_s \mbox{\bf div} U_s =  \kappa \v {\bf\Delta}T_s + 2\mu \v |S(U_s)|^2 + (\lambda-\mu) \v  |\mbox{\bf div}U_s|^2 + R_{s3}
	\end{cases}
\end{align}}
where we have  denoted $U_s:=(u_s, v_s)$. The remainders $R_{si}, i=0,1,2,3$ are determined by the Euler layers and Prandtl layers constructed in Theorem \ref{thm3}.  The exact expression  of $R_{s0}$ is:
{\small
	\begin{align}\label{B.2}
		R_{s0}
		:&= \sum_{k\geq N+1} \sqrt{\v}^{k} \pa_x\Big(\sum_{\substack{i+j=k\\ 0\leq i,j\leq N}} \rho_p^j u^i_p + \sum_{\substack{i+j+l=k\\ 0\leq i,j,l\leq N}} \f{y^{l}}{l!} \big[  \pa_Y^{l}\rho^i_e(x,0)u^j_p + \pa_Y^{l} u^i_e(x,0) \rho^j_p\big]\Big) \nonumber\\
		&\quad + \sum_{k\geq N+1} \sqrt{\v}^{k} \pa_y \bigg(\sum_{\substack{i+j=k\\ 0\leq i,j\leq N}} \rho_p^j v^i_p + \sum_{\substack{i+j+l=k\\ 0\leq i,j,l\leq N}} \f{y^{l}}{l!} \pa_Y^{l}\rho^i_e(x,0) v^j_p + \sum_{\substack{i+j+l=k\\ 0\leq  j,l\leq N\\ 0\leq i\leq N-1}} \f{y^{l}}{l!} \pa_Y^{l} v^{i+1}_e(x,0) \rho^j_p \bigg)\nonumber\\
		&\quad + \sum_{\substack{i+j\geq N+1\\1\leq i,j\leq N}} \sqrt{\v}^{i+j} \mbox{\bf div}(\rho_{e}^i  U_{e}^j) +  \sqrt{\v}^{N+1}  \pa_y\Big(y^{N+1} \Big\{\sum_{i,j=0}^{N} \sqrt{\v}^{i+j} O(1) v^j_p + \sum_{i=0}^{N-1} \sum_{j=0}^{N} \sqrt{\v}^{i+j} O(1)  \rho^j_p \Big\}\Big) \nonumber\\
		&\quad  + \sqrt{\v}^{N+1}  \sum_{i,j=0}^{N}  \sqrt{\v}^{i+j} \pa_x \Big(y^{N+1}\big[O(1)  u^j_p + O(1) \rho^j_p\big]\Big)  - \sqrt{\v}^{N+1} \big(\bar{\rho}^0_p \mathfrak{v}^N_p \bar{\chi}(\sqrt{\v}y)\big)_y.
\end{align}}
where $O(1)$ in \eqref{B.2} is determined by the remainder term of Taylor expansions in $Y$ at $Y=0$ of $(\rho^i_e, u^i_e, v^i_e, T^i_e)$, i.e.,
\begin{align*}
O(1)=\pa_Y^{N+1}\rho_e^i(x, \tilde{Y}) , \quad \pa_Y^{N+1}u_e^i(x, \tilde{Y}), \quad  \pa_Y^{N+1}v_e^i(x, \tilde{Y}), \quad \pa_Y^{N+1}T_e^i(x, \tilde{Y}),
\end{align*}
for some $\tilde{Y}\in [0,Y]$.
We remark that the remainders $R_{si}, i=1,2,3$ have similar but more complicate expressions, we omit the details here for simplicity of presentation. A direct calculation shows that 
\begin{align}\label{B.3}
\sum_{i+j\leq \f12K_0} \|\pa^i_x\pa^j_y(R_{s0}, R_{s1}, R_{s2}, R_{s2}) \la y\ra^{\f{\mathfrak{b}}4}\|_{L^\infty} \lesssim \sqrt{\v}^{N}. 
\end{align}

For the mass equation $\eqref{B.1}_1$, noting \eqref{B.2}, we can write $R_{s0}$ in the following conservation form 
\begin{align}
	R_{s0}&=\pa_x A_{s1} + \pa_y A_{s2},
\end{align}
where
{\small
\begin{align}
A_{s1}:&=\sum_{k\geq N+1} \sqrt{\v}^{k} \Big(\sum_{\substack{i+j=k\\ 0\leq i,j\leq N}} \rho_p^j u^i_p + \sum_{\substack{i+j+l=k\\ 0\leq i,j,l\leq N}} \f{y^{l}}{l!} \big[  \pa_Y^{l}\rho^i_e(x,0)u^j_p + \pa_Y^{l} u^i_e(x,0) \rho^j_p\big]\Big) \nonumber\\
&\quad + \sum_{\substack{i+j\geq N+1\\1\leq i,j\leq N}} \sqrt{\v}^{i+j}  \rho_{e}^i  U_{e}^j + \sqrt{\v}^{N+1} \sum_{i,j=0}^{N}  \sqrt{\v}^{i+j} y^{N+1}\big[O(1)  u^j_p + O(1) \rho^j_p\big],\\
A_{s2}:&=\sum_{k\geq N+1} \sqrt{\v}^{k} \bigg\{\sum_{\substack{i+j=k\\ 0\leq i,j\leq N}} \rho_p^j v^i_p + \sum_{\substack{i+j+l=k\\ 0\leq i,j,l\leq N}} \f{y^{l}}{l!} \pa_Y^{l}\rho^i_e(x,0) v^j_p + \sum_{\substack{i+j+l=k\\ 0\leq  j,l\leq N\\ 0\leq i\leq N-1}} \f{y^{l}}{l!} \pa_Y^{l} v^{i+1}_e(x,0) \rho^j_p \bigg\}\nonumber\\
&\quad + \sum_{\substack{i+j\geq N+1\\1\leq i,j\leq N}} \sqrt{\v}^{i+j-1} \rho_{e}^i v_{e}^j  +  \sqrt{\v}^{N+1}  y^{N+1} \Big\{\sum_{i,j=0}^{N} \sqrt{\v}^{i+j} O(1) v^j_p + \sum_{i=0}^{N-1} \sum_{j=0}^{N} \sqrt{\v}^{i+j} O(1)  \rho^j_p \Big\}\nonumber\\
&\quad - \sqrt{\v}^{N} \bar{\rho}^0_p \mathfrak{v}^N_p \bar{\chi}(\sqrt{\v}y).
\end{align}}

\medskip

\noindent{\it 2. Boundary conditions:} Noting \eqref{3.6-01}, \eqref{2.4-1}, \eqref{3.6-2A} \eqref{0.23},  and \eqref{05.1}, the following boundary conditions hold:
\begin{align}
\begin{cases}
(u_s, v_s)|_{Y=0}=0,\\
\pa_YT_s|_{Y=0}= \sqrt{\v}^{N}\pa_YT^{N}_e(x,0)\,\, \mbox{for NBC},\\
T_s|_{Y=0}=\mathscr{T}_b(x)\,\, \mbox{for DBC}.
\end{cases}
\end{align}

\smallskip

\noindent{\it 3. Adjustment:} Noting that in general $R_{s0}(x,0) \neq 0$ and $\pa_YT_e^{N}(x,0)\neq0$, we know that
\begin{align*}
	v_{sy}|{y=0} = O(\sqrt{\v}^{N}) \neq 0 \quad \mbox{and}\quad \pa_YT{s}|_{y=0}=O(\sqrt{\v}^N)\neq 0,
\end{align*}
so it is necessary to make some adjustments for later use in \cite{Guo-Wang}.
We define 
\begin{align}
\begin{cases}
	\dis \tilde{\rho}_s:=\rho_s,\qquad \dis\tilde{u}_s:= u_s  \\
\dis\tilde{v}_s:=v_s - \sqrt{\v}\f{1}{\rho_s} y\chi(y)\big(\pa_xA_{s1} + \pa_yA_{s2}\big)(x,0),\\
\dis\tilde{T}_s:= T_s - \sqrt{\v}^N\chi(y) T_e^N\, {\bf 1}_{\{NBC\}},
\end{cases}
\end{align}
which, together with \eqref{B.1}, yields that 
{\small 
\begin{align}\label{B.10}
\begin{cases}
\dis	\pa_x(\tilde{\rho}_s \tilde{u}_s) + \pa_Y(\tilde{\rho}_s \tilde{v}_s) = \tilde{R}_{s0},\\
\dis	\tilde{\rho}_s  \big[\tilde{u}_s \pa_x \tilde{u}_s + \tilde{v}_s \pa_Y \tilde{u}_s \big]  + \pa_x \tilde{p}_s  = \mu \, \v {\bf\Delta} \tilde{u}_s  + \lambda\, \v\, \pa_x\mbox{\bf div} \tilde{U}_s +  \tilde{R}_{s1},\\
\dis	\tilde{\rho}_s  \big[\tilde{u}_s \pa_x \tilde{v}_s + \tilde{v}_s \pa_Y \tilde{v}_s \big]  + \pa_Y \tilde{p}_s  = \mu \, \v {\bf\Delta} \tilde{v}_s  + \lambda\, \v\, \pa_Y\mbox{\bf div} \tilde{U}_s + \sqrt{\v}\tilde{R}_{s2},\\
\dis	\tilde{\rho}_s  \big[\tilde{u}_s \pa_x \tilde{T}_s + v_s \pa_Y \tilde{T}_s \big] + \tilde{p}_s \mbox{\bf div} \tilde{U}_s =  \kappa \v {\bf\Delta}\tilde{T}_s + 2\mu \v |S(\tilde{U}_s)|^2 + (\lambda-\mu) \v  |\mbox{\bf div}\tilde{U}_s|^2 + \tilde{R}_{s3}
\end{cases}
\end{align}}
where $p_s=\tilde{\rho}_s \tilde{T}_s$, and 
\begin{align}\label{B.11}
	\sum_{i+j\leq \f12K_0} \|\pa^i_x\pa^j_y(\tilde{R}_{s0}, \tilde{R}_{s1}, \tilde{R}_{s2}, \tilde{R}_{s2}) \la y\ra^{\f{\mathfrak{b}}4}\|_{L^\infty} \lesssim \sqrt{\v}^{N}. 
\end{align}
The exact expressions  of $(\tilde{R}_{s1}, \tilde{R}_{s2}, \tilde{R}_{s2})$ can be obtained easily, we ignore this for simplicity.

For later use in \cite{Guo-Wang}, we note that 
\begin{align}
\tilde{R}_{s0}=\pa_x\tilde{A}_{s1}(x,y) + \pa_y \tilde{A}_{s2}(x,y),
\end{align}
with
\begin{align}\label{B.12}
\begin{split}
\tilde{A}_{s1}:&=A_{s1} - \pa_y\big(y\chi(y)\big)  A_{s1}(x,0) +\chi'(y)\int_0^x A_{s2}(t,0) dt,\\
\tilde{A}_{s2}:&=A_{s2} - \chi(y) A_{s2}(x,0) - y \chi(y) (\pa_yA_{s2})(x,0).
\end{split}
\end{align}

\smallskip

With above construction,  it is clear that 
\begin{align}\label{B.13}
\tilde{A}_{s1}(x,y) \cong \sqrt{\v}^{N} y,\quad \tilde{A}_{s2}(x,y) \cong \sqrt{\v}^{N} y^2\quad \mbox{for}\,\, y\in [0,1],
\end{align}
which, together with  $\eqref{B.10}_1$, yields that 
\begin{align}
\pa_y\tilde{v}_s|_{y=0}=0\quad \Longrightarrow\quad \tilde{v}_s\cong y^2\quad \mbox{for}\,\, y\in [0,1].
\end{align}
Moreover, it holds
\begin{align}
\pa_Y\tilde{T}_s|_{Y=0}= 0\quad \mbox{for NBC}.
\end{align}

\bigskip	

\noindent{\bf Acknowledgments.} 
Yong Wang’s research  is partially supported by  National Key R\&D Program of China grant No. 2021YFA1000800, the National Natural Science Foundation of China grants No. 12421001 \& 12288201, and CAS Project for Young Scientists in Basic Research, grant No. YSBR-031. 

\medskip

\noindent{\bf Data availability:} The manuscript has no associated data.

\medskip

\noindent{\bf Conflict of Interest:} None of the authors is in a situation which gives rise to a conflict of interest.

\end{document}